\documentclass[lecture]{pcms-l}

\makeatletter

\newif\iffirstlecture\firstlecturefalse

\newcommand{\lectureseries}{\firstlecturetrue
              \secdef\@lectureseries\@slectureseries} 

\newcommand{\@lectureseries}[2][default]{\chapter*{#2}%
              \gdef\thelectureseries{#1}} 

\newcommand{\@slectureseries}[1]{\chapter*{#1}}

\renewcommand{\auth}{\secdef\@auth\@sauth}

\newcommand{\@auth}[2][default]{\vspace{-1pc}{\raggedleft
        \Large\bf\noindent
        #2\endgraf
        \vspace*{2pc}
        }
        \def\@author{#1}%
}

\newcommand{\@sauth}[1]{\vspace{-1pc}{\raggedleft
        \Large\bf\noindent
        #1\endgraf
        \vspace*{2pc}
        }
        \def\@author{#1}%
}

\def\lecture#1{\global\Lecturetrue\global\Monographfalse
\iffirstlecture\else\chapter*{}\fi\firstlecturefalse
  \global\let\sectionmark\@gobble 
  \addtocounter{lecture}1\relax
  \refstepcounter{chapter}%
\gdef\thelecturename{#1\unskip}
  {\Large\bfseries
   \raggedleft
   \@xp\uppercase\@xp{\thelecturelabel} {\LARGE\thelecturenum}\\
   \vspace*{3pt}%
   \thelecturename
   \endgraf}%
  \let\@secnumber=\thelecturenum
  \@xp\lecturemark\@xp{\thelecturename}%
  \addcontentsline{toc}{chapter}{%
    \thelecturelabel\ \thelecturenum.\ \thelecturename}%
  \vspace{34\p@}\noindent}
  
\def\lecture{\global\Lecturetrue\global\Monographfalse
\iffirstlecture\else\chapter*{}\fi%
  \global\let\sectionmark\@gobble 
\secdef\@lecture\@slecture}

\def\@lecture[#1]#2{%
  \addtocounter{lecture}1\relax
  \refstepcounter{chapter}%
\gdef\thelecturename{#1\unskip}\firstlecturefalse
  {\Large\bfseries
   \raggedleft
   \@xp\uppercase\@xp{\thelecturelabel} {\LARGE\thelecturenum}\\
   \vspace*{3pt}%
    #2\unskip
   \endgraf}%
  \let\@secnumber=\thelecturenum
  \@xp\lecturemark\@xp{\thelecturename}%
  \addcontentsline{toc}{chapter}{%
    \thelecturelabel\ \thelecturenum.\ #2}%
  \vspace{34\p@}\noindent}
  
\def\slecturerunhead#1#2#3{%
    \let\@tempa\chaptername
    \uppercasenonmath{\@tempa}%
    \def\@tempb{#3\unskip}%
    \uppercasenonmath{\@tempb}%
    {\normalfont\@tempb}
    }
\def\slecturemark{
    \@secmark\markright\slecturerunhead\chaptername}%

\def\@slecture#1{%
\iffirstlecture
\gdef\thelecturename{#1\unskip}\firstlecturefalse
  {\Large\bfseries
\noindent\thelecturename
   \endgraf}%
  \let\@secnumber=\thelecturenum
  \@xp\slecturemark\@xp{\thelecturename}%
  \addcontentsline{toc}{chapter}{%
    \thelecturename}%
 \vspace{-6\p@}\noindent
\else
\gdef\thelecturename{#1\unskip}\firstlecturefalse
  {\Large\bfseries
   \raggedleft
   \@xp\uppercase\@xp{\thelecturename}
   \endgraf}%
  \let\@secnumber=\thelecturenum
  \@xp\slecturemark\@xp{\thelecturename}%
  \addcontentsline{toc}{chapter}{%
    \thelecturename}%
  \vspace{34\p@}\noindent
\fi}

\ifLecture
  \def\chapterrunhead#1#2#3{%
    \let\@tempa\@author
    \uppercasenonmath{\@tempa}%
    \uppercasenonmath{\thelectureseries}%
    \textmd{\@tempa, \thelectureseries}
    }
  \def\lecturerunhead#1#2#3{%
    \let\@tempa\chaptername
    \uppercasenonmath{\@tempa}%
    \def\@tempb{#3\unskip}%
    \uppercasenonmath{\@tempb}%
    \textmd{\@tempa\ #2. \@tempb}
    }
\else
  \let\chapterrunhead\partrunhead
\fi

\newif\ifBibliographyIsASection\BibliographyIsASectionfalse

  \def\bibliomark{
    \@secmark\markright\bibliorunhead\chaptername}%

  \def\bibliorunhead#1#2#3{%
    \let\@tempa\chaptername
    \uppercasenonmath{\@tempa}%
    \def\@tempb{#3\unskip}%
    \uppercasenonmath{\@tempb}%
    \textmd{\@tempb}
    }

\def\thebibliography#1{%
  \ifBibliographyIsASection
    \section*\refname
    \if@backmatter
      \markboth{\refname}{\refname}%
    \fi
  \else
\chapter*{}
  {\Large\bfseries
   \raggedleft
   \@xp\uppercase\@xp{\bibname} \\
   \endgraf}%
  \let\@secnumber=\thelecturenum
  \@xp\bibliomark\@xp{\bibname}%
  \addcontentsline{toc}{chapter}{%
    \bibname}%
  \vspace{34\p@}\noindent
  \fi
  \normalsize\labelsep .5em\relax
  \list{\@arabic\c@enumi.}{\settowidth\labelwidth{\@biblabel{#1}}%
  \leftmargin\labelwidth
  \advance\leftmargin\labelsep
        \usecounter{enumi}}\sloppy
  \clubpenalty9999 \widowpenalty\clubpenalty  \sfcode`\.\@m}

  \def\indexmark{
    \@secmark\markright\indexrunhead\chaptername}%

  \def\indexrunhead#1#2#3{%
    \let\@tempa\chaptername
    \uppercasenonmath{\@tempa}%
    \def\@tempb{#3\unskip}%
    \uppercasenonmath{\@tempb}%
    \textmd{\@tempb}
    }

\ifLecture
\def\theindex{\cleardoublepage
\@restonecoltrue\if@twocolumn\@restonecolfalse\fi
\columnseprule \z@ \columnsep 35pt
\def\indexchap{\@startsection
                {chapter}{1}{\z@}{8pc}{34pt}%
                {\raggedleft
                \Large\bfseries}}%
 \twocolumn[\indexchap[{\indexname}]{\@xp\uppercase\@xp{\indexname}}]
  \@xp\indexmark\@xp{\indexname}%
        \thispagestyle{plain}\let\item\@idxitem\parindent\z@
         \footnotesize\parskip\z@ plus .3pt\relax\let\item\@idxitem}
\fi

\def\@makefntext{\noindent\@makefnmark}

\def\setaddress{%
  {\let\@makefnmark\relax  \let\@thefnmark\relax
        \nobreak
        \addressnum@=\z@
        \loop\ifnum\addressnum@<\addresscount@\advance\addressnum@\@ne
           \footnote{$^{\hbox{\tiny\number\addressnum@}}$%
           \csname @address\number\addressnum@\endcsname
           \csname @curraddr\number\addressnum@\endcsname
           \csname @email\number\addressnum@\endcsname}\repeat
  \ifx\@empty\@date\else \@footnotetext{\@setdate}\fi
  \ifx\@empty\@subjclass\else \@footnotetext{\@setsubjclass}\fi
  \ifx\@empty\@keywords\else \@footnotetext{\@setkeywords}\fi
  \ifx\@empty\thankses\else \@footnotetext{%
    \def\par{\let\par\@par}\@setthanks}\fi
    }%
  \@setcopyright
}

\def\@tmpevenhead{\relax}

\def\cleardoublepage{\clearpage\if@twoside \ifodd\c@page\else
    \let\@tmpevenhead\@evenhead \let\@evenhead\relax\hbox{}\eject 
    \let\@evenhead\@tmpevenhead\if@twocolumn\hbox{}\newpage\fi\fi\fi}

\def\@setcopyright{%
  \let\copyrightyear\currentyear             
  \insert\copyins{\hsize\textwidth
    \parfillskip\z@ \leftskip\z@\@plus.9\textwidth
    \fontsize{6}{7\p@}\normalfont\upshape
    \everypar{}%
    \vskip-\skip\copyins \nointerlineskip
    \noindent\vrule\@width\z@\@height\skip\copyins
    \copyright\copyrightyear\ GMZ\par
    \kern\z@}%
}

\renewcommand{\@auth}[2][default]{{\raggedleft
        \begingroup
  \fontsize{\@xivpt}{18}\bfseries
  #2\par \endgroup
        \vspace*{2pc}
        }
        \def\@author{#1}%
}

\renewcommand{\@sauth}[1]{{\raggedleft
        \begingroup
  \fontsize{\@xivpt}{18}\bfseries
  #1\par \endgroup
        \vspace*{2pc}
        }
        \def\@author{#1}%
}

\def\@lecture[#1]#2{%
  \addtocounter{lecture}1\relax
  \refstepcounter{chapter}%
\gdef\thelecturename{#1\unskip}\firstlecturefalse
  {\Large\bfseries
   \raggedleft
   \@xp\uppercase\@xp{\thelecturelabel} {\LARGE\thelecturenum}\\
   \vspace*{3pt}%
    #2\unskip
   \endgraf}%
  \let\@secnumber=\thelecturenum
  \@xp\lecturemark\@xp{\thelecturename}%
  \addcontentsline{toc}{chapter}{%
    \thelecturelabel\ \thelecturenum.\ #2}%
  \vspace{10\p@}\noindent}
  
\def\@slecture#1{%
\iffirstlecture
\gdef\thelecturename{#1\unskip}\firstlecturefalse
  {\Large\bfseries
\noindent\thelecturename
   \endgraf}%
  \let\@secnumber=\thelecturenum
  \@xp\slecturemark\@xp{\thelecturename}%
  \addcontentsline{toc}{chapter}{%
    \thelecturename}%
 \vspace{-6\p@}\noindent
\else
\gdef\thelecturename{#1\unskip}\firstlecturefalse
  {\Large\bfseries
   \raggedleft
   \@xp\uppercase\@xp{\thelecturename}
   \endgraf}%
  \let\@secnumber=\thelecturenum
  \@xp\slecturemark\@xp{\thelecturename}%
  \addcontentsline{toc}{chapter}{%
    \thelecturename}%
  \vspace{10\p@}\noindent
\fi}

\makeatother
\usepackage{amssymb,amsmath,amsthm,bbm,paralist,url,umlaut,graphicx,color}
 \usepackage{srcltx,epsfig,psfrag}
\graphicspath{{EPS/}}
\def\currentvolume{14}
\def\currentyear{2004}
\numberwithin{section}{chapter}
\numberwithin{equation}{chapter}
\numberwithin{figure}{chapter}

\theoremstyle{plain}
\newtheorem{lemma}{Lemma}[chapter]
\newtheorem{theorem}[lemma]{Theorem}
\newtheorem{corollary}[lemma]{Corollary}
\newtheorem{proposition}[lemma]{Proposition}
\newtheorem{conjecture}[lemma]{Conjecture}
\theoremstyle{definition}
\newtheorem{definition}[lemma]{Definition}
\newtheorem{example}[lemma]{Example}
\newtheorem*{example*}{Example}
\newtheorem*{LBT}{``Lower Bound Problem''}
\newtheorem*{UBT}{``Upper Bound Problem''}
\theoremstyle{remark}

\newcommand\Z{{\mathbb Z}}
\newcommand\R{{\mathbbm R}}

\newcommand\eps{\varepsilon}

\newcommand\bb{b}
\newcommand\Izero{I_0}

\DeclareMathOperator{\DVT}{DVT}
\DeclareMathOperator{\stack}{stack}
\DeclareMathOperator{\Stack}{Stack}
\DeclareMathOperator{\proj}{proj}
\DeclareMathOperator{\cone}{cone}
\DeclareMathOperator{\conv}{conv}
\DeclareMathOperator{\aff}{aff}
\newcommand\inter{\textrm{int}}
\newcommand\bdy{\textrm{bdy}}
\newcommand\exercise[1]{Exercise~\thechapter.\ref{#1}}
\newcommand\kdiamond{\,\raisebox{-2.5pt}[0pt][0pt]{\includegraphics[height=9pt]{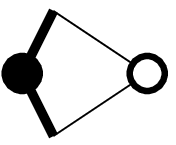}}\,}
\begin{document}

%
%



\mainmatter
\setcounter{page}{1}

\LogoOn

\lectureseries[Convex Polytopes]
{Convex Polytopes:\\{E}xtremal Constructions\\and {\itshape f}\/-Vector Shapes}

\auth{Günter M. Ziegler} 

\address{Institute of Mathematics, MA 6--2, TU Berlin,
D-10623 Berlin, Germany}
\email{ziegler@math.tu-berlin.de}


\thanks{Partially supported by the Deutsche
  Forschungs-Gemeinschaft (DFG), via the Research Center \textsc{Matheon}
  ``Mathematics in the Key Technologies'', 
  the Research Groups ``Algorithms, Structure, Randomness'' 
  and ``Polyhedral Surfaces'', and a Leibniz grant.}

\setaddress

\mbox{~}\\[20mm]

\lecture*{Introduction}

These lecture notes treat some current aspects of two closely interrelated
topics from the theory of convex polytopes: the
shapes of $f$-vectors, and extremal constructions.

The study of \emph{$f$-vectors} has had huge
successes in the last forty years.
The most fundamental one is undoubtedly the
``$g$-theorem,'' conjectured by McMullen in 1971
and proved by Billera \& Lee and Stanley in 1980,
which characterizes the $f$-vectors of simplicial and
of simple polytopes combinatorially.
See also Section~\ref{forman:sect.5.2} of Forman's article in this volume,
where $h$-vectors are discussed in connection with the
Charney--Davis conjecture.
Nevertheless, on some fundamental problems 
embarassingly little progress was made;
one notable such problem concerns the shapes of
$f$-vectors of $4$-polytopes.

A number of striking and fascinating polytope
\emph{constructions} has been proposed and
analyzed over the years. In particular,
the Billera--Lee construction produces
``all possible $f$-vectors'' of simplicial polytopes.
Less visible progress was made  outside the range
of simple or simplicial polytopes --- where our
measure of progress is that new 
polytopes ``with interesting $f$-vectors'' should be produced.
Thus, still
``it seems that overall, we are short of examples.
The methods for coming up with useful examples in mathematics
(or counterexamples to commonly believed conjectures)
are even less clear than the methods for proving mathematical
statements'' (Gil Kalai, 2000). 

These lecture notes are meant to display a 
fruitful interplay of these two areas of study:
The discussion of $f$-vector shapes suggests
the notion of ``extremal'' polytopes, that is,
of polytopes with ``extremal $f$-vector shapes.''
Our choice of constructions to be discussed here is guided by this:
We will be looking at constructions 
that produce interesting $f$-vector shapes.

After treating $3$-polytopes in the first lecture
and the $f$-vector shapes of very high-dimensional polytopes in
the second one, we will start to analyze the
case of $4$-dimensional polytopes in detail.
Thus the third lecture will explain a surprisingly simple
construction for $2$-simple $2$-simplicial $4$-polytopes,
which have symmetric $f$-vectors. 
Lecture four sketches the geometry of the cone
of $f$-vectors for $4$-polytopes, and thus identifies
the existence/\allowbreak construction of $4$-polytopes
of high ``fatness'' as a key problem.
In this direction,
the last lecture presents a very recent construction
of ``projected products of polygons,'' whose
fatness reaches $9-\eps$.
This shows that, on the topic of $f$-vectors of $4$-polytopes,
there is a narrowing gap between ``the constraints we know'' and
``the examples we can construct.''

\subsection*{Sources and Acknowledgements}

The main sources, and the basis for the
presentation in these lecture notes, are as follows.
The proof of Steinitz' theorem in Lecture~1 is due to
Alexander Bobenko and Boris Springborn \cite{BobenkoSpringborn}.
A detailed report about the work by Ludwig Danzer,
Anders Björner, Carl Lee, Jürgen Eckhoff and many others
(Lecture~2) appears in \cite[Sect.~8.6]{Z35}; see
also \cite{Bjo1} and \cite{Bjo6}.
The ``deep vertex truncation'' construction presented
in Lecture~3 is from joint work with Andreas Paffenholz \cite{Z89};
my understanding of $2$-simple $2$-simplicial polytopes benefits
also from my previous work with David Eppstein and
Greg Kuperberg \cite{Z80}.
Lecture~4 draws heavily on my 2002 Beijing ICM report \cite{Z82},
which was relying on previous studies by Marge Bayer \cite{Bay},
and with Andrea Höppner \cite{Z59}.
Finally, the construction presented in Lecture 5 was announced
in~\cite{Z97}; the intuition for it was built in previous joint
work with Nina Amenta \cite{Z51a} and Michael Joswig \cite{Z62}.

Nikolaus Witte and Thilo Schröder have forcefully directed the
problem sessions  for my Utah
lectures, and suggested a number of exercises.
Nikolaus Witte has prepared many, and the nicest, figures 
for these notes. The construction of many of the
examples and the beautiful Schlegel diagram graphics
are based on the \texttt{polymake} system
by Ewgenij Gawilow and Michael Joswig \cite{GawrilowJoswig}
\cite{GawrilowJoswig2} \cite{GawrilowJoswig4},
which everyone is invited and recommended to try out, and use.
An introduction to \texttt{polymake}, by Nikolaus Witte and Thilo
Schröder, appears as an appendix, 
pp.~\pageref{appendix:POLYMAKE}--\pageref{appendix:POLYMAKE2}.

I am grateful to all these colleagues for their work, for their
explanations and critical comments, and for support on these
lectures as well as in general.
I have benefitted a lot from the lively discussions at the
PCMI after my lectures, and from the many interesting questions and
diverse feedback I got. Boris Springborn, Nikolaus Witte,
Günter Rote, and many others provided very helpful comments on
the draft version of these lecture notes. Thank you all!

More than usually, for the trip to Utah I have depended on
the support, care, and love of Torsten Heldmann.
Without him, I wouldn't have been able to go.

\newpage

\lecture{Constructing 3-Dimensional Polytopes}

All the \emph{polytopes} considered in these lecture notes are
convex.
A \emph{$d$-polytope} is a $d$-dimensional polytope;
thus the $3$-dimensional polytopes to be discussed
in this lecture are plainly \emph{$3$-polytopes}.\footnote{%
We assume that the readers
are familiar with the basic terminology and discrete geometric 
concepts; see e.g.~\cite[Lect.~0]{Z35} or~\cite{Z49-2}.}

How many $3$-dimensional polytopes ``do we know''?
When pressed for examples, we will perhaps start with the
platonic solids: the regular tetrahedron, cube and octahedron,
icosahedron and dodecahedron. 

\begin{figure}[ht]
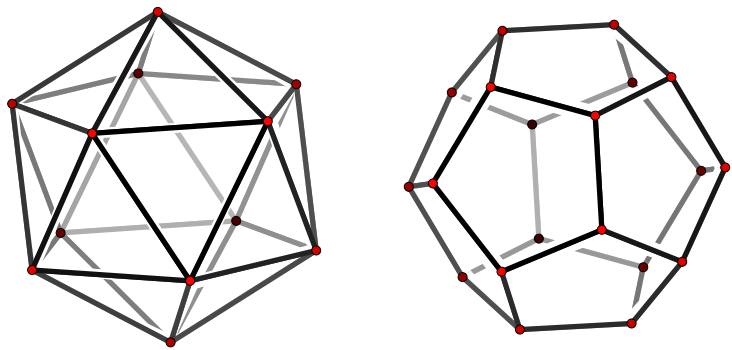

\begin{center}
\includegraphics[width=50mm]{EPS/icosahedron.eps}\quad
\includegraphics[width=50mm]{EPS/dodecahedron.eps}
\end{center}
\caption{The regular icosahedron and dodecahedron}
\label{figure:icos+dodec}
\end{figure}

\begin{figure}[ht]
\begin{center}
\includegraphics[width=66mm]{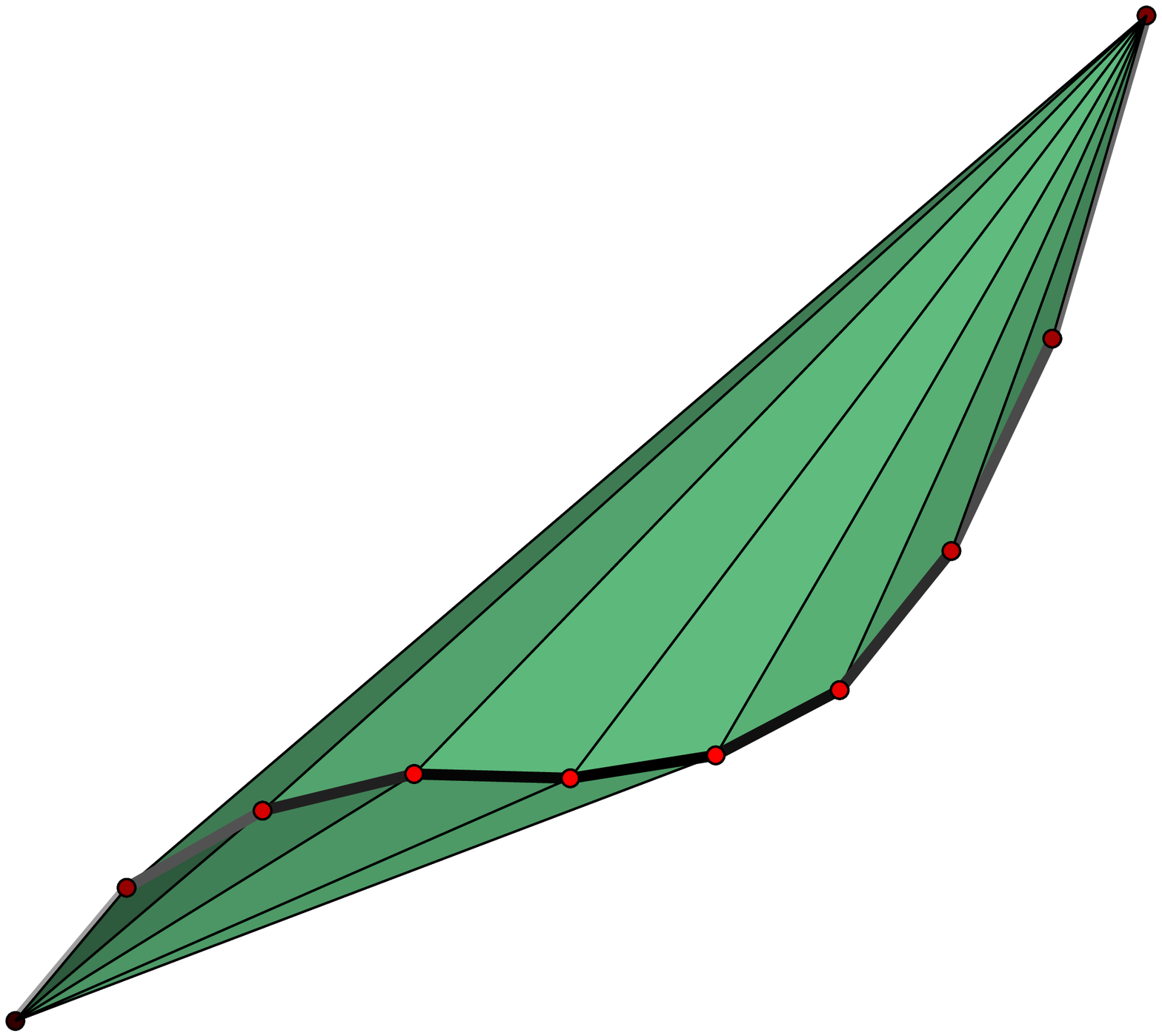}\hspace{-8mm}
\includegraphics[width=66mm]{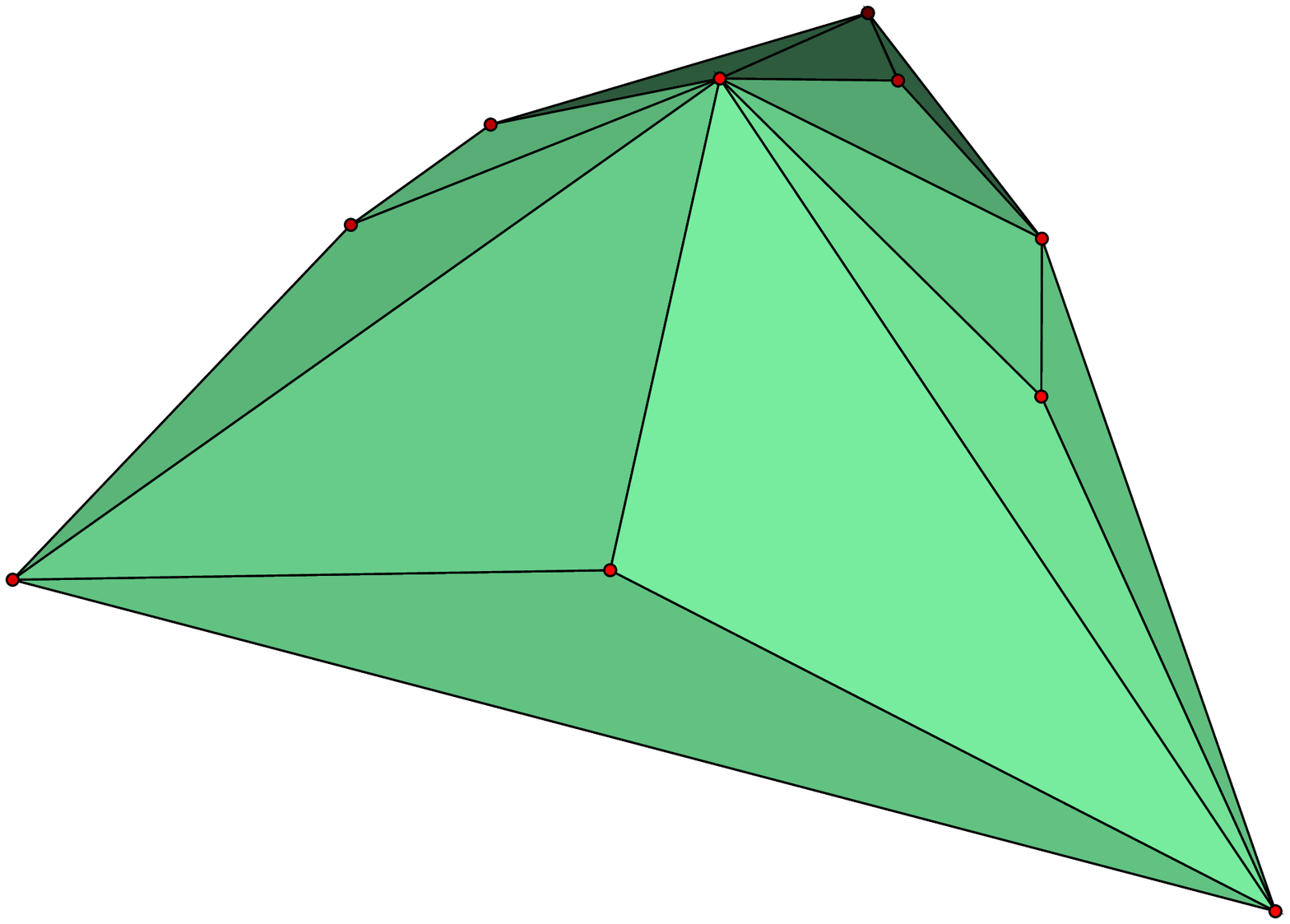}
\end{center}
\caption{A cyclic $3$-polytope $C_3(10)$
and a stacked $3$-polytope, with $10$ vertices each}
\label{figure:cyclic+stacked}
\end{figure}

The classes of \emph{stacked} and \emph{cyclic}
polytopes are of great importance for high-dimensional
polytope theory because of their extremal $f$-vectors (according to 
the \emph{lower bound theorem} and the \emph{upper bound theorem}):
Stacked polytopes arise from a simplex by repeatedly stacking pyramids onto
the facets (cf.\ Lecture~2); 
cyclic polytopes are constructed as the convex hull of
$n>d$ points on a curve of order $d$.
However, neither of these constructions produces particularly impressive
objects in dimension~$3$ 
(compare Figure~\ref{figure:cyclic+stacked}, and  \exercise{ex:CyclicStacked3d}). 

\begin{figure}[ht]
\begin{center}
\includegraphics[width=60mm]{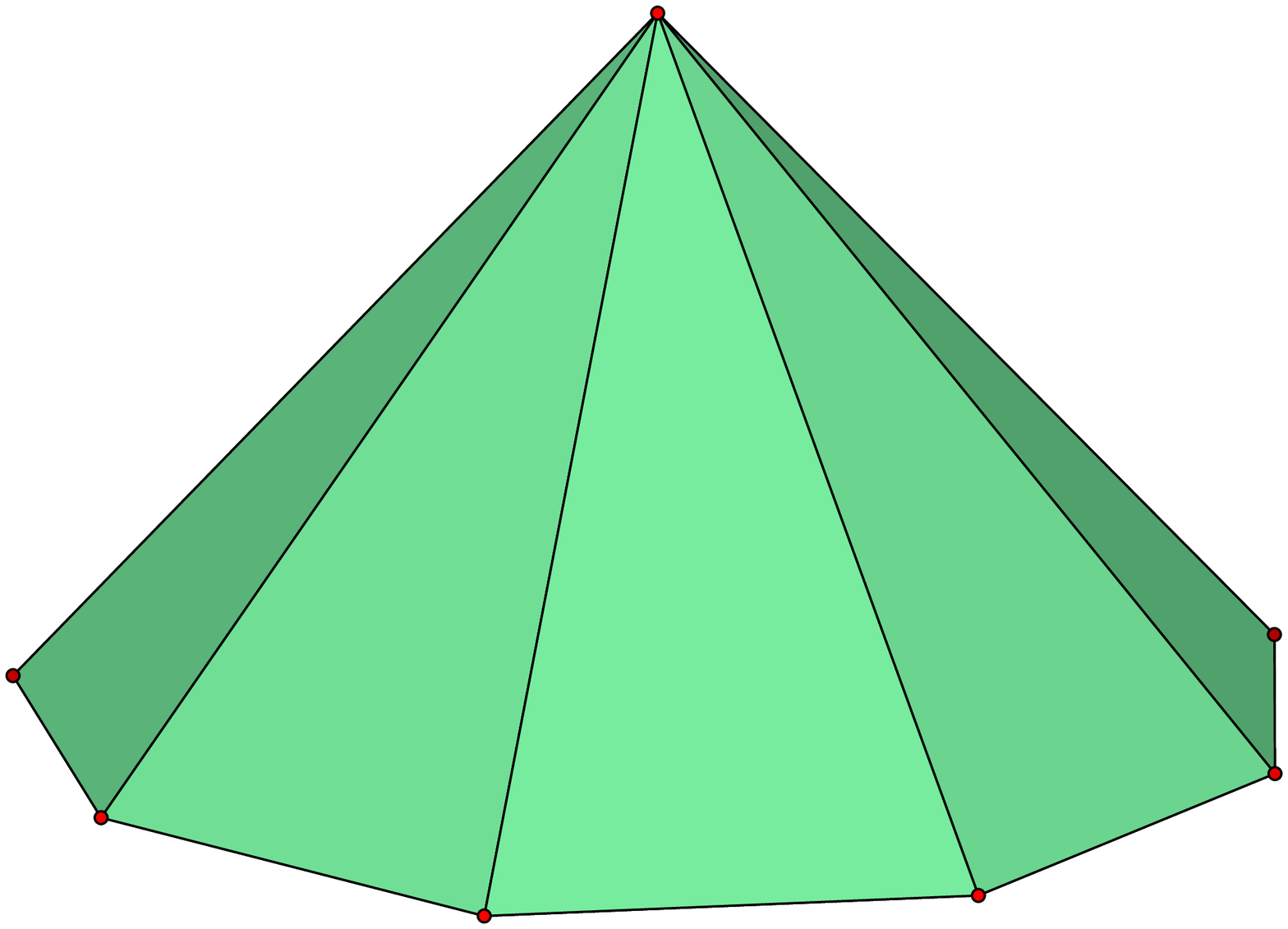}
\includegraphics[width=60mm]{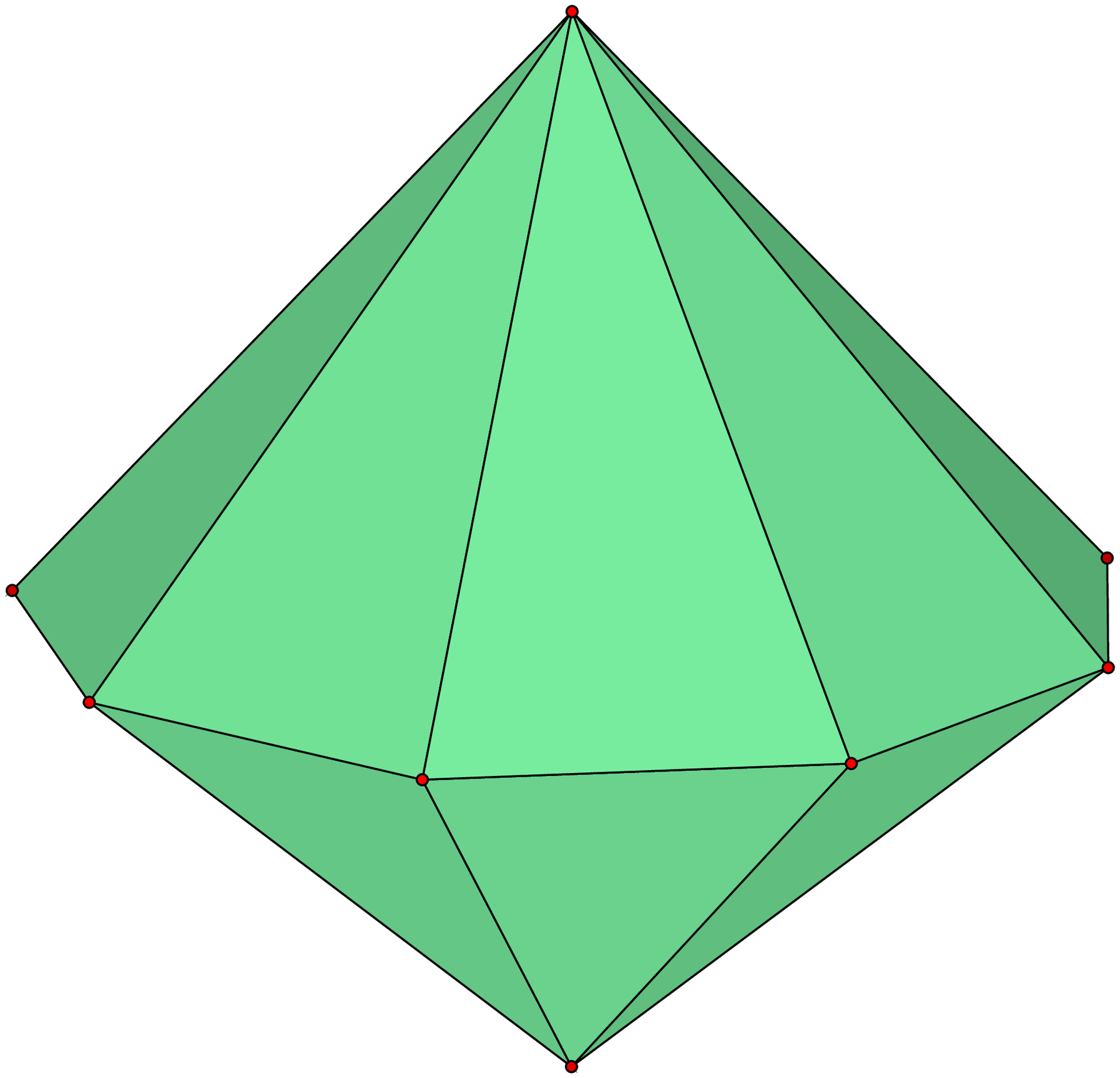}
\end{center}
\caption{The pyramid and the bipyramid over a regular $10$-gon}
\label{figure:py+bip}
\end{figure}

The same must be said about pyramids and bipyramids over $n$-gons
$(n\ge3$) --- see Figure~\ref{figure:py+bip}. 

How do we get a ``random'' $3$-polytope with lots of 
vertices? An obvious thing to look at is the convex
hull of $n$ random points on a $2$-sphere.  

\begin{figure}[ht]
\begin{center}
\includegraphics[width=75mm]{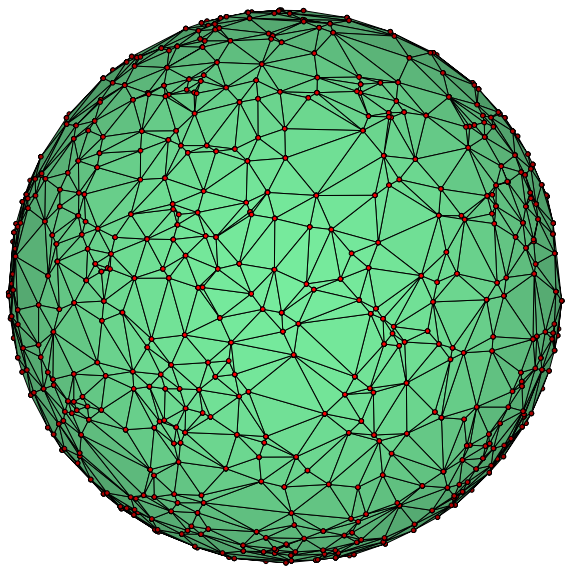}
\end{center}
\caption{A random $3$-polytope, with 1000 vertices on a sphere}
\label{figure:random3dpoly}
\end{figure}

Why is this not satisfactory? First, it produces only simplicial
polytopes (with probability~$1$), and secondly it does not even
produce all possible combinatorial types
of simplicial $3$-polytopes --- see \cite[Sect.~13.5]{Gr1-2}.
It is a quite non-trivial problem to randomly produce all
combinatorial types of polytopes of specified size (say, 
with a given number of edges). 
With the Steinitz theorem discussed below this reduces
to a search for a random planar 3-connected graph with a given
number of edges, say. See Schaeffer \cite{Schaeffer99}
for a recent treatment of this problem.

\section{The cone of \boldmath$f$-vectors}

\begin{figure}[t]
  \begin{center}
\input{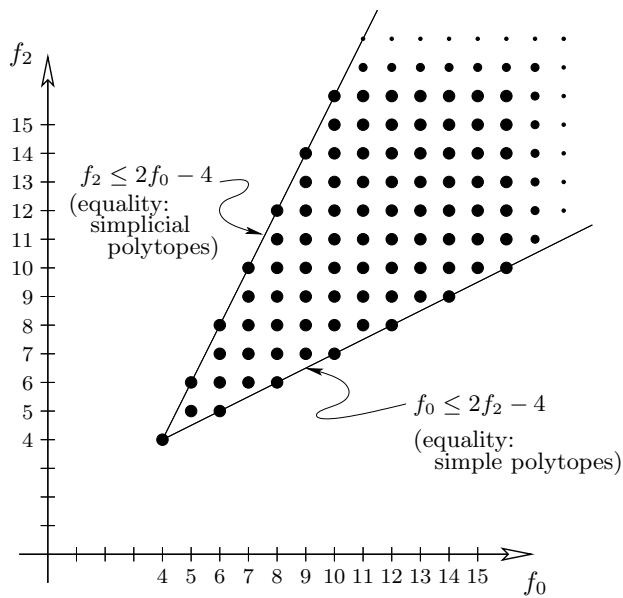}
  \end{center}
\caption{The set $\bar{\mathcal{F}}_3$, 
according to Steinitz' Lemma~\ref{lemma:Steinitz06}}
\label{figure:Steinitz06}
\end{figure}

The \emph{$f$-vector} of a $3$-polytope $P$ is the triplet of integers
\[
f(P)\ =\ (f_0,f_1,f_2)\ \in\ \Z^3,
\]
where $f_0$ is the number of vertices, $f_1$ is the number of
edges, and $f_2$ denotes the number of facets ($2$-dimensional faces).
In view of Euler's equation $f_0-f_1+f_2=2$\break
(which we take for granted here; but see
Federico \cite{Federico:Descartes}, Eppstein \cite{Eppstein:Euler}, 
and \cite[Chap.~11]{Z58-3}), the set of all $f$-vectors of $3$-polytopes,
\[
\mathcal{F}_3\ := \{(f_0,f_1,f_2)\ \in\ \Z^3: f(P)=(f_0,f_1,f_2)
\textrm{ is the $f$-vector of a $3$-polytope~$P$}\}
\]
is a $2$-dimensional set. Thus $\mathcal{F}_3$ is 
faithfully represented by the $(f_0,f_2)$-pairs of $3$-polytopes,
\[
\bar{\mathcal{F}}_3\ := \{(f_0,f_2)\ \in\ \Z^3: f(P)=(f_0,f_1,f_2)
\textrm{ for some $3$-polytope~$P$}\},
\]
as shown in Figure~\ref{figure:Steinitz06}:
The missing $f_1$-component is given by $f_1=f_0+f_2-2$.

The set of all $f$-vectors of $3$-polytopes was
completely characterized by a young Privatdozent at 
the Technische Hochschule Berlin-\allowbreak Charlottenburg (now TU Berlin),
Ernst Steinitz, in 1906: In a simple two-and-a-half-page paper
he obtained the following result, whose proof we leave to you
(\exercise{ex:Steinitzlemma}).

\begin{lemma}[Steinitz' lemma \cite{Stei3}]\label{lemma:Steinitz06}
The set of all $f$-vectors of $3$-polytopes
is given by
\[
\mathcal{F}_3\ := \{(f_0,f_1,f_2)\ \in\ \Z^3: 
f_0-f_1+f_2=2,\ f_2\le2f_0-4,\ f_0\le2f_2-4\}.
\]
\end{lemma}

This answer to the $f$-vector problem for $3$-polytopes is
remarkably simple: $\mathcal{F}_3$ is the set of \emph{all}
integral points in a $2$-dimensional convex polyhedral cone.
The three constraints that define the cone have clear
interpretations:
They are the Euler equation $f_0-f_1+f_2=2$, the upper bound
inequality $f_2\le2f_0-4$,\break which is tight exactly for the $f$-vectors
of simplicial polytopes, and its dual, $f_0\le2f_2-4$, which in the
case of equality characterizes the $f$-vectors of simple
$3$-polytopes.

For the centennial of Steinitz' lemma, in 2006, let's strive for a
characterization of the \emph{cone} spanned by the $f$-vectors of
$4$-dimensional polytopes, $\cone(\mathcal{F}_4)$. As we will see at
the beginning of Lecture~4, this is a much more modest goal than a
characterization of $\mathcal{F}_4$, which is not the set of all
integral points in a convex set: It has ``concavities'' and even
``holes.''

Steinitz' lemma, as graphed in Figure~\ref{figure:Steinitz06},
also shows that all ($f$-vectors of) convex $3$-polytopes
lie between the extremes of simple and of simplicial polytopes.
And indeed, there seems to be the misconception that 
an analogous statement should be true in higher dimensions as
well --- it isn't. As we will see, there are additional interesting
extreme cases in dimension $4$, which are by far not as
well understood as the simple and simplicial cases.

For any $3$-polytope that is not a simplex, we may compute
the ``slope''
\[
\phi(P)\ :=\ \frac{f_2-4}{f_0-4}
\]
it generates in the graph of Figure~\ref{figure:Steinitz06},
with respect to the apex $(4,4)$ of the cone, which corresponds to
a simplex.
This slope satisfies
\[
\tfrac12\ \le\ \phi(P)\ \le\ 2,
\]
where the lower bound characterizes simple polytopes,
while the upper bound is tight for simplicial polytopes.
Another interpretation of the parameter $\phi$ is that it
is a homogeneous coordinate for the cone, where the denominator
$f_0-4$ measures the ``size'' of the $f$-vector.
($\phi$ is homogeneous, so it yields $\frac00$ for
the $f$-vector of a simplex, which is the apex of the cone.
Compare \exercise{ex:fslope}.)

\section{The Steinitz theorem}\label{SteinitzTheorem}

While Steinitz' lemma from 1906 is a very simple result,
his theorem from 1922, characterizing the graphs of 
$3$-polytopes, is substantial and deep.
He knew that: He called it the ``Fundamentalsatz der konvexen Typen,''
the fundamental theorem of convex types.
Here is an informal version of it.

\begin{theorem}[Steinitz' theorem \cite{Stei1,StRa}]\label{theorem:Steinitz}
There is a bijection
\[
\{\textrm{$3$-connected planar graphs}\}
\ \ \longleftrightarrow\ \ 
\{\textrm{combinatorial types of $3$-polytopes}\}.
\]
\end{theorem}

\begin{figure}[ht]
\[
\includegraphics{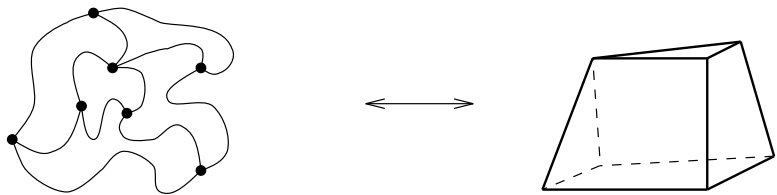}
\]
\caption{Graphs $\longleftrightarrow$ polytopes, according to
Steinitz' Theorem~\ref{theorem:Steinitz}}
\label{figure:SteinitzTheorem}
\end{figure}

\noindent
The direction ``$\longleftarrow$'' of Steinitz' theorem 
is not hard to establish.

Indeed, we do get a graph for any $3$-polytope, namely
the abstract graph whose nodes are
the vertices of the polytope, and whose arcs are
given by the edges of the polytope.
This graph is indeed planar:
To see this,
one may first produce a radial projection of the polytope
boundary (and thus of the vertices and edges) onto a sphere
that contains the polytope, and then apply a stereographic
projection \cite[\S36]{HilbertCohnVossen} to the plane. 
Or one may directly generate the
``Schlegel diagram'' and thus a straight-edge drawing of
the graph in the plane. (In Lecture~3 we will see more of this tool,
which shows its true power in the visualization of $4$-polytopes.)

To see that the graph of any $3$-polytope is $3$-connected
is also easy, using Menger's characterization of a $d$-connected
graph as a graph that cannot be disconnected by removing or blocking
less than $d$ of its vertices.
A powerful extension of this result is Balinski's theorem
\cite{Bali} \cite[Thm.~3.14]{Z35}, that the graph of any $d$-polytope
is $d$-connected.
\medskip

\noindent
Thus the hard and interesting part of Steinitz' theorem is 
the direction ``$\longrightarrow$.'' It poses a non-trivial
\emph{construction problem}: To produce a convex $3$-polytope with a
prescribed graph (a geometric object)
from an abstract planar graph (that is, from purely combinatorial data).

The first (easy) step for this is to convince oneself that
the graph characterizes the complete combinatorial structure of the
polytope. This follows from the simple observation (due to
Whitney) that the faces of the polytope correspond exactly to the
\emph{non-separating induced cycles} in the graph.

Thus we have to construct convex $3$-polytopes with
prescribed combinatorics (face lattice), as
given by a $3$-connected planar graph.
The importance of this step may be seen from the fact that
three completely different types of proofs (and construction methods!)
have been designed for it:
Let's call them Steinitz type proofs,
Tutte--Maxwell type proofs,
and Koebe--Thurston type proofs.

\subsubsection*{Steinitz type proofs} 
Such proofs (of which Steinitz gave details on one in \cite{Stei1},
and three are given in the Steinitz--Rademacher book \cite{StRa} 
that appeared after Steinitz' death),
are based on the following principle.
Any planar $3$-connected graph can be ``reduced'' to the complete
graph $K_4$ by local operations, which yields a sequence
\[
G=G_0\ \rightarrow\ 
G_1\ \rightarrow\ 
G_2\ \rightarrow\ \ \dots\ \ \rightarrow\ 
G_{N-1}\ \rightarrow\ G_N=K_4.
\]
of $3$-connected planar graphs.

This reduction sequence should then be reversed:
Starting with a simplex $\Delta_3$ (with graph $K_4$) we build up
a sequence of polytopes, 
\[
P=P_0\ \leftarrow\ 
P_1\ \leftarrow\ 
P_2\ \leftarrow\ \ \dots\ \ \leftarrow\ 
P_{N-1}\ \leftarrow\ P_N=\Delta_3,
\]
where $P_i$ is a $3$-polytope with graph~$G_i$,
again by simple/local construction steps.

Such a proof is presented in detail in~\cite[Lect.~4]{Z35},
so there is no need to do this here.
We just mention that a number of interesting extensions
and corollaries may be derived from Steinitz type proofs.
Indeed, Barnette \& Grünbaum~\cite{BaG2} proved that 
in the construction of the polytope $P$, the shape
of one face of the polytope may be prescribed.
For example, some hexagon face may be required to be a regular 
hexagon, which imposes a non-trivial additional constraint.
Similarly, Barnette~\cite{Bar4} proved with a Steinitz type
argument that a ``shadow boundary'' may be prescribed:
$P$ may be constructed in such a way that 
from some view-point outside the polytope, the edges
that bound the visible part of the surface of the polytope 
correspond to a prescribed simple cycle in the graph of the polytope
(which need not be induced).
Equivalently, we may construct $P\subset\R^3$
so that the image $\pi(P)$ of~$P$ under the orthogonal projection
$\pi:\R^3\rightarrow\R^2$ is a polygon whose edges are given
exactly by the edges of~$P$ that realize the prescribed cycle.
Indeed, the edges must be ``strictly preserved'' by the projection,
in the terminology that we will develop and use in Lecture~5.

\subsubsection*{Tutte--Maxwell type proofs}
The Tutte--Maxwell approach to realizing $3$-polytopes
works in two stages: First one gets a ``correct''
drawing of the graph in the plane, then this drawing is
lifted to $3$-space.

For the first stage, one may assume that the graph contains
a triangle face (if not, one dualizes; see \exercise{ex:triangle}).
Then the vertices of this triangle are fixed in the plane,
the edges are interpreted as ideal rubber bands, and the
other vertices are placed according to the unique and
easy-to-compute energy minimum, for which the sum of all
squared edge lengths is minimal.
This produces a correct, planar drawing of the graph without
intersections --- this is the (non-trivial) claim of 
Tutte's (1963) ``rubber band method'' \cite{Tutte:Howto};
moreover, any such drawing can be lifted to three-space
according to Maxwell--Cremona theory, which may be
traced back to work by Maxwell \cite{Maxw}
nearly one hundred years earlier (1864).
We refer to Richter-Gebert \cite[Sect.~13.1]{Rich4} for a modern treatment,
with all the proofs.

The Tutte--Maxwell proofs also buy us non-trivial corollaries: Indeed,
each combinatorial type of $3$-polytope can be realized with rational
coordinates, and thus even with integral vertex coordinates (by
clearing denominators). One can derive from a Tutte--Maxwell proof
that singly-exponential vertex coordinates suffice for this: After a
number of improvements on the original estimates by Onn \&
Sturmfels \cite{OnSt} we now know that each type of an $n$ vertex
$3$-polytope with a triangle face can be represented with vertex
coordinates in $\{0,1,2,\dots,\lfloor28.45^n\rfloor\}$ (see \cite{stein00:_realis_polyt}, 
\cite{rote:_quant} and \cite{Ribo-diss}). 
It is not clear whether polynomial-size vertex coordinates
can be achieved.

\begin{figure}[ht]
\begin{center}
\includegraphics[width=90mm]{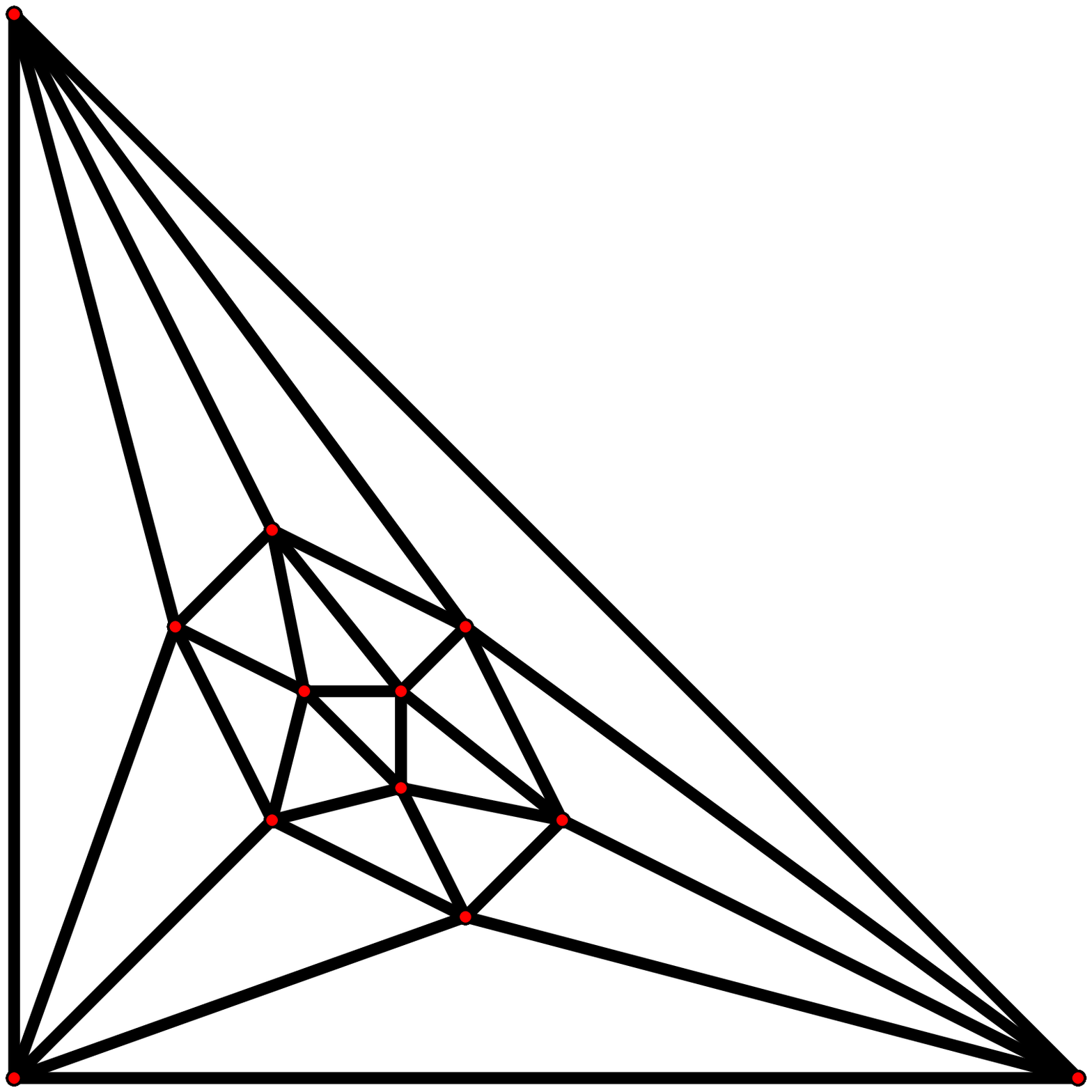}
\includegraphics[width=120mm]{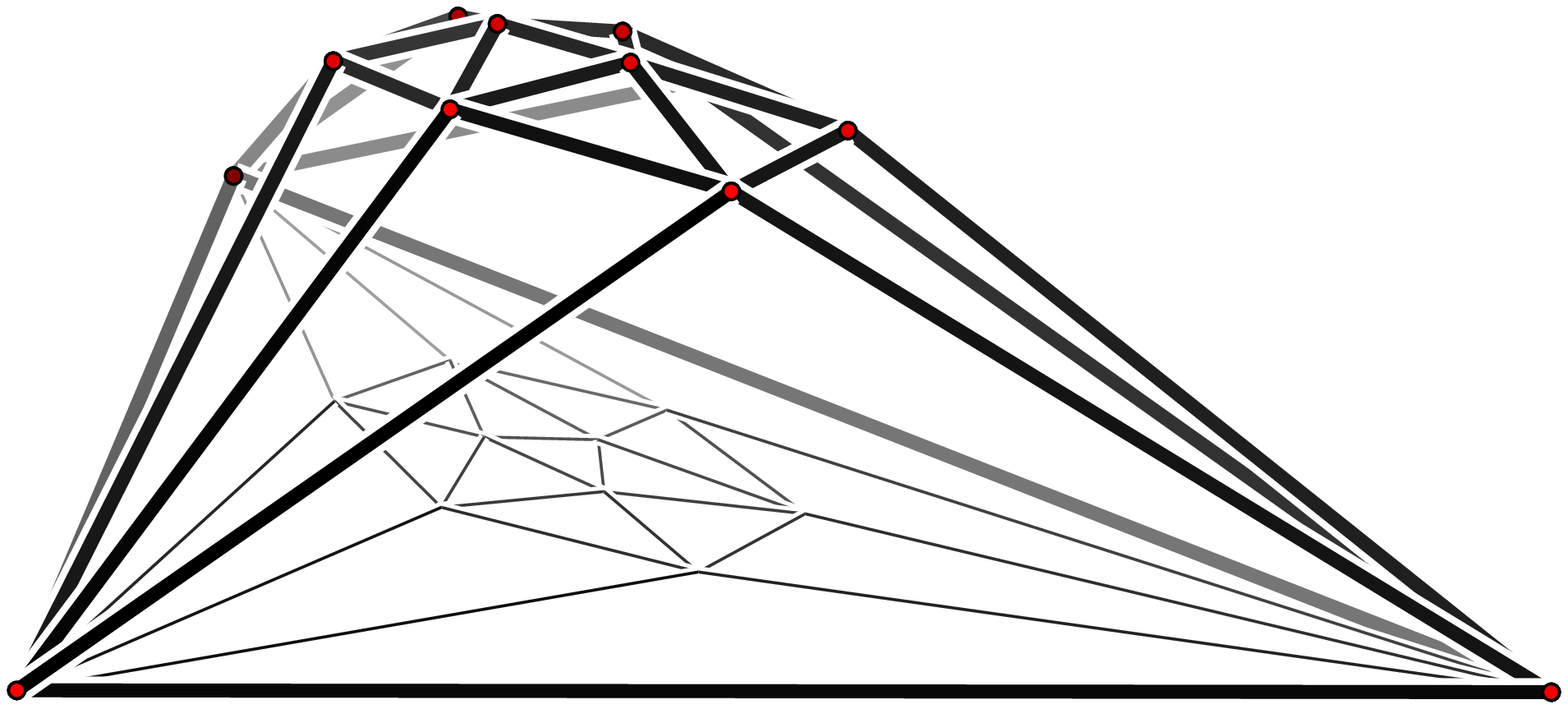}
\end{center}
\caption{A Tutte drawing of the icosahedron graph, and
the corresponding Maxwell--Cremona lifting}
\label{figure:Tutte}
\end{figure}

\subsubsection*{Koebe--Thurston type proofs}
Geometric realizations of $3$-polytopes
\emph{with all edges tangent to the sphere} may be derived
from planar circle packings. Moreover, such a representation
is essentially unique.

This seems to be essentially due to Bill Thurston \cite{Thur} --- who
traces it back to Paul Koebe's \cite{Koeb} work on complex functions,
and to work by E. M. Andreev \cite{Andr} from the sixties on
hyperbolic polyhedra.
Thurston's insight was followed up, explained, extended and 
generalized by a number of authors. Pach \&
Agarwal \cite[Chap.~8]{PachAgarwal} describe the ``standard'' proof,
based on a (non-constructive) fixed point argument.
However, 
Mohar \cite{Mohar} described an effective construction algorithm, and
Colin de Verdi\`ere \cite{CdV}
was the first to prove that the circle packings 
in question can be derived from a variational principle
(that is, an energy functional).
In this line of work,
Bobenko \& Springborn \cite{BobenkoSpringborn}
have quite recently discovered
an explicit, elegant and quite general variational principle for
the construction of circle patterns with prescribed intersection
angles.
In the following, we prove the Steinitz theorem based
on their functional --- taking advantage
of all the simplifications that occur in their proof and
formulas if one wants to 
``just'' get the orthogonal circle patterns needed for
the Steinitz theorem.
(See also Springborn~\cite{Springborn04} for an additional
discussion of uniqueness.)

\section{Steinitz' theorem via circle packings}

\begin{theorem}[The Koebe--Andreev--Thurston theorem]
Each $3$-connected planar graph can be realized by
a $3$-polytope which has all edges tangent to the unit sphere.

Moreover, this realization is unique up to Möbius transformations
(projective transformations that fix the sphere).
The edge-tangent realization for which the barycenter of the
tangency points is the center of the sphere is unique
up to ortho\-gonal transformations.
\end{theorem}

\begin{figure}[h]
\begin{center}
\includegraphics[height=105mm]{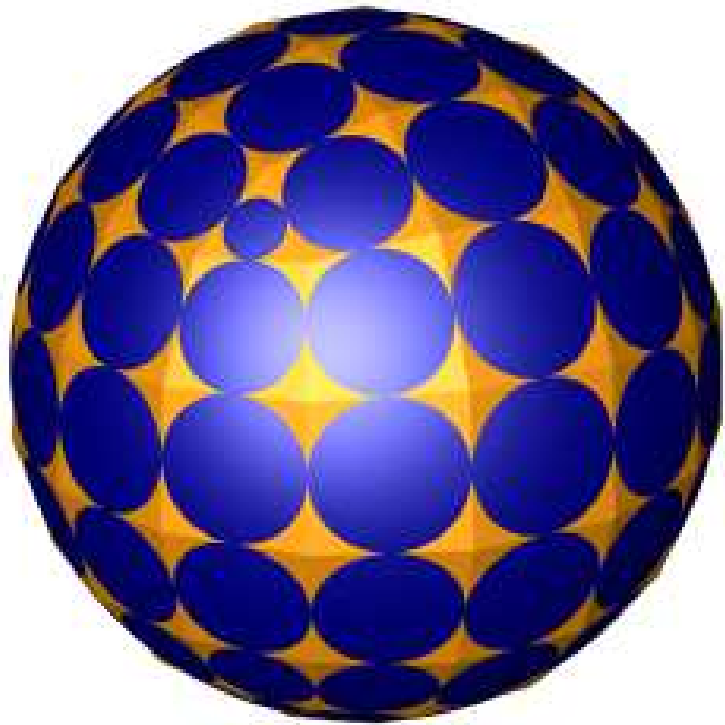}
\end{center}
\caption{Edge-tangent representation of a polyhedron,
according to the Koebe--Andreev--Thurston theorem
[Graphics by Boris Springborn, \textsc{Matheon}]}
\label{figure:Tutte2}
\end{figure}

In our presentation of the proof, we first explain
how any edge-tangent representation of a polytope $P$
induces a circle pattern on the sphere,
which in turn yields a planar circle pattern, 
and the combinatorics of the planar circle pattern yields
a quad graph (a planar graph whose faces are
quadrilaterals), which has $G(P)$ as a subdivided subgraph.
This yields steps (1) to (4) in the following scheme:
\bigskip

\begin{center}
\input{EPS/SBplan.pstex_t}
\end{center}
\bigskip

Our plan is to then reverse this four-step process, 
in order to construct an edge-tangent polytope from the given graph~$G$.
In step (5), 
the quad graph is derived directly from the graph $G=G(P)$,
by superposing the graph with its dual.
Then, in step (6), we construct the rectangular circle pattern
with the combinatorics of the quad graph, and then
proceed to construct~$P$ from it. 

The steps (5), (7), and (8) are quite straightforward:
The key, non-trivial step is~(6), the construction of
the (unique) rectangular circle pattern, which we achieve
via the ``euclidean Bobenko--Springborn functional.''

\begin{proof}
We start with a detailed description of the four-step process from 
edge-tangent polytopes
to planar $3$-connected graphs, via circle packings and quad graphs.
\medskip

\noindent
\textbf{(1).}
Assume that $P\subset\R^3$ is a $3$-polytope whose edges are tangent to
the unit sphere $S^2\subset\R^3$. 
Then the facet planes of~$P$ intersect the unit sphere $S^2$ 
in circles that we call the \emph{facet circles}: We get one
circle for each facet, and the circles are disjoint, but they
touch exactly if the corresponding facets are adjacent.
We also get a second set of circles which we
call the \emph{vertex horizon circles}: Each such circle is the boundary
of the spherical cap consisting of all the points on the sphere that are
``visible'' from the respective vertex. We get one vertex horizon
circle for each vertex, and the circles are disjoint, but they
touch exactly if the corresponding vertices are adjacent.

Moreover, at each edge tangency point, the two touching facet circles
and the two touching vertex horizon circles intersect orthogonally;
see Figure~\ref{figure:BS2} for an example.
(The vertex horizon circles of~$P$
are the facet circles of the dual polytope~$P^*$,
whose edges have the same tangency points as the edges of~$P$;
the facet circles for $P$ are also the vertex horizon circles for~$P^*$;
corresponding edges $e\subset P$ and $e^*\subset P^*$ intersect orthogonally
at the respective tangency point.)

\noindent
\textbf{(2).} We perform a stereographic projection to the plane,
using one of the edge tangency points $p_0$ as the projection center,
and mapping all the facet and vertex horizon circles to the equator
plane corresponding to the projection point.
In the resulting planar figure, the two facet circles through $p_0$
yield two parallel lines (and after a rotation we may assume
that these are horizontal); the two vertex horizon circles
through $p_0$ also yield  two parallel lines, orthogonal
to the first two (and thus vertical). So we get a planar pattern
that consists of four lines bounding an axis-parallel rectangle, and circles
that touch resp.\ intersect orthogonally in the plane.
This is the \emph{rectangular circle pattern}.

If the faces adjacent to the edge $f$ through~$p_0$ are an
$h_1$-gon and an $h_2$-gon, then we get $h_1-2$ resp.\ $h_2-2$ circles
along the horizontal edges of the rectangle.  Similarly, if the end
vertices of~$f$ have degrees $v_1$ and $v_2$, then we get $v_1-2$
resp.\ $v_2-2$ circles along the vertical edges of the rectangle. 
The example that one obtains from the cube (Figure~\ref{figure:BS2})
is displayed in Figure~\ref{figure:BS3}.

\begin{figure}[h]
\begin{center}
\input{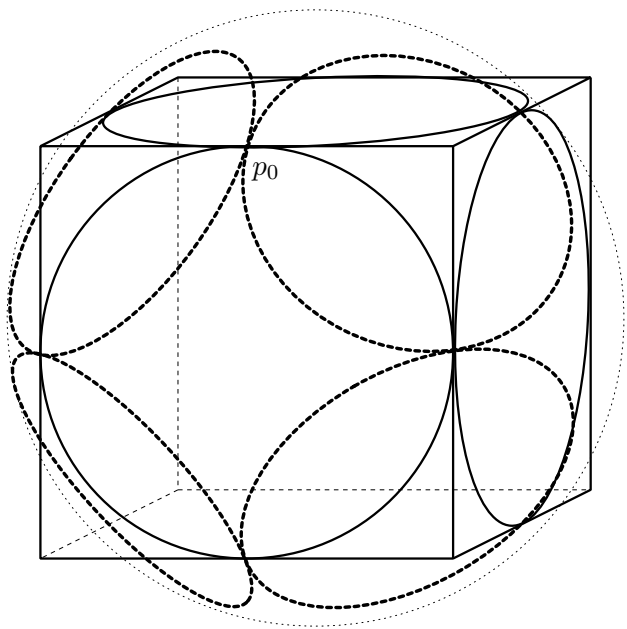}
\end{center}
\caption{The facet circles and the vertex horizon circles (dashed) for
an edge-tangent representation of a regular cube.}
\label{figure:BS2}
\end{figure}

\begin{figure}[h]
\begin{center}
\input{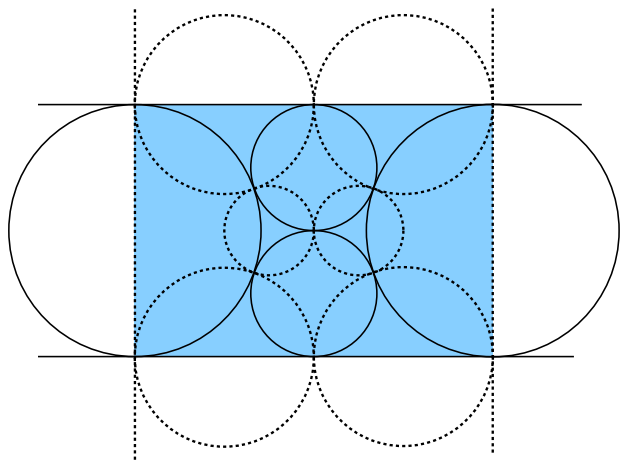}
\end{center}
\caption{The rectangular circle pattern derived from an edge-tangent
  $3$-cube (with $h_1=h_2=4$, $v_1=v_2=3$)}
\label{figure:BS3}
\end{figure}

\noindent\textbf{(3).}
Any rectangular circle pattern yields a quad graph drawing as
follows: The vertex set consists of the centers of all the circles,
with four additional vertices ``far out'' representing the 
four lines that bound the rectangles (as in Figure~\ref{figure:BS4}).
We obtain drawings of both~$G$ and~$G^*$ by 
connecting the centers of touching facet circles resp.\ vertex
horizon circles.
This includes one horizontal edge $f$ of~$G$ ``going through infinity,''
while dual graph $G^*$ has the corresponding edge $f^*$
going through infinity vertically.

From the rectangular circle pattern, we obtain 
a decomposition of a rectangle into
quadrilaterals by connecting the centers of adjacent
facet circles, and the centers of adjacent 
vertex horizon circles.
See the example of Figure~\ref{figure:BS4}, where the rectangle is
shaded. The graph of this rectangle decomposition is
the \emph{quad graph}: Its vertices correspond to (the centers of)
the facet circles that don't contain $p_0$, the vertex horizon circles
that don't contain $p_0$, and intersection points of edges
$e$ and $e^*$ of~$G$ and $G^*$, other than the edges $f,f^*$ that
contain~$p_0$.

\begin{figure}[ht]
\begin{center}
\input{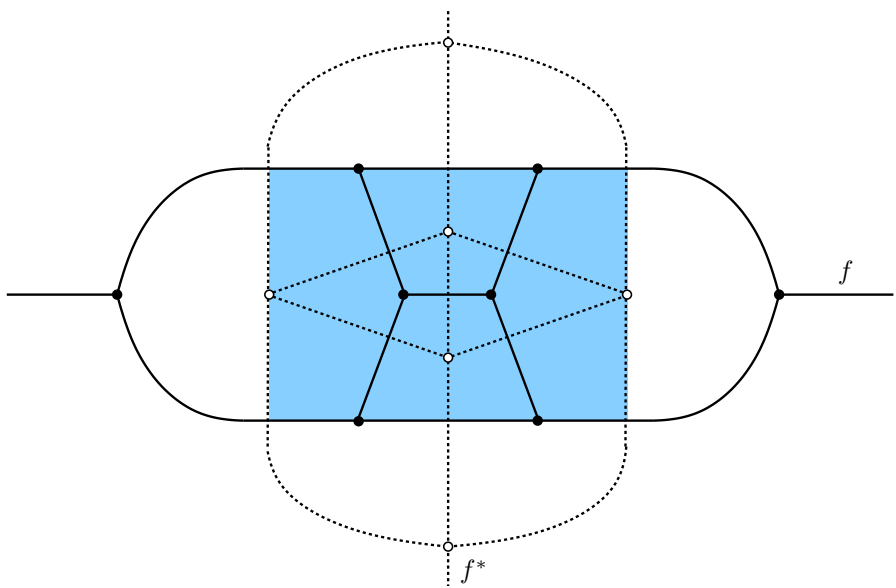} 
\end{center}
\caption{The quad graph for the cube, generated from Figure~\ref{figure:BS3}:
  The white vertices are given by the facet circle centers,
  while the black vertices correspond to the vertex horizon circles;
  the dashed edges connect the centers of adjacent facet circles,
  and the straight edges correspond to adjacent vertex horizon circles.}
\label{figure:BS4}
\end{figure}

\noindent\textbf{(4).} 
In particular, the graph $G$ may be derived from the 
quad graph, by ``deleting the dashed edges.''
\medskip

\noindent
This ends the description of the passage from an edge tangent
polytope to the planar graph drawing. Now we start the way back:
Another four-step process leads us from graphs via quad
graphs and circle patterns to edge-tangent $3$-polytopes.
\medskip

\noindent\textbf{(5).} 
The quad graph may be derived from knowledge
of the graph $G$ alone, plainly by overlaying $G$ and $G^*$. 
For our cube example, the result may look 
like the drawing given in Figure~\ref{fig:quadgraph}.

\begin{figure}[ht]
\begin{center}\vskip-5mm
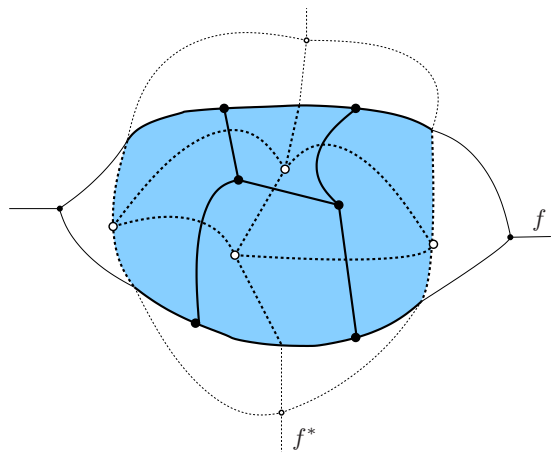\vskip-4mm
\end{center}
\caption{The quad graph for the cube, generated from an overlay
  of the cube graph (black edges) layed out with the edge $f$ ``at
  infinity'' and the dual graph (dashed edges), with the
  dual edge $f^*$ ``at infinity.''
The shaded part defines the restricted quad graph.}
\label{fig:quadgraph}
\end{figure}

The input for the next step will be the \emph{restricted quad graph}:
It is obtained from the full quad graph 
by deleting everything that is adjacent to the original 
edges $f$ and~$f^*$. 
Its bounded faces are quadrilaterals (\emph{quads} for short),
with two black and two dashed edges each.
Each quad has
\begin{compactitem}[~$\bullet$~]
\item a \emph{black} vertex and a \emph{white} vertex\\
  (the black vertex, where the two black edges meet, corresponds to the center
  of a face circle; the white one, where the two dashed edges meet,
  corresponds to the center of a  horizon circle),
\item and two more   vertices where a black and a dashed edge meet\\
  (they correspond to edge tangency points).
\end{compactitem}
For the following, we use $\Izero$ as an indexing set for
the black and white vertices in the restricted quad graph. 
It is in bijection with the vertices of~$G$ and of~$G^*$,
except for the vertices of the edges $f$ and~$f^*$,
which yield lines rather than circles. That is, we have
\[
\Izero\ :=\ V(G-f)\cup V(G^*-f^*).
\]

\medskip
\noindent
The following step, which takes us from combinatorics
(a graph drawing) to geometry (a circle pattern), is the crucial one. 
\medskip

\noindent
\textbf{(6).}
In the ``correct'' realization of the restricted quad graph, which
would yield a circle packing, each quad is drawn as a \emph{kite}
in which
\begin{compactitem}[~$\bullet$~]
\item the two black  edges have the same length\\
  (radius $r_i$ of the corresponding vertex horizon circle), 
\item the two dashed edges have the same length\\
  (radius $r_j$ of the corresponding facet circle), 
\item and there are two right angles 
  between black and dashed edges\\
  (where facet and vertex horizon circles intersect).
\end{compactitem}

\noindent
The kites have to look like the one in Figure~\ref{figure:BS5}.

\begin{figure}[ht]
  \begin{center}
\input{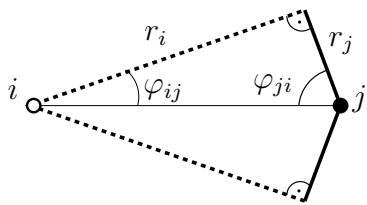}
\end{center}
\vskip-4mm
\caption{A kite, with radii $r_i=e^{\rho_i}$, $r_j=e^{\rho_j}$, and angles
$\varphi_{ij}$ and $\varphi_{ji}$}
\label{figure:BS5}
\end{figure}

\noindent
Hence, we have to solve the following construction problem:
\begin{quote}
Given a quad graph decomposition of a rectangle, derived from the
overlay of a $3$-connected planar graph $G$ and its dual $G^*$,
construct a geometric drawing, with straight edges,
as a kite decomposition of a rectangle.
\end{quote}
The kites are completely determined if we know their
edge lengths: If the edge lengths in a kite
are $r_i,r_j>0$, then the angles are given by
\[
\varphi_{ij}\ =\ \arctan\Big(\frac{r_j}{r_i}\Big)\quad\textrm{and}\quad
\varphi_{ij}\ =\ \arctan\Big(\frac{r_i}{r_j}\Big),
\]
with $\varphi_{ij}+\varphi_{ji}=\frac\pi2$ (see Figure~\ref{figure:BS5}).
Thus all we have to do is to determine radii~$r_i$ corresponding to
the black and white vertices of the quad graph, such that 
the following system of equations is satisfied:
\begin{equation}\label{eq:nonlinear}
\fbox{$\quad\displaystyle
\sum^{\strut}_{\textstyle j: i\kdiamond j}2\arctan\Big(\frac{r_j}{r_i}\Big)
\ \ =\ \ \Phi_i
\qquad\textrm{for all vertices $i\in \Izero$},$\quad}
\end{equation}
where the right-hand-sides are given by
\[
\Phi_i\ :=\ 
\begin{cases}
\ \ \pi& \textrm{if $i$ is on the boundary,}\\
\  2\pi& \textrm{if $i$ is in the interior.}
\end{cases}
\]
In the equation whose right hand side is $\Phi_i$, the sum 
on the left hand side is taken over all vertices $j\in\Izero$
that are opposite to $i$ in one of the kites.
(If $i$ is a white vertex, then $j$ will be black,
and vice versa.)

Indeed, if \eqref{eq:nonlinear} is satisfied, then we can easily
construct the kites and piece them together to get a flat rectangle
and the circle packing.
Badly enough, (\ref{eq:nonlinear}) is a non-linear system
of equations, which we have to solve in positive variables $r_i>0$.
We want to know that this has a solution, which is unique up
to multiplying all the $r_i$s with the same factor, 
and which can be computed efficiently.
Luckily, we can do this, since the system is solved by 
minimizing an explicit and easy-to-write-down ``energy'' functional
which will turn out to be convex, with a unique minimum.
For this, we first do a change of variables,
\[
\rho_i\ :=\ \log r_i.
\]
Then we normalize by the condition $\prod_ir_i=1$, that is, 
\[
\sum_i \rho_i\ =\ 0.
\]
Furthermore, we define 
\[
f(x)\ :=\ \arctan(e^x).
\]
This auxiliary function is graphed in Figure~\ref{figure:f}.
Note that $f(-x)=\frac\pi2-f(x)$.

\begin{figure}[ht]
\begin{center}  
  \psfrag{0.5}{$0.5$}
  \psfrag{1.0}{$1.0$}
  \psfrag{1.5}{$1.5$}
  \psfrag{2}{$\raisebox{-3pt}{2}$}
  \psfrag{4}{$\raisebox{-3pt}{4}$}
  \psfrag{-2}{$\raisebox{-3pt}{$-2$}$}
  \psfrag{-4}{$\raisebox{-3pt}{$-4$}$}
  \psfrag{auxiliary}{\qquad\qquad$\arctan(e^x)$}
  \psfrag{x}{$ x$}
\includegraphics[height=40mm]{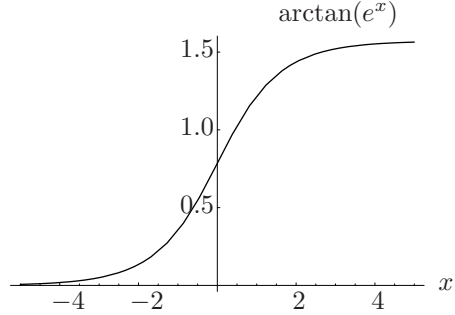}
\end{center}
\vskip-5mm
\caption{$f(x)=\arctan e^x$}
\label{figure:f}
\end{figure}

\noindent
We differentiate ~$f$,
\[
f'(x)\ =\ \frac1{1+e^{2x}}e^x\ =\ \frac1{2\cosh x},
\] 
which yields $f'(-x)=f'(x)>0$ for all $x\in\R$.
We also integrate $f$, and define
\[
F(x)\ :=\ \int_{-\infty}^x f(t) \mathrm{d}t.
\]
This function satisfies $F(x)\ge0$ for all $x$,
but also $F(x)\ge\frac\pi2 x$.
Thus we get that 
\begin{equation}
  \label{eq:F-F}
  F(x)+F(-x)\ \ge\ \tfrac\pi2 |x|.
\end{equation}
The system (\ref{eq:nonlinear}) we have to solve may be rewritten 
in terms of $f(x)$ as
\begin{equation}\label{eq:nonlinear2}
\sum_{\textstyle i: i\kdiamond j} 2 f(\rho_j-\rho_i)\ \ =\ \ 
\Phi_i \qquad\textrm{for all black or white vertices, $i\in \Izero$.}
\end{equation}
To solve this, Bobenko \& Springborn \cite{BobenkoSpringborn}
present the functional
\begin{equation}\label{eq:BSfunctional}
\fbox{\quad$\displaystyle
\mathrm{BS}(\mathbf{\rho})\ :=\ 
\sum^{\strut}_{\strut\textstyle i\kdiamond j}
\Big\{F(\rho_j{-}\rho_i) + F(\rho_i{-}\rho_j) - \tfrac\pi2(\rho_i+\rho_j)\Big\}
\ +\sum_{\textstyle i\in\Izero}\Phi_i\rho_i,
$\quad}
\end{equation}
where the first sum is over all \emph{unordered} pairs $\{i,j\}$ of
vertices $i,j\in\Izero$ that are opposite in one of the kites.
The claim is now that
\begin{compactenum}[(A)]
\item
the critical points of $\mathrm{BS}(\mathbf{\rho})$ are exactly the
solutions to our system~(\ref{eq:nonlinear2}),
\item
the functional is convex: Restricted to $\sum_i\rho_i=0$ 
it is strictly positive definite, 
so the critical point is unique if it exists, and
\item
the functional gets large if any of the differences $\rho_i-\rho_j$
gets large: Thus the functional must have a critical point (a minimum) ---
the solution we are looking for.
\end{compactenum}
For (A), 
a simple computation yields the gradient of~$\mathrm{BS}({\rho})$:
\[
\frac{\partial\mathrm{BS}({\rho})}{\partial\rho_i}\ =\ 
\Phi_i - \sum_{\textstyle i\kdiamond j} 2f(\rho_j-\rho_i).
\]
Thus the critical points of~$\mathrm{BS}({\rho})$ are exactly the
solutions to (\ref{eq:nonlinear2}).

For (B), we compute the Hessian (the matrix of second derivatives)
for $\mathrm{BS}({\rho})$, and find that
\[
{x}^T \mathrm{BS}({\rho})'' {x}\ =\ 
2 \sum_{\textstyle i\kdiamond j} f'(\rho_j-\rho_i)\,(x_j-x_i)^2.
\]
We know that $f'(\rho_j-\rho_i)>0$, so this quadratic form can vanish
only if all the differences $x_j-x_i$ vanish for 
``adjacent'' $i,j\in\Izero$ (that is, for black/white
vertices that share a kite). But the graph
we consider is connected, so this implies that all variables
$x_i$ are equal. Restricted to $\sum_i x_i=0$ this
yields that all $x_i$ vanish, so the Hessian is positive definite
on the restriction hyperplane,
and the solution we are striving for is unique if it exists.

To prove the existence claim (C), 
we have to find that $\mathrm{BS}({\rho})$
grows large if any difference of variables $\rho_k-\rho_i$ gets large.
With the same argument we just used this implies that some
difference of ``adjacent'' variables
will become large.
Then also 
$F(\rho_j{-}\rho_i)+F(\rho_i{-}\rho_j)\ge\frac\pi2|\rho_j-\rho_i|$ gets large,
but it will grow only linearly in $|\rho_j-\rho_i|$, and it is not
obvious that the growing positive terms in (\ref{eq:BSfunctional})
will ``outrun'' the negative terms. This will require a 
careful ``matching'' between positive and negative terms.

To achieve this, we use the existence of a \emph{coherent angle system},
that is, an assignment of angles $\varphi_{ij},\varphi_{ji}>0$
to the kites that satisfies the conditions
\begin{equation}  \label{eq:coherent}
\varphi_{ij}+\varphi_{ji}\ \ =\ \ \tfrac\pi2 \qquad\textrm{and}\qquad
\sum_{\textstyle j:i\kdiamond j}   2\varphi_{ij}\ \ =\ \ \Phi_i.
\end{equation}
Any solution to (\ref{eq:nonlinear}) would give us 
a coherent angle system, but the existence of such a coherent
angle system is much weaker, far from solving the system
(\ref{eq:nonlinear}): If we have a coherent angle system,
then we could construct kites from this --- whose
angles would fit together at the black and white vertices, but
whose side lengths might not. (Compare
Figure~\ref{figure:coherentanglesystem}.)

For any coherent angle system,
$\eps_0:=\min\limits_{k,\ell}\varphi^{}_{k\ell}$ is a positive
number.

\subsubsection*{If there is a coherent angle system, then the minimum exists.}
Let's assume for now that a coherent angle system exists (this will
be proved below). Then 
\begin{eqnarray*}
\mathrm{BS}(\mathbf{\rho})
&= &
\sum_{i\kdiamond j}
\Big\{F(\rho_j{-}\rho_i) + F(\rho_i{-}\rho_j) - \tfrac\pi2(\rho_i+\rho_j)\Big\}
\ +\ 
\sum_{i}\Phi_i\rho_i \notag
\\
&\stackrel{\textrm{(i)}}>  &
\sum_{i\kdiamond j}
\Big\{\tfrac\pi2\big|\rho_i-\rho_j\big| - \tfrac\pi2(\rho_i+\rho_j)\Big\}
\hspace{18mm}\ +\  
\sum_{i}\Phi_i\rho_i 
\\
&\stackrel{\rm(ii)}=  &
\sum_{i\kdiamond j}
\Big\{\tfrac\pi2\big|\rho_i-\rho_j\big| - \tfrac\pi2(\rho_i+\rho_j)\Big\}
\ +\ 
\sum_{i\kdiamond j}
2(\varphi_{ij}\rho_i + \varphi_{ji}\rho_j)
\\
&\stackrel{\rm(iii)}=  &
\sum_{i\kdiamond j} -\pi \min\{\rho_i,\rho_j\}
\hspace{12mm}\ +\ 
\sum_{i\kdiamond j}
2(\varphi_{ij}\rho_i + \varphi_{ji}\rho_j)
\\
&\stackrel{\rm(iv)}\ge&
\sum_{i\kdiamond j} -\pi \min\{\rho_i,\rho_j\}
\hspace{2mm}\ +\ 
\sum_{i\kdiamond j}
\pi\min\{\rho_i,\rho_j\}+ 2\min\{\varphi_{ji},\varphi_{ij}\}|\rho_i-\rho_j|
\\[1mm]
&=  & 
\sum_{i\kdiamond j}
2\min\{\varphi_{ji},\varphi_{ij}\}|\rho_i-\rho_j|
\ \ \ge\ \ 
2\eps_0\,\sum_{i\kdiamond j}|\rho_i-\rho_j|.
\end{eqnarray*}

\noindent
Here 
\begin{compactitem}[$\bullet$]
  \item the estimate for (i) uses
    $F(x)+F(-x)\ge\tfrac\pi2|x|$,
    which is (\ref{eq:F-F}).
  \item (ii) is obtained by substituting \eqref{eq:coherent}. 
  We need the second term in the second sum in~(ii)
  since the sums over ``$i\kdiamond j$'' are 
  sums over unordered pairs; there is no extra summand for~``$j\kdiamond i$.'' 
  \item (iii) follows from
$|x-y|-(x+y)= -2\min\{x,y\}$, 
  \item For (iv), in the case $\rho_j\ge\rho_i$ we compute
\begin{eqnarray*}
2(\varphi_{ij}\rho_i + \varphi_{ji}\rho_j)&=&
\pi\rho_i - 2\varphi_{ji}\rho_i + 2\varphi_{ji}\rho_j\\
&=&
\pi\min\{\rho_i,\rho_j\} + 2\varphi_{ji}|\rho_i-\rho_j|\\
&\ge&
\pi\min\{\rho_i,\rho_j\}+2\min\{\varphi_{ji},\varphi_{ij}\}|\rho_i-\rho_j|,
\end{eqnarray*}
and analogously for $\rho_i\ge\rho_j$.
\end{compactitem} 
We are dealing with a connected quad graph.
Thus if the norm of the vector $\rho$ gets large, 
while the sum of the $\rho_i$ is zero, then
also for two 
$i,j\in\Izero$ in the same quadrilateral the difference $|\rho_i-\rho_j|$
gets large. Thus by the computation above,
$\mathrm{BS}(\rho)>2\eps_0|\rho_i-\rho_j|$ gets large.
This is sufficient to prove that the strictly convex function 
$\mathrm{BS}(\rho)$ does have a (unique) minimum --- the solution
to our problem.

\begin{figure}[t]
\begin{center}
\input{EPS/quadgraph3b.pstex_t}%
\end{center}
\vskip-3mm
\caption{%
The assignment in this
figure is a coherent angle system -- but not one that corresponds
to a correct circle pattern.
\newline
(Note that the construction of the coherent angle system proceeds
from the plane graph without use of a straight edge drawing.
In the figures further down we draw the graphs with straight
edges for simplicity, but this structure is not used in the proof.
Rather, it is produced by the proof.)}
\label{figure:coherentanglesystem}
\end{figure}

\subsubsection*{A coherent angle system exists}
Finally, we have to verify the existence of a coherent angle system.
We will see here that via some simple network flow
theory, this follows from an expansion property in
the ``diagonal graph'' $D(G\cup G^*)$.
After that, we will prove the expansion property.

Let $G$ be a $3$-connected planar graph, $G^*$ its dual,
both of them again drawn into the plane with dual edges $f,f^*$
intersecting ``at infinity.''
Then the \emph{diagonal graph} $D=D(G\cup G^*)$ has the same vertex set
as $G\cup G^*$. Its edges correspond to the diagonals
in the quad graph given by $G\cup G^*$.

Equivalently, the diagonal graph $D$ has black vertices
corresponding to the vertices of~$G$, and white vertices
corresponding to the faces of~$G$. The edges of $D$ correspond to
the vertex--face incidences of~$G$. 
See Figure~\ref{figure:diagonalgraph} for an example.

The \emph{reduced diagonal graph} $D'=(V',E')$ is obtained from the
diagonal graph $D=(V,E)$ by removing the
two vertices of~$f$, the two vertices of~$f^*$, and the four edges that
connect them, but none of the others. So
indeed, $D'$ does have pending edges (half-edges)
which have lost one of their end-vertices.%
\footnote{%
I am sure you won't be troubled too much by the fact that
this is not a graph in the usual technical sense, since
it does have half-edges with only one end-point.}
See Figure~\ref{figure:diagonalgraph-reduced} for an example.

\begin{figure}[ht]
\begin{center}
\input{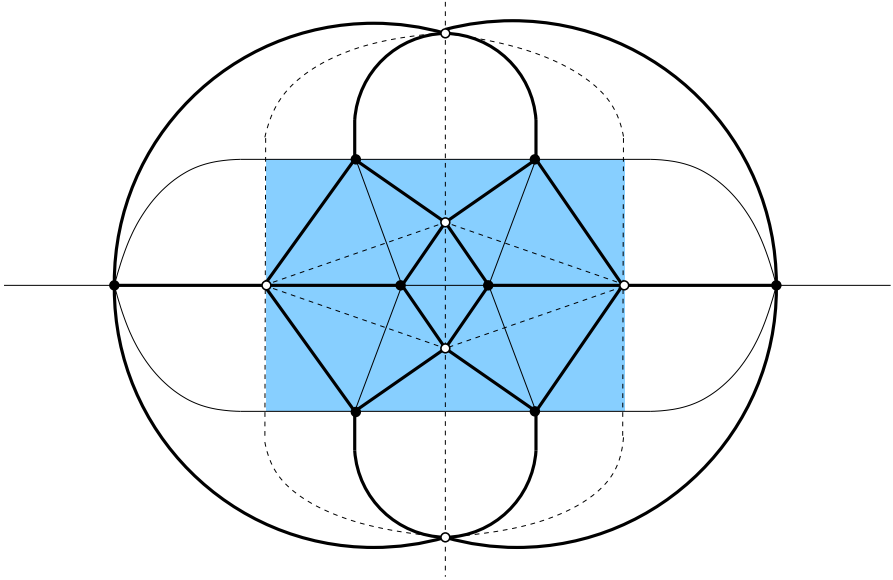}
\end{center}
\caption{The diagonal graph $D=D(G\cup G^*)$, given by the
  fat edges, where $G$ is the graph of the cube}
\label{figure:diagonalgraph}
\end{figure}

\begin{figure}[ht]
\begin{center}
\input{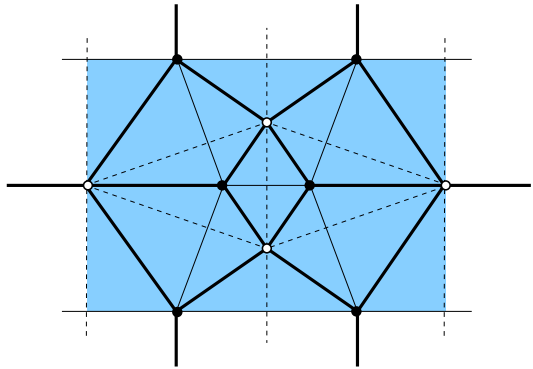}
\end{center}
\caption{The fat edges in this figure display the
  reduced diagonal graph $D'=D'(G\cup G^*)$ in the case
  where $G$ is the graph of the cube, derived from
  Figure~\ref{figure:diagonalgraph}. Note that the 
  fat edges leaving the rectangle are included in $D'$, their
  vertices at the other end are not. So in this example $D'$ has
  $10$ vertices and $20$ edges, including $6$ half-edges with only
   one end-vertex.}
\label{figure:diagonalgraph-reduced}
\end{figure}

The diagonal graph $D=D(V,E)$ is a quad graph: All its faces,
including the ``unbounded'' face (if we draw it in the plane)
are quadrilaterals. From this, we get by double counting
that $|F|=2|E|$ and thus $|V|=2|E|-4$ by Euler's relation.
The reduced quad graph $D'=(V',E')$
has $|V'|=|V|-4$ vertices and $|E'|=|E|-4$ edges. Hence we get $|E'|=2|V'|$:
The reduced quad graph has exactly double as many edges as vertices.


The concept of a coherent angle system has a very nice
interpretation in terms of the restricted diagonal graph:
Each vertex $v_i$ gets a weight of $2\pi$, and this has
to be distributed to the edges $e$ incident to $v_i$ such that
\begin{compactitem}[~$\bullet$~]
  \item each edge $e$ incident to $v_i$ gets a positive part of the
    weight $2\pi$ of~$v$, 
  \item all of the weight $2\pi$ of $v_i$ is distributed to its incident
    edges, and
  \item the weights assigned to each edge sum to~$\pi$.
\end{compactitem}
Indeed, in such an assignment any half-edge clearly gets 
a weight of~$\pi$ from its only end-vertex, which 
corresponds to a boundary vertex of the restricted quad graph;
thus the boundary vertex $v_i$ distributes a weight of exactly~$\Phi_i=\pi$
to its other incident edges, that is, to the (diagonals of the) 
kites it is incident to. The vertices of~$D'$ without
an incident half-edge correspond to interior vertices $v_j$
of the restricted quad graph, so they have a weight/angle
of~$\Phi_j=2\pi$ to distribute to the incident edges/kites.
\medskip

The ``weight distribution problem'' for the reduced diagonal
graph~$D'=(V',E')$ may also be interpreted as a flow problem 
(cf.~\cite{AhujaMagnantiOrlin}):
We have to find a maximal flow, of weight $2\pi|V'|=\pi|E'|$,
in a two-layer network as depicted in Figure~\ref{figure:network}.
It consists of a source node~$s$, then a layer of nodes formed 
by the vertex set $V'$ of~$D'$,
then a layer of nodes in bijection to the
the edge set $E'$, and then the sink node~$t$.
There are three groups of arcs:
The arcs $(s,v')$ emanating from the source all have
an upper bound of $2\pi$; 
the arcs of type $(v',e')$, where the edge $e'$ is incident to~$v'$,
get an upper bound of~$\infty$, 
while the arcs at the sink, $(e',t)$, have an upper bound of~$\pi$.

We need a \emph{positive} flow in this network;
to get this, we put a small lower bound of $\eps>0$ on each edge
of type~$(v',e')$, and $0$ on all other edges.
There is a \emph{feasible flow} in this network with upper and lower
bounds on each edge: For this, 
let the flow value be~$\eps$ on each $(v',e')$-arc,
and a suitable multiple of~$\eps$ on the other arcs.

\begin{figure}[ht]
\begin{center}
\input{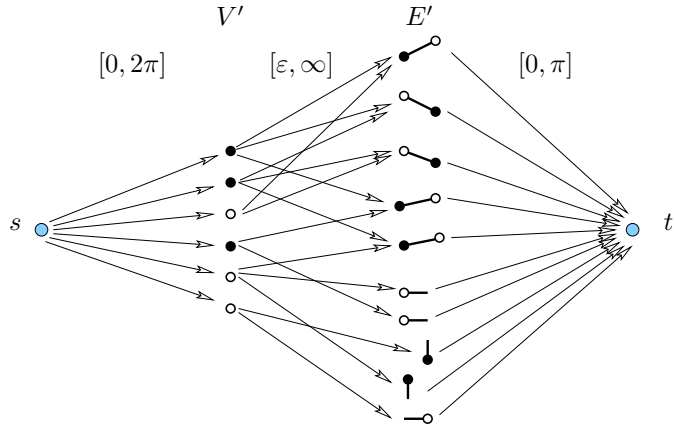}
\end{center}
\vskip-2mm
\caption{Construction of a coherent angle system 
from a network flow problem with lower and upper bounds,
which are indicated by intervals like~$[0,\pi]$.}
\label{figure:network}
\end{figure}

\begin{figure}[ht]
\begin{center}
\input{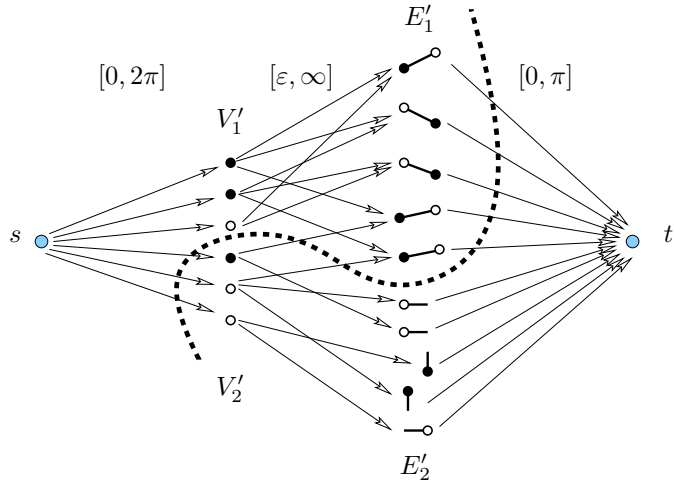}
\end{center}
\vskip-2mm
\caption{The dashed line indicates the cut 
$[\{s\}\cup V'_1\cup E'_1,\, V'_2\cup E'_2\cup \{t\}]$ in our network}
\label{figure:network2}
\end{figure}

We need a positive flow of value $2\pi|V'|=\pi|E'|$
in this network. There is a feasible flow,
and no flow with a larger value than $2\pi|V'|$ can exist
due to the cuts that separate $s$ or $t$ from the rest of
the network. Thus we can apply the following generalization of the
Max-Flow Min-Cut Theorem on network flows.
(You should prove this yourself: See \exercise{ex:MFMC}.)
\smallskip

\begin{theorem}[Generalized Max-Flow Min-Cut Theorem; 
cf.~{\cite[Sect.~6.7]{AhujaMagnantiOrlin}}]~\\
If an $(s,t)$-network with lower and upper bounds 
has a feasible flow, then the value of a maximal $(s,t)$-flow
is the capacity of a minimal $(s,t)$-cut.
\end{theorem}

The \emph{capacity} of an $(s,t)$-cut in a network with upper and lower
bounds is the sum on the upper bounds of the forward arcs,
minus the sum of the lower bounds on the backward arcs
across the cut. So in our example the cuts
$[\{s\},V'\cup E'\cup \{t\}]$ and
$[\{s\}\cup V'\cup E',\{t\}]$ have capacity $2\pi|V'|=\pi|E'|$.
Could there be a cut of smaller capacity?
Any $(s,t)$-cut is of the form
\[
[\{s\}\cup V'_1\cup E'_1,\, 
           V'_2\cup E'_2\cup \{t\}]
\]
for partitions $V'=V'_1\uplus V'_2$ and $E'=E'_1\uplus E'_2$.
Such a cut has finite capacity if there are no arcs
$(v',e')$ from~$V'_1$ to~$E'_2$; compare Figure~\ref{figure:network2}.
That is, we should take $E'_1$ to include all the edges
that are incident to a vertex in $V'_1$.

The capacity of the cut
$[\{s\}\cup V'_1\cup E'_1,\, V'_2\cup E'_2\cup \{t\}]$
is 
\[
  2\pi|V'_2| + \pi|E'_1| - \eps|A(V_2',E'_1)|
\ \ =\ \ 2\pi|V'| 
- 2\pi|V'_1| + \pi|E'_1| - \eps|A(V_2',E'_1)|,
\]
where $|A(V_2',E'_1)|$ denotes the number of arcs from~$V'_2$ to~$E'_1$.
For small enough $\eps$, say $\eps=1/|A(V',E')|$, we have
$\eps|A(V_2',E'_1)|<1$. Thus the following ``expansion property'' for the
diagonal graph implies that 
all cuts $[\{s\}\cup V'_1\cup E'_1,\, V'_2\cup E'_2\cup \{t\}]$
have capacity larger than $2\pi|V'|$,
except in the two trivial cases given as examples above,
where the capacity is exactly $2\pi|V'|$.
Thus the maximal flow, of value $2\pi|V'|$, exists; it is positive,
and yields the coherent angle system.

\subsubsection*{Expansion in the diagonal graph.}
It remains to verify the following:
  Let $V'_1\subseteq V'$ be a set of vertices
  in the reduced diagonal graph $D'(G\cup G^*)=(V',E')$,
  and assume that $E'_1\subseteq E'$ includes
  all edges of~$D'$ that are incident to a vertex in~$V'_1$.
  Then
  \begin{equation}
    \label{eq:matching}
    |E'_1|\ \ \ge\ \  2|V_1'|,
  \end{equation}
  with equality only in the trivial cases $V'_1=\emptyset$ and $V'_1=V'$.

  For this we may assume that the 
  subgraph induced by $V'_1$ is connected, because
  we can consider its components separately.
  We may also assume that $|V_1|\ge2$, so $V'_1$ contains both
  a black and a white vertex.

  Now let $U$ be an open subset of the plane (or of $S^2$)
  whose boundary curves separate $V'_1$ from the $h+1$ components
  of the graph $D\setminus V'_1$,
  as illustrated in Figure~\ref{figure:DiscWithHoles}.
  Topologically, $U$ is an open disk with $h\ge0$ holes.

  The diagonal graph yields a cell decomposition of~$U$, consisting
  of $f_0=|V'_1|$ vertices,
  $f_1^\inter$ interior edges, $f_1^\bdy$ other (half-)edges,
  $q$ quadrilateral faces, and $b_1+b_2+b_3$ boundary faces,
  where $b_i$ counts the faces with $i$ vertices in~$I'$. 
  In particular the total number of edges is
  $f_1=f_1^\inter+f_1^\bdy=|E'_1|$,

\begin{figure}[ht]
\begin{center}
\input{EPS/cube-circ6c.pstex_t}
\end{center}\vskip-5mm
\caption{An example of five vertices in the
 reduced diagonal graph of Figure~\ref{figure:diagonalgraph-reduced}.
 The neighborhood $U$ is shaded.\newline 
  $f_0=|V'_1|=5$, $f_1=|E'_1|=14$, $f_1^\inter=4$, $f_1^\bdy=10$, 
 $h=0$, $q=0$, $b_1=4$, $b_2=4$, $b_3=2$.}
\label{figure:DiscWithHoles}
\end{figure}

  Double counting the edge-face incidences yields
  \begin{equation}
    \label{eq:EdgeDoubleCount}
    2f_1^\inter\ =\ 4q  + b_2+2b_3\qquad\textrm{and}\qquad
    2f_1^\bdy\ =\ 2b_1+2b_2+2b_3.
  \end{equation}

The Euler characteristic of $U$ is
  \begin{equation}
    \label{eq:EulerChar}
    1-h\ \ =\ \ f_0-f_1+q+b_1+b_2+b_3
       \ \ =\ \ f_0-f_1^\inter+q.
  \end{equation}
With this we get
\begin{eqnarray*}
  |E'_1|-2|V'_1|\ =\  f_1 - 2f_0 
&\stackrel{\eqref{eq:EulerChar}}=&
      (f_1^\inter+f_1^\bdy)\ -\ 2(f_1^\inter-q+1-h) \\
&=&
      f_1^\bdy-f_1^\inter+2q+2h-2 \\
&\stackrel{\eqref{eq:EdgeDoubleCount}}=&
      (b_1+b_2+b_3) - (2q+\tfrac12b_2+b_3)\ +\ 2q+2h-2 \\
&=&
      \tfrac12(2b_1+b_2-4)\ +\ 2h.
\end{eqnarray*}
To conclude that $|E'_1|-2|V'_1|\ge0$, with equality
only if $V'_1=V'$, we use $h\ge0$, and need to verify
that $2b_1+b_2\ge4$ holds, with equality only in the trivial case~$V'_1=V'$.

For this we count the vertices $v$ of $D\setminus V'_1$ which are
\emph{adjacent to~$V'_1$}, that is, such that
some quad in
the full quad-graph $D$ contains both $v$ and a vertex from~$V'_1$.
Walking along the boundary curves of~$U$, and exploring the
quads that we traverse that way, we see 
that there are not more than $2b_1+b_2$ such vertices~$v$:
We find at most two new vertices in any quad that contains a
boundary cell with $1$ vertex in~$V'_1$,
and at most one new vertex in the quad of
a boundary cell with $2$ vertices in~$V'_1$.
The vertices found during the walk need not be all distinct,
and some may not even lie outside~$V'_1$
(compare Figure~\ref{figure:DiscWithHoles2}). Thus we get only an inequality,
\[
2b_1+b_2\ \ \ge\ \ \#\{\textrm{vertices of $D\setminus V'_1$ 
adjacent to~$V'_1$}\}.
\]
In the boundary of each ``hole'' of $U$ 
we will discover at least one vertex of $D'\setminus V'_1$. 
In the outer face during our walk we even discover a cycle of~$D$
(see Figure~\ref{figure:DiscWithHoles2}).
Since $D$ is bipartite, this cycle has even length.
In the trivial case of $V_1'=V'$ this is exactly the
$4$-cycle $C'$ given by $D\setminus D'$. 
If $V_1'\neq V'$, then the vertices we discover either
yield the cycle $C'$ plus additional vertices,
or we find a different cycle. But any cycle other than $C'$
must have at least $6$ vertices: Indeed, it is an even cycle,
on which black and white vertices alternate. The black vertices
on the cycle either include both the vertices of~$f$, or
with respect to the original graph $G$ they separate a 
black vertex in $V'_1$ from a vertex of~$f$;
from the $3$-connectivity of $G$ we thus get that the cycle
contains at least three black vertices, that is, at least
$6$ vertices in total. The same holds for the white vertices,
the dual graph $G^*$, which is also $3$-connected, and the vertices of~$f^*$.
Thus
\[
\#\{\textrm{vertices of $D$ in the boundary of~$U$}\}\ \ \ge\ 4,
\]
with equality only if $V'_1=V'$.
This completes the proof for the expansion property,
and thus for the existence of a coherent angle system, and of the 
circle packing.

\begin{figure}[h]
\begin{center}
\input{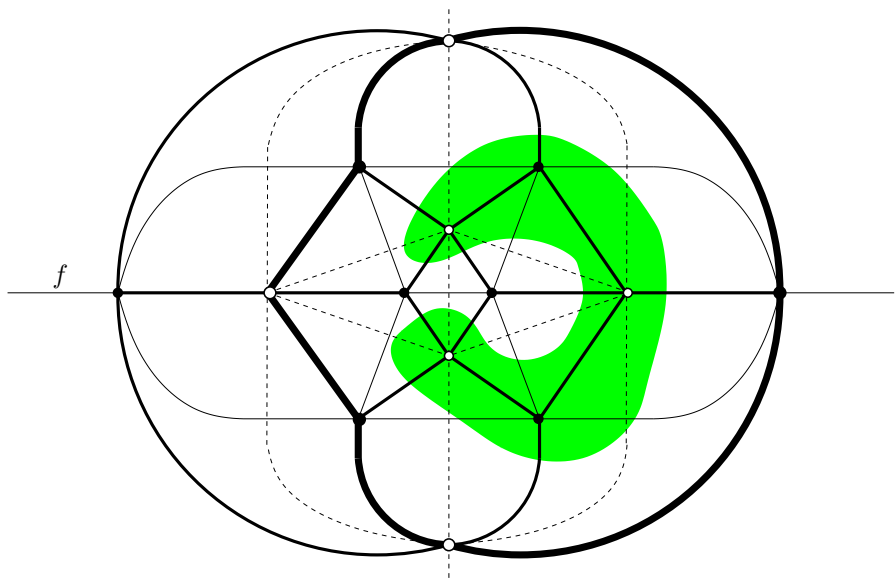}
\end{center}
\vskip-5mm
\caption{%
The cycle in the outer face to be discovered during the walk along the
boundary curve of~$U$ is drawn with fat edges; it is a $6$-cycle. 
With respect to $G$, which is drawn
in thin black lines, the
three vertices of the $6$-cycle separate
a vertex of $f$ from the two black vertices in $V'_1$.}
\label{figure:DiscWithHoles2}
\end{figure}

\noindent
\textbf{(7), (8).}
Given a correct rectangular circle pattern,
it is easy to reconstruct the spherical circle pattern (via an
inverse stereographic projection). From this, we
obtain the edge-tangent polytope: Its face planes are given
by the facet circles (and its vertices are given by the
cone points for which the vertex horizon circles do indeed
appear on the horizon).
Thus construction steps (7) and (8) are easy --- the hard part was~(6).
\end{proof}

Is this the perfect proof? I think it is really nice,
but still one could dream of a proof that avoids the 
stereographic projection, and produces the circle packing
directly from some functional on the sphere~\ldots.

\subsection*{Exercises}
\begin{compactenum}[\thechapter.1.]
\item\label{ex:triangle}
Show that each $3$-polytope has a triangle face, or a simple vertex
(a vertex of degree~$3$), or both.
Even stronger, show that the number of triangle faces plus 
the number of simple vertices is at least eight, so there 
are at least four triangle faces, or at least four simple vertices.
\\
Hint: Use the Euler equation.
\item
Prove that each $3$-polytope has two faces with the same number of vertices.
\\
Hint: Do not use the Euler equation.
\item\label{ex:Steinitzlemma}
Prove the Steinitz Lemma \ref{lemma:Steinitz06}:
 \begin{compactitem}[--]
 \item Prove the ``upper bound theorem'' for dimension~$3$, that is,
   that $f_2\le2f_0-4$ (you may use Euler's equation), and derive
   $f_0\le2f_2-4$ by duality.
 \item Compute the $f$-vectors of the pyramids over $n$-gons.
 \item How does $(f_0,f_2)$ change if you stack a pyramid onto a
   triangle $2$-face, or if you truncate a simple vertex?
 \end{compactitem}
\item
If a $3$-dimensional polytope has $f_1=23$ edges,
how many vertices/faces can it have?
Construct an example for each possible pair~$(f_0,f_2)$.
\item\label{ex:fslope}
Alternative homogeneous coordinates for the
cone of $f$-vectors are given by the ``imbalance''
$\sigma:=\tfrac{f_2-f_0}{f_1-6}$, where
the self-dual term $f_1-6$ measures the ``size.''
\hfill\break
Show that $-\tfrac13\ \le\ \sigma\ \le\ +\tfrac13$,
where $\sigma=\pm\frac13$ characterizes simple resp.\ simplicial polytopes.
\item
Characterize the possible $(f_0,f_2)$-pairs for \emph{cubical} $3$-polytopes,
that is, for all polytopes with quadrilateral $2$-faces only.
Where are the $(f_0,f_2)$-pairs of cubical $3$-polytopes
in Figure~\ref{figure:Steinitz06}?

How about $3$-polytopes with pentagon faces only? 
Hexagon faces only?
\item
Construct quad graphs and the planar circle patterns for 
\begin{compactenum}[(a)]
\item a square pyramid,
\item a cube/octahedron,
\item a cube with vertex cut off,
\item a dodecahedron.
\end{compactenum}
Which of the circle patterns do you get with rational coordinates?
\item\label{ex:CyclicStacked3d}
 Show that every $3$-dimensional cyclic polytope $C_3(n)$
is a stacked polytope.
(However, $C_d(n)$ is not stacked, for $d\ge4$ and $n\ge d+2$.)
\item\label{ex:MFMC}
Describe a computational procedure
to construct a coherent angle system:
For this use a scheme to augment flows along undirected paths
in the network with lower and upper bounds
(increasing the value along forward arcs, decreasing the 
values on backward arcs).
Your procedure should also imply a proof for the
Generalized Max-Flow Min-Cut Theorem 
\cite[Thm.~6.10, p.~193]{AhujaMagnantiOrlin}.
\end{compactenum}

\lecture{Shapes of {\itshape {f}}\/-Vectors}

Let's look at the $f$-vectors of $d$-dimensional convex polytopes $P$,
where the dimension $d$ is really large.  Any such $f$-vector
\begin{eqnarray*}
f(P)&=&(f_0, f_1,f_2,
\qquad\quad\ldots\qquad\quad\ldots\quad\qquad\ldots\quad\qquad,
f_{d-3},f_{d-2}, f_{d-1})\\
&=&  
\textrm{(\,\#vertices,\,\#edges,\,\#2-faces,\ \ldots\ ,
\,\#subridges,\,\#ridges,\#facets)}
\end{eqnarray*}
is a long sequence of large numbers, which we may graph just like
a continuous function, and ask for its ``shape.''
Indeed, we might look at a shape function
$\varphi:[0,1]\rightarrow\R$ that is defined by
$\varphi (x) := f_{x(d-1)}$;
this is defined for any $x=\frac k{d-1}$ that is a multiple of~$\frac1{d-1}$,
and these values are rather dense if $d$ is large.
We might interpolate if we want. But what 
types of $f$-vector shape functions $\varphi$ do we get that way?

\begin{figure}[ht]
\begin{center}
\input{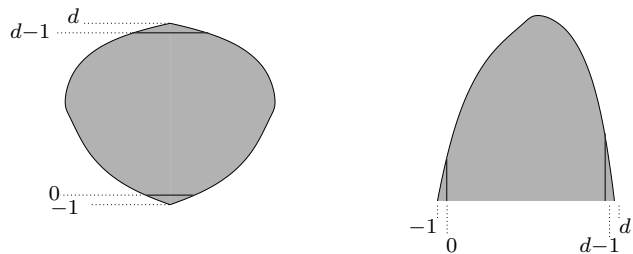}%
\end{center}
\caption{A rough, ``naive'' picture of the shape of the face lattice,
and the $f$-vector, for a high-dimensional polytope}
\label{figure:fshape}
\end{figure}

Figure~\ref{figure:fshape} shows two ``naive'' views, of
the shape of an $f$-vector, and --- equivalently --- of the shape of a
typical face lattice (displayed as a \emph{Hasse diagram}, so the
sizes of rank levels are the $f_i$-values).

A very simple observation is that each vertex
of a $d$-polytope has degree at least $d$,
so double counting yields $f_1\ge\frac{d}{2} f_0>f_0$;
dually, we have $f_{d-2}\ge\frac{d}{2} f_{d-1}> f_{d-1}$.
So in the first step, the $f$-sequence increases, in the last step it
decreases. Does this mean that the $f$-vector 
``first goes up, then comes down,'' that it is \emph{unimodal},
with no ``dip'' in the middle?

\section{Unimodality conjectures}

Unimodality conjectures and theorems abound in combinatorics 
\cite{Stanley-unimodal} \cite{Brenti-unimodal}:
for binomial coefficients, Stirling numbers and their generalizations,
matroids and geometric lattices, etc.\ \ldots
The basic unimodality conjecture for convex polytopes was posed 
at least twice, by Theodore Motzkin in the late fifties, and by 
Dominic Welsh in 1972 (see \cite{Bjo1}). 
Apparently it was disproved dramatically by Ludwig Danzer,
already in the early sixties (presented in a lecture in Graz in 1964,
according to J\"urgen Eckhoff), 
but this is ``lost mathematics,'' no published account exits.

\begin{conjecture}\label{conj:unimodality}
  The $f$-vectors of convex polytopes are unimodal, that is, 
  for each $d$-polytope $P$ there is an $\ell=\ell(P)$ such that
\[
f_0\le f_1\le\ \cdots\ \le f_{\ell}\ge\ \cdots\ \ge f_{d-2}\ge f_{d-1}.
\]
\end{conjecture}

The main point of this lecture will be to 
see that this is dead wrong, even for simplicial polytopes.
Moreover, we want to see this ``asymptotically,''
without substantial amounts of computation, 
without having to list explicit $f$-vectors.

This asymptotic view is also motivated by the fact that the conjecture
fails only in high dimensions.  For example, for simplicial polytopes,
it is true up to $d=19$, and fails beyond this dimension.  For
general polytopes, we will see a counterexample for $d=8$, but none are
known for a smaller dimension.
The conjecture holds in full for $d\le4$ 
(\exercise{ex:unimodal4}),
and also for $d=5$, according to Werner \cite{Werner:f5}.

Since the conjecture is so badly wrong, it might pay off
to explicitly state what remains from it:

\begin{conjecture}[Björner \cite{Bjo1} \cite{Bjo6}]\label{conj:unimodality2}
  The $f$-vectors of convex polytopes increase on the first quarter,
  and they decrease on the last quarter:
\[f_0<f_1<\ \cdots\ <f_{\lceil\frac{d-1}{4}\rceil},\qquad\qquad
  f_{\lfloor\frac{3(d-1)}{4}\rfloor}>\ \cdots\ >f_{d-2} >f_{d-1}.
\]
\end{conjecture}

This is trivially true for $d\le5$.
It also is true for \emph{simplicial} $d$-polytopes
(the $f$-vectors of simplicial polytopes indeed increase up to the
middle, and they decrease in the last quarter), 
but the available proof for this depends on
the necessity part of the $g$-theorem, so it is
quite non-trivial; see \cite{Bjo6}. 

To demonstrate our ignorance on such basic $f$-vector shape
matters, here is a suspiciously innocuous conjecture.
Apparently no one has an idea for a proof, up to now.

\begin{conjecture}[B\'ar\'any]%
\label{conj:unimodality3} 
For any $d$-polytope, $f_k\ge \min\{ f_0, f_{d-1}\}$.
\end{conjecture}

B\'ar\'any's conjecture holds for $d\le6$ \cite{Werner:f5}.
However, not even 
\[
f_k\ \ \ge\ \ \tfrac{1}{10000}\min\{ f_0, f_{d-1}\}
\]
is proven for large dimensions~$d$\,! 
We know so little~\ldots

\section{Basic examples}

Let's compute the $f$-vector shapes for the
most basic high-dimensional polytopes that
we can come up with.
For rough estimates, we use a very crude version
of Stirling's formula, 
\[
n!\ \sim\ \left(\frac ne\right)^n.
\]

\begin{figure}[ht]
\begin{center}
  \psfrag{1}{$1$}
  \psfrag{2}{$2$}
  \psfrag{3}{$3$}
  \psfrag{4}{$4$}
  \psfrag{5}{$5$}
  \psfrag{6}{$6$}
  \psfrag{10}{$10$}
  \psfrag{15}{$15$}
  \psfrag{20}{$20$}
  \psfrag{25}{$25$}
  \psfrag{k}{$k$}
  \psfrag{simplex}{$f_k/10^7$}
  \psfrag{cross}{$f_k/10^{12}$}
  \psfrag{cyclic}{$f_k/10^{18}$}
  \includegraphics[height=55mm]{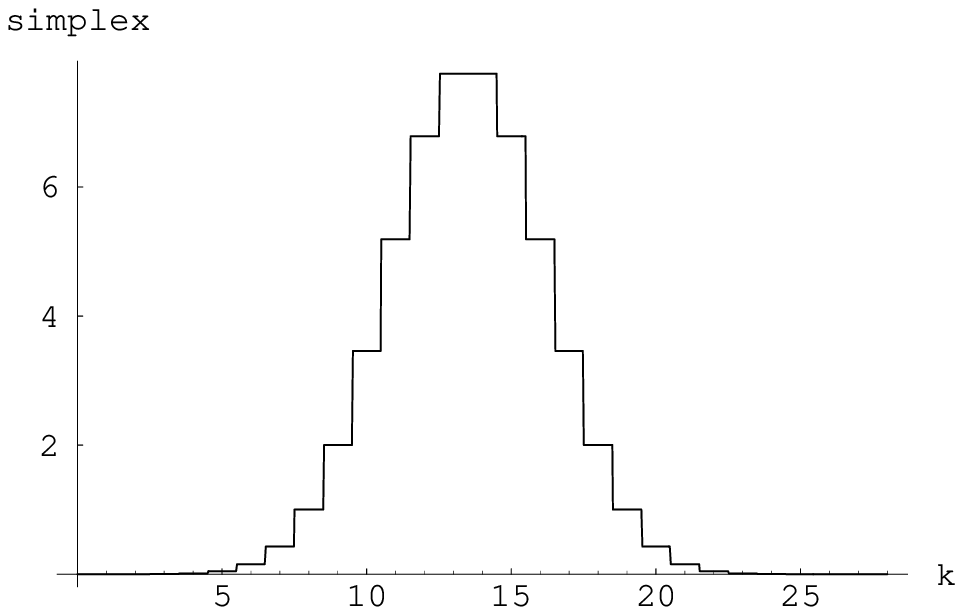} 
  \includegraphics[height=55mm]{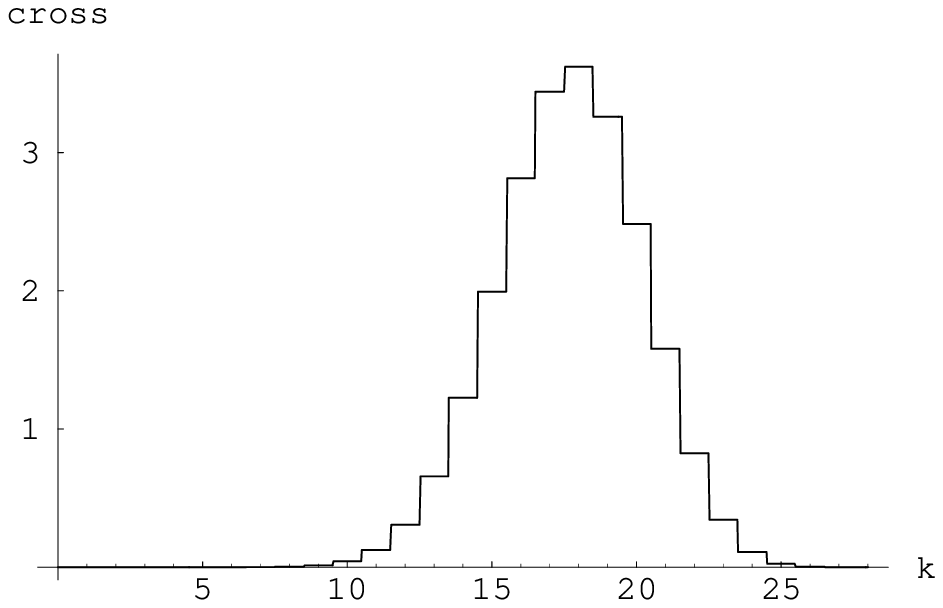} 
  \includegraphics[height=55mm]{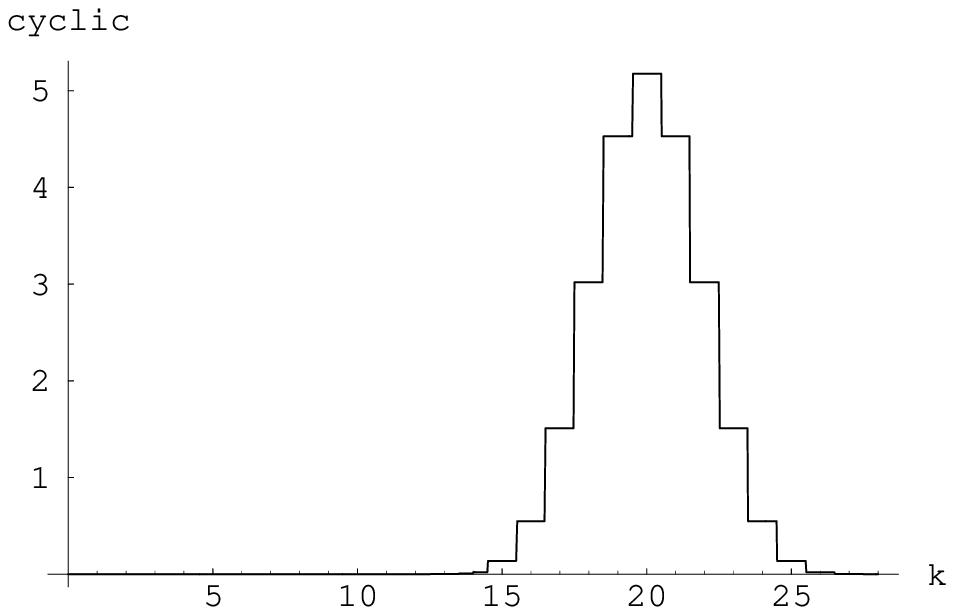} 
\end{center}
\caption{The $f$-vector shapes of the $28$-dimensional simplex $\Delta_{28}$,
  the cross polytope $C_{28}^*$, and a cyclic polytope with 80 vertices
  $C_{28}(80)$}
\label{figure:fshape-simplex}%
\label{figure:fshape-crosspoly}%
\label{figure:fshape-cyclic}%
\end{figure}

\begin{example}[The simplex]
  For the $(d-1)$-simplex $\Delta_{d-1}$ we have 
\[
f_{k-1}(\Delta_{d-1})\ =\ \binom dk.
\]
  With logarithms taken with base~$2$, $x:=\frac kd$, and
$\varphi(x)=f_{xd-1}$, we get 
\[
\log\varphi(x)\ =\ \log\binom{d}{xd}\ \sim\ -x\log x- (1-x)\log (1-x).
\]
  A little bit of analysis shows from this 
  that the $f$-vector is symmetric, with a sharp peak 
  in the middle (at $x=\frac12$), of width $\sim\frac1{\sqrt d}$.
  Figure~\ref{figure:fshape-simplex} displays a realistic example.

  Of course, this is a well-known property of binomial coefficients,
  and the strong limit theorems of probability theory depend on it.  
(In this context $\varphi(x)$ is known as the ``entropy function.'') 

\end{example}

\begin{example}[Cross polytopes]
For the $d$-dimensional cross polytope 
$C_d^*=\conv\{\pm e_1, \ldots, \pm e_d\}$
we have 
\[
f_k(C_d^*)\ =\ \binom{d}{k+1}2^{k+1}.
\]
Again, approximating crudely and taking logarithms base~$2$,
we  get
\[
\log \varphi(x)\ \sim\ -x\log x- (1-x)\log (1-x)+x.
\]
The derivative
\[
\frac{\textrm{d}}{\textrm{d}x}
\varphi(x)\ \sim\ 
 -\log x - \frac{1}{\ln 2}+\log (1-x)+\frac{1}{\ln 2}+1
\]
vanishes at $x=\frac23$: That's where $\log\varphi(x)$ has its maximum,
and where $\varphi(x)$ has a sharp peak 
(compare Figure~\ref{figure:fshape-crosspoly}).

Thus the $f$-vector of a $d$-dimensional cross polytope,
for large $d$, has a sharp peak at $k=\frac23d$.
By duality, this means that the $f$-vector of the $d$-cube
peaks at $k=\frac13d$, for large~$d$.
\end{example}

\begin{example}[Cyclic polytopes]\label{example:cyclic}
Let's look at cyclic polytopes $C_d(n)$ with many vertices,
$n\gg d$. For simplicity, we assume that the dimension $d$ is {even}.

A curve in $\R^d$ has degree $d$ if no $d+1$ points on the curve
lie on a hyperplane.
The convex hull of any $n>d$ points on such a curve
is a cyclic polytope~$C_d(n)$.
\emph{Gale's evenness criterion} \cite{Gale1} gives a combinatorial 
description for the facets, 
which is easy to visualize (see Figure~\ref{figure:Gale}):
Any $d$ points on a degree $d$ curve span a hyperplane~$H$.
If the $d$ points are supposed to span a facet of the polytope,
then all the other $n-d$ points must lie on the same side of~$H$.
Since the curve crosses $H$ \emph{only} in these $d$ points,
this means that the $d$ points split into $\frac d2$ adjacent pairs.
So, if we number the points $1,2,\dots,n$ along the curve, then
the facets of their convex hull (the cyclic polytope)
are given by $\frac d2$ pairs $i,i+1\bmod n$.
The ($k-1$)-faces are given by the $k$-subsets of such a $d$-set:
For $k\le\frac d2$ \emph{any} such subset will do
(the cyclic polytopes are \emph{neighborly}), while
for $k>\frac d2$ the faces consist of $k-\frac d2$ pairs, and 
$d-k$ singletons. Thus the ($k-1$)-faces
may be obtained by choosing $\frac d2$ vertices $i_j$ arbitrarily, 
and also taking $i_j+1$ for $k-\frac d2$ of these 
(see Figure~\ref{figure:Gale2}). 
Thus, with a bit of an over-count, we get
\[
f_{k-1}(C_d(n))\ 
\begin{cases}
\ \ =\    \binom nk             &\textrm{ for } k\le\frac d2,\\
\ \ \sim\ \binom n{\frac d2}\binom{\frac d2}{k-\frac d2}
                                             &\textrm{ for } k > \frac d2.
\end{cases}
\]

\begin{figure}[ht]
\begin{center}
\input{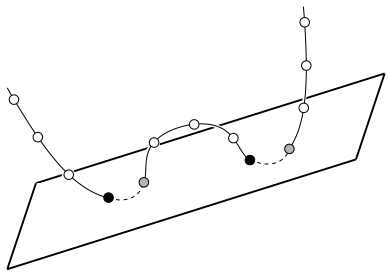}
\end{center}
\vskip-3mm
\caption{Sketch for Gale's evenness criterion.}
\label{figure:Gale}
\end{figure}

\begin{figure}[ht]
\begin{center}
\input{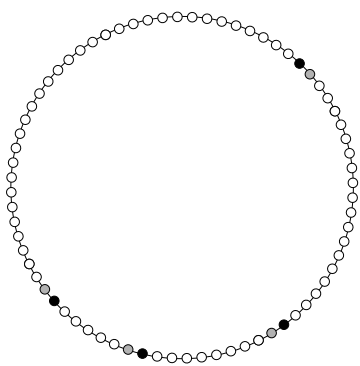}
\end{center}
\caption{An estimate for the number of facets of 
$C_d(n)$, for $n\gg d$, with $d$ even:
There are $\binom n{d/2}$ choices for the black points;
with high probability, they are non-adjacent; the
$\frac d2$ pairs can be completed by taking the gray points.}
\label{figure:Gale2}
\end{figure}

Clearly this peaks at $x=\frac34$:
We get the larger entries in the case $k>\frac d2$, and then
the maximum is achieved when $\binom{d/2}{k-d/2}$
is maximal, that is, for
$k=\frac34d$. Figure~\ref{figure:fshape-cyclic}
gives a realistic impression of the $f$-vector shape of a cyclic polytope.
\end{example}

An explicit,
exact formula for $f_{k-1}(C_d(n))$ is available (\exercise{ex:cyclic}), but
this doesn't answer all the questions. In particular,
is it really true that the $f$-vector is unimodal?
As far as I know, the Unimodality Conjecture~\ref{conj:unimodality}
has not been established in full for the cyclic polytopes.
It does hold for small $n>d$, and certainly also if
$n\gg d$ is sufficiently large compared to~$d$ 
(with the $f$-vector peak at~$k=\lfloor\frac{3(d-1)}4\rfloor$), but in
an intermediate range for $n$ a challenge remains~\ldots

\section{Global constructions}

We have seen classes of simplicial $d$-polytopes whose 
normalized $f$-vector functions $\varphi(x)=f_{x(d-1)}$
peak at $x=\frac12$, at $x=\frac23$, or at~$x=\frac34$.
By dualization we get simple polytopes with peaks at 
$x=\frac13$, and at~$x=\frac14$. 
The ``global constructions'' of products and joins now
yield examples with peaks in the whole range between $x =\frac14$
and $x=\frac34$.
(The product construction is elementary, well-known, and well-understood,
but a review perhaps can't harm, also in view of our needs for
Lecture~5.
Joins are similarly elementary and well-understood,
but perhaps not that well-known.)

\begin{example}[Products]\label{ex:products}
Let $P$ and $Q$ be polytopes of dimensions $d$ and $e$.
Then the product 
\[
P\times Q\ :=\ \{(x,y):x\in P,\, y\in Q\}
\]
is a polytope of dimension $\dim (P\times Q)=\dim P+\dim Q=d+e$.

The nonempty faces of $P\times Q$ 
are the products of nonempty faces of~$P$ and non-empty faces of~$Q$:
In particular, 
the vertices of~$P\times Q$ are of the form ``vertex times vertex,''
the edges are of the form ``edge times vertex'' or ``vertex times edge,''
and the facets are ``$P$ times facet~of~$Q$'' or ``facet of~$P$ times~$Q$.''
With the convention $f_d(P)=f_e(Q)=1$ this yields the formula
\begin{equation}\label{eq:f-product}
f_m(P\times Q)\ =\ \sum_{{k+\ell=m}\atop{k,\ell\ge 0}}  f_k (P)\, f_{\ell}(Q)
\end{equation}
for $m\ge0$.
\end{example}

The product construction is dual to the ``free sum'' construction,
$P\oplus Q$: For this let $x_0\in P\subset\R^d$ and $y_0\in Q\subset\R^e$
be interior points, and take the convex hull
\[
P\oplus Q\ :=\ \conv\big(P\times\{y_0\}\ \cup\ \{x_0\}\times Q\big).
\]
The \emph{proper} faces of~$P\oplus Q$ (that is, faces other than the polytope
itself) arise as joins of proper faces of~$P$ and of~$Q$.

\begin{figure}[ht]
  \begin{center}
  \psfrag{P}{$P$}
  \psfrag{Q}{$\!Q$}
  \includegraphics[width=50mm]{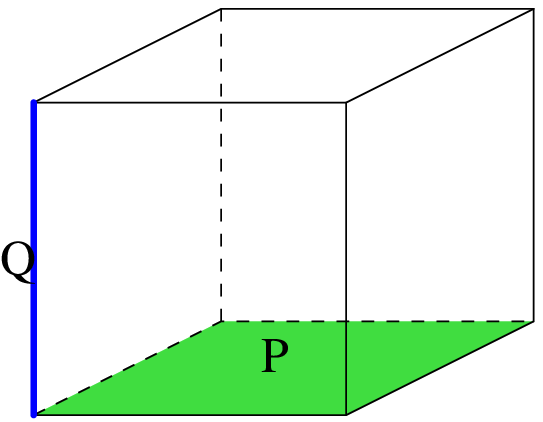}\qquad
  \includegraphics[width=50mm]{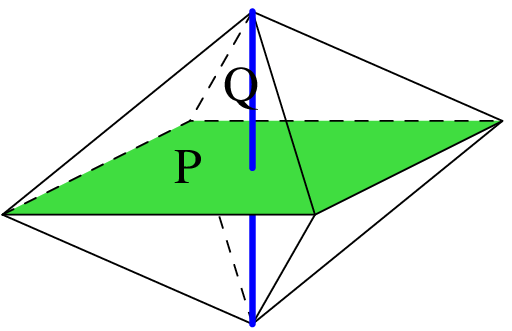}
  \end{center}
  \caption{Product and free sum, for $P=I^2$ and $Q=I$}
  \label{figure:product+free_sum}
\end{figure}

The product and the free sum construction are illustrated
in Figure~\ref{figure:product+free_sum}.

Since joins come up as faces of free sums, let's briefly talk about joins.

\begin{example}[Joins]
Let again $P$ and $Q$ be polytopes of dimensions $d$ and $e$.
Then the \emph{join} $P* Q$ is obtained by positioning
$P$ and $Q$ into skew affine subspaces, and taking the convex hull.
Thus the join is a polytope of dimension
$\dim (P* Q)=\dim P+\dim Q+1=d+e+1$.
  
The faces of~$P*Q$ are the joins of faces of~$P$ and faces of~$Q$:
This refers to \emph{all} faces, including the empty face and the
polytope itself. The corresponding formula, with 
$f_{-1}(P)=f_{-1}(Q)=1$, is
\begin{equation}\label{eq:f-join}
f_m(P*Q) \ =\  \sum_{{k+\ell=m-1}\atop{k,\ell\ge -1}} f_k (P)\, f_{\ell} (Q),
\end{equation}
valid for all $m$, that is, for $-1\le m\le d+e+1$.
\end{example}

\begin{figure}[h]
  \begin{center}
  \psfrag{P}{$P$}
  \psfrag{Q}{\,$Q$}
  \includegraphics[height=50mm]{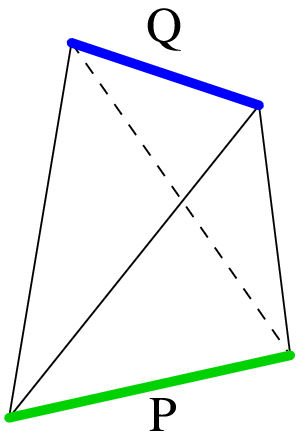}\qquad
  \includegraphics[height=50mm]{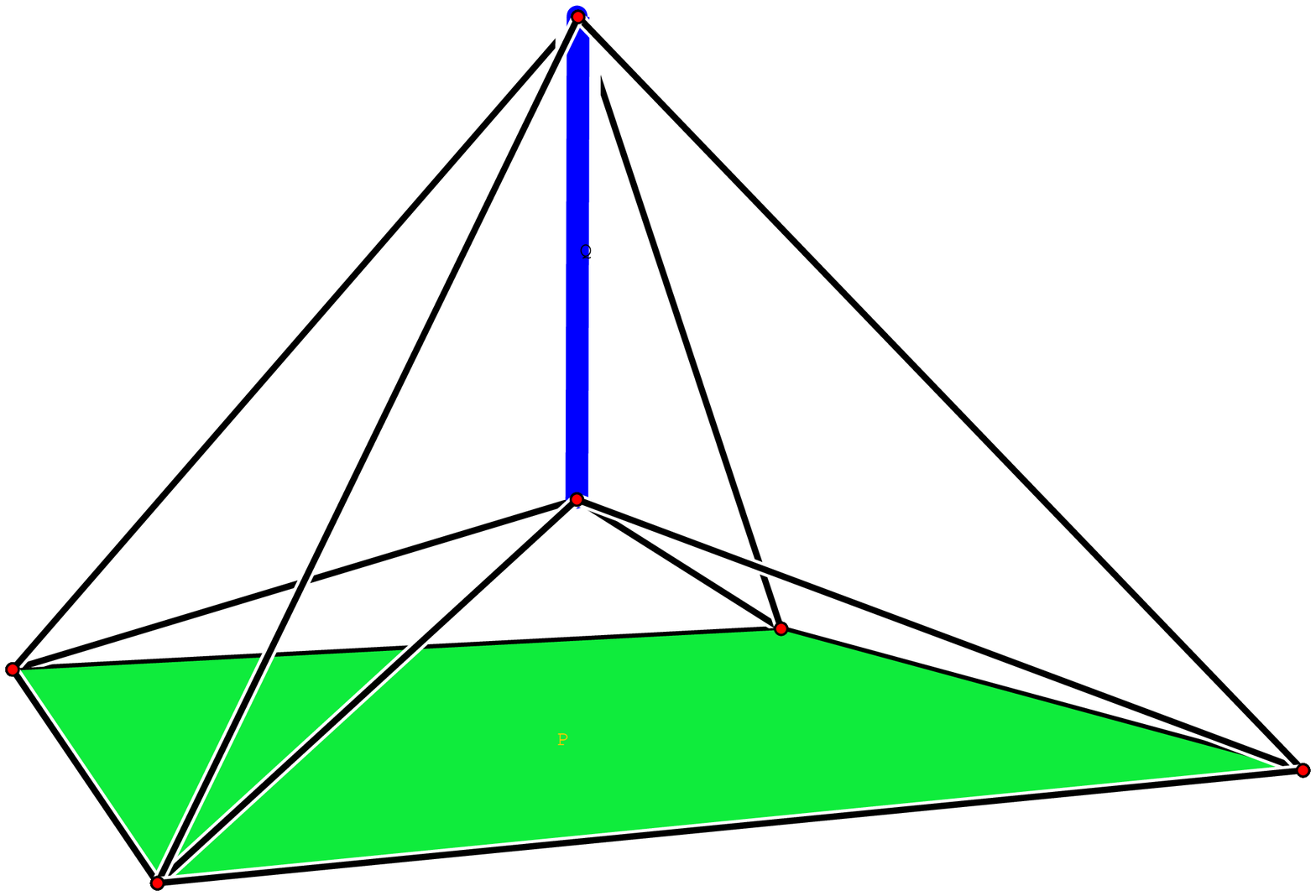}
\end{center}
\caption{Joins $I*I$ and $I^2*I$, of an edge with an edge, resp.\ of
  a square with an edge}
  \label{figure:joins}
\end{figure}

Joins are illustrated in Figure~\ref{figure:joins}.
The dual construction to taking joins is the join construction again.

Product and join are two distinct constructions, and they
do yield different polytopes, of different dimensions (by~$1$).
However, in a birds' eye view, asymptotically, they do behave
quite similarly, and indeed, their effects on $f$-vector shapes
are almost the same. Namely, the formulas (\ref{eq:f-join})
and (\ref{eq:f-product}) describe finite convolutions,
and the only difference is whether the entry $f_{-1}=1$
is counted. For large dimensions, and large $f$-vectors,
this does not make much of a difference, and in both cases we
get a convolution of $f$-vector shapes.
Thus, in particular, if the $f$-vectors of~$P$ and of~$Q$ have
sharp peaks, then the product or join will have a peak as well:
\begin{center}
(peak at $x$) $*$ (peak at $y$)\ \ $\longrightarrow$\ \ 
(peak at $\tfrac d{d+e}x + \tfrac e{d+e}y$).
\end{center}
In particular, for $d=e$ this yields
\begin{center}
(peak at $x$) $*$ (peak at $y$)\ \ $\longrightarrow$\ \ 
(peak at $\tfrac{x+y}{2}$).
\end{center}
To see this, just compute that if the peak 
(or, just the largest $f$-vector entry) for $P_1$ is 
at $x=\frac kd$ and for $P_2$ at $y=\frac\ell e$,
then the peak for $P_1*P_2$ will be at 
\[
\tfrac{k+\ell}{d+e}
\ =\ \tfrac kd\tfrac d{d+e} + \tfrac\ell e\tfrac e{d+e}
\ =\ x\tfrac  d{d+e} + y\tfrac e{d+e}.
\]
This also yields a convolution formula for the $f$-vector shape
of~$P_1\times P_2$ or $P_1*P_2$, for large dimensions:
\[
\varphi(x)\ =\ 
\int_0^1
\varphi_1\big(t\tfrac  d{d+e}\big)\,\varphi_2\big((1-t)\tfrac e{d+e}\big)\,
\mathrm{d}t
\]
Thus, by just taking products of sums of suitable cyclic polytopes and
their duals, we do get polytopes with $f$-vector peaks in 
the whole range between $\frac{1}{4}$ and~$\frac{3}{4}$.

\section{Local constructions}

Perhaps the simplest local operation that can be
applied to a polytope is to ``stack a pyramid onto a simplicial facet.''
To perform such a \emph{stacking} operation geometrically, the
new vertex of course has to be chosen carefully 
(\emph{beyond} the simplicial facet, and 
\emph{beneath} all other facets, in 
Grünbaum's terminology \cite[Sect.~5.2]{Gr1-2}),
but the combinatorial description is easy enough.
In particular, we get the following $f$-vector equation:
\[
f_k(\stack P)\ =\ f_k(P)+f_k(\Delta_d)-f_k(\Delta_{d-1}).
\]
This is valid for $k<d-1$, the rest is ``boundary effects''
that we may safely ignore.
Furthermore, the usual binomial recursion yields
$f_k(\Delta_d)-f_k(\Delta_{d-1})=f_{k-1}(\Delta_{d-1})$, 
and we get
\[
f_k(\stack P)\ =\ f_k(P)\ +\ f_{k-1}(\Delta_{d-1}).
\]
So, the effect of stacking on the $f$-vector is to add
a bump at $x=\frac12$. The effect may be negligible if
the $f$-vector of~$P$ is large, and has large slopes.
However, 
any stacking operation destroys a simplicial facet and creates $d$ new ones,
so it can be repeated. We write $\stack^NP$ for a polytope
that is obtained from $P$ by $N$ subsequent stacking operations.
Thus we get
\[
f_k(\stack^N P)\ =\ f_k(P)\ +\ N\,f_{k-1}(\Delta_{d-1}),
\]
where we may choose $N\ge0$ freely. Thus we are adding a function
with peak at $\frac12$ to a function whose peak may be, for example,
at~$\frac23$. 

\begin{corollary}[Danzer 1964] 
For large enough $d$ and suitable $N$, the
Unimodality Conjecture \ref{conj:unimodality}
fails for ``$N$-fold stacked crosspolytopes'' $\stack^N C_d^*$.
\end{corollary}

Indeed, Danzer apparently also derived that the $f$-vector
of a simplicial $d$-polytope may have not only one dip (between two
peaks), but arbitrarily many dips and peaks! 

Also, dualization yields that a suitable number of vertex truncations applied
to a high-dimensional cube leads to a simple polytope with a
non-unimodal $f$-vector.

However, cross polytopes are not the most effective
starting points for non-unimodal examples: If we use cyclic
polytopes, then the peak (at~$\frac34$) is further away from the
peak for a simplex that we can ``add'' by stacking (at~$\frac12$).
Moreover, in cyclic polytopes we can control the number of vertices
in fixed dimension as well, and thus make the peak at~$\frac34$
as sharp as we want.

\begin{theorem}[Björner \cite{Bjo1} \cite{Bjo6}, Lee \cite{Lee5} \cite{BiLe2}, Eckhoff \cite{eckhoff:_combin}]
\label{thm:UCfalse-simplicial}
The Unimodality Conjecture \ref{conj:unimodality} holds for
simplicial $d$-polytopes of dimensions $d\le19$, but it
fails for $d\ge 20$.
\end{theorem}

Specifically: 
Stacking $N=259\cdot10^{11}$ 
times onto the cyclic polytope $C_{20}(200)$,
one obtains a polytope 
with a dip $f_{11}>f_{12}<f_{13}$ in the $f$-vector,
\begin{eqnarray*}
f_{11} &=& 5049794068451336750\\[-4pt]
\vee\                         \\[-4pt]
f_{12} &=& 5043828885028647000\\[-4pt]
\wedge\                       \\[-4pt]
f_{13} &=& 5045792044986529500.
\end{eqnarray*}



The proof of the first part of 
Theorem~\ref{thm:UCfalse-simplicial} utilizes the 
$g$-theorem (see Stanley \cite{Sta4} and Björner~\cite{Bjo2}),
which explicitly describes the $f$-vectors of
the simplicial polytopes, plus a substantial amount of
``binomial coefficient combinatorics.'' 
See \cite{Bjo6} for $d\le 16$; the extension to
$d\le19$, due to Eckhoff, unfortunately is still not published.

If we leave the realm of simplicial polytopes, then it
becomes even easier to construct polytopes with a non-unimodal
$f$-vector. Then we can try to add the 
$f$-vectors of two polytopes with peaks at $\frac14$
and at $\frac34$, say a cyclic polytope and its dual.
And indeed, just as we can glue a pyramid onto a simplicial facet,
we can glue any polytope with a simplicial facet onto another one
--- after a projective transformation, if needed \cite[p.~274]{Z35}.
The $f$-vector effect of such a glueing is 
essentially
\[
f(P\# P')\ =\ f(P)+f(P') -f(\Delta_{d-1});
\]
if the $f$-vector components of~$P$ and of~$P'$ are 
large, then the simplex may be neglected, and
we are essentially just ``adding the $f$-vectors.''

We can even do this with cyclic polytopes:
For example, $C_d(n)$ is simplicial;
its dual, $C_d(n)^*$ is simple (without simplicial facets),
but if we cut off (``truncate'') one of the simple vertices, then
a simplicial facet results. Write $C_d(n)'$ for the 
``dual with a vertex cut off.'' 

\begin{corollary}[Eckhoff \cite{eckhoff:_combin}]
The Unimodality Conjecture \ref{conj:unimodality} fails for
$d$-polytopes of dimensions $d\ge8$. In particular, 
\[
f(C_8(25) \# C_8(25)')\ =\ (7149, 28800,46800,46400,46400,46800,28800,7149).
\]
\end{corollary}
This $f$-vector has a nice ``1\% dip'' in the middle!  
We don't know whether the Unimodality
Conjecture~\ref{conj:unimodality}
is true for dimensions $d=6$ or~$7$.

\subsection*{Exercises}
\begin{compactenum}[\thechapter.1.]
\item For $d=3,4,5,\dots$ construct a $d$-polytope with $12$ vertices
  and $13$ facets. How far do you get?
\item\label{ex:unimodal4}
  Show that $f$-vectors of $4$-polytopes are unimodal.
\item\label{ex:cyclic}
  Derive an \emph{exact} formula for $f_{d-1}(C_d(n))$,
  and for $f_k(C_d(n))$,   for even $n$.
\item Compute $f_i(C_8(25))$. How bad is the approximation given
  in Example~\ref{example:cyclic}?
\item Count and describe the $2$-faces of a product of
  a pentagon and a heptagon, $P_5\times P_7$.
\item Compute $f((C_{10})^{10})$, for the product of ten $10$-gons.
  Where is the peak?
\item Estimate/compute $d$ and $N$ such that the ``$N$-fold truncated
  $d$-cube'' has a non-unimodal $f$-vector.
\item If you stack ``too often'' onto $C_{20}(200)$, then unimodality
  is restored. How often?
\end{compactenum} 

\lecture{2-Simple 2-Simplicial 4-Polytopes}

The boundary complex of a $4$-polytope is a $3$-dimensional
geometric structure. So, in contrast to the high-dimensional
polytopes discussed in the previous lecture, we can hope
to approach $4$-polytopes via explicit visualization and
geometric constructions. 
\emph{Schlegel diagrams} are a key tool for this.\footnote{%
These were apparently introduced by Dr.\ Victor Schlegel,
a highschool (Gymnasium) teacher from Waren an der Müritz,
in his paper~\cite{Schl1} from~1883.
The plates for the paper include a Schlegel diagram
(``Zellgewebe'') of a $4$-cube, as well as two quite insufficient 
drawings representing the $24$-cell.
Classical, beautiful drawing may be found in 
Hilbert \& Cohn-Vossen \cite[p.~135]{HilbertCohnVossen}.}
Another one,
which we will also depend on in a key moment of this lecture, is
\emph{dimensional analogy}: To describe a construction of
$4$-polytopes, we phrase a key step as a statement that it is valid
``for all $d\ge3$,'' where the visualization is done for the special
case $d=3$, while the most interesting results are obtained for~$d=4$.

The geometry and combinatorics of polytopes in dimension~$4$
is much more interesting, rich, and difficult than in
$3$ dimensions, because $4$-polytopes aren't constrained between only two
extremes, simple and simplicial.
Some of the most fascinating examples around,
such as Schläfli's 24-cell, are neither simple nor simplicial, but
\emph{$2$-simple $2$-simplicial}.  This property was thought to be
rare until recently: Only a few years ago, exactly $8$ such polytopes
were known.
(Unfortunately, a claim by Shephard from 1967 did not work out:
In \cite[p.~82]{Gr1-2} it had been claimed that 
Shephard could produce infinite families, and
that each $4$-dimensional convex body could be approximated by
$2$-simple $2$-simplicial $4$-polytopes, which would
have established a conjecture by David Walkup. 
Compare \cite[p.~96b]{Gr1-2})

The main goal for this lecture is to describe a simple, explicit,
geometric construction that produces rich infinite families of
{$2$-simple $2$-simplicial} $4$-polytopes.  
The first infinite families, obtained by Eppstein, Kuperberg \& Ziegler
in 2001 \cite{Z80}, relied on rather subtle constructions, via Koebe--Thurston
type edge-tangent realizations of $4$-polytopes (which exist only in
rare cases), and hyperbolic angle measurements.
In contrast to this, the
\emph{deep vertex truncation} construction to be described here is
remarkably simple; it appears in Paffenholz \& Ziegler \cite{Z89},
while special instances (for semi-regular polytopes)
can be traced back to Coxeter's classic \cite[Chap.~VIII]{Cox},
who refers to Cesàro (1887) for the construction of 
the 24-cell by what we here call a ``deep vertex truncation''
of the regular $4$-cube.

\section{Examples}\label{sec:Examples}

Let's start with examples of well-known $4$-polytopes --- and for 
each of those let's look at a Schlegel diagram, and record the
$f$-vector 
\[
(f_0,f_1,f_2,f_3)\ =\ 
\textrm{(\,\# vertices, \# edges, \# $2$-faces (=\,ridges), \# facets\,)}.
\]

A \emph{Schlegel diagram} is a way to visualize a $4$-polytope
in terms of a $3$-dimensional complex.
We can't develop the theory of Schlegel diagrams here
(see \cite[Sect.~3.3]{Gr1-2} and \cite[Lect.~5]{Z35}),
but we can offer two interpretations, both in terms of
dimensional analogy.
\begin{compactitem}[~$\bullet$~]
\item Assume that one face of a $3$-polytope is transparent
  (a ``window''), press your nose to the window,
  and look inside: Then you will see all the other faces of the polytope
  through the window. If you now close one eye (and thus lose
  the spatial impression, or depth view), then you will see how the other
  faces tile the window; you can see how they fit together,
  and thus the whole combinatorial structure of the $3$-polytope
  is projected into a $2$-dimensional window.
  This is the Schlegel diagram of a $3$-polytope.
\item Any $3$-polytope can be projectively deformed in such a way that
  looking at it from a suitable point, you \emph{see} all faces
  except for one single face, which is on the back.
  What you see is a polytopal complex which has the same
  shape as the back face, but this is broken into all the many
  faces that you see on the front side. What you see is the
  $2$-dimensional Schlegel diagram of a $3$-polytope.
\end{compactitem}
The Schlegel diagram of a $4$-polytope, analogously, is 
a $3$-dimensional complex that represents all the faces
of the polytope, except for one facet (the window resp.\ back facet).
The whole combinatorial structure of the polytope may be read
from such a visualization. Thus, for example, one can tell 
whether the polytope is simple, or simplicial, or cubical, etc.

The pictures of Schlegel diagrams as presented in the following are
generated automatically in the \texttt{polymake} system 
by Gawilow \& Joswig \cite{GawrilowJoswig}, with the
\texttt{javaview} back-end by Polthier et~al.~\cite{polthier-javaview04}.
They have three limitations:
They show only a $2$-dimensional projection of an object that
you should see rotating, $3$-dimensionally, on a screen;
they depict only the edges, so in some examples it is hard
to tell/imagine where the faces and facet-boundaries go;
and we don't have color available here.
Nevertheless, I think they are impressive, and you should be
able to ``see'' in them what the (boundary complexes of) some 
$4$-polytopes look like.

\begin{figure}[ht]
\begin{center}
\includegraphics[height=60mm]{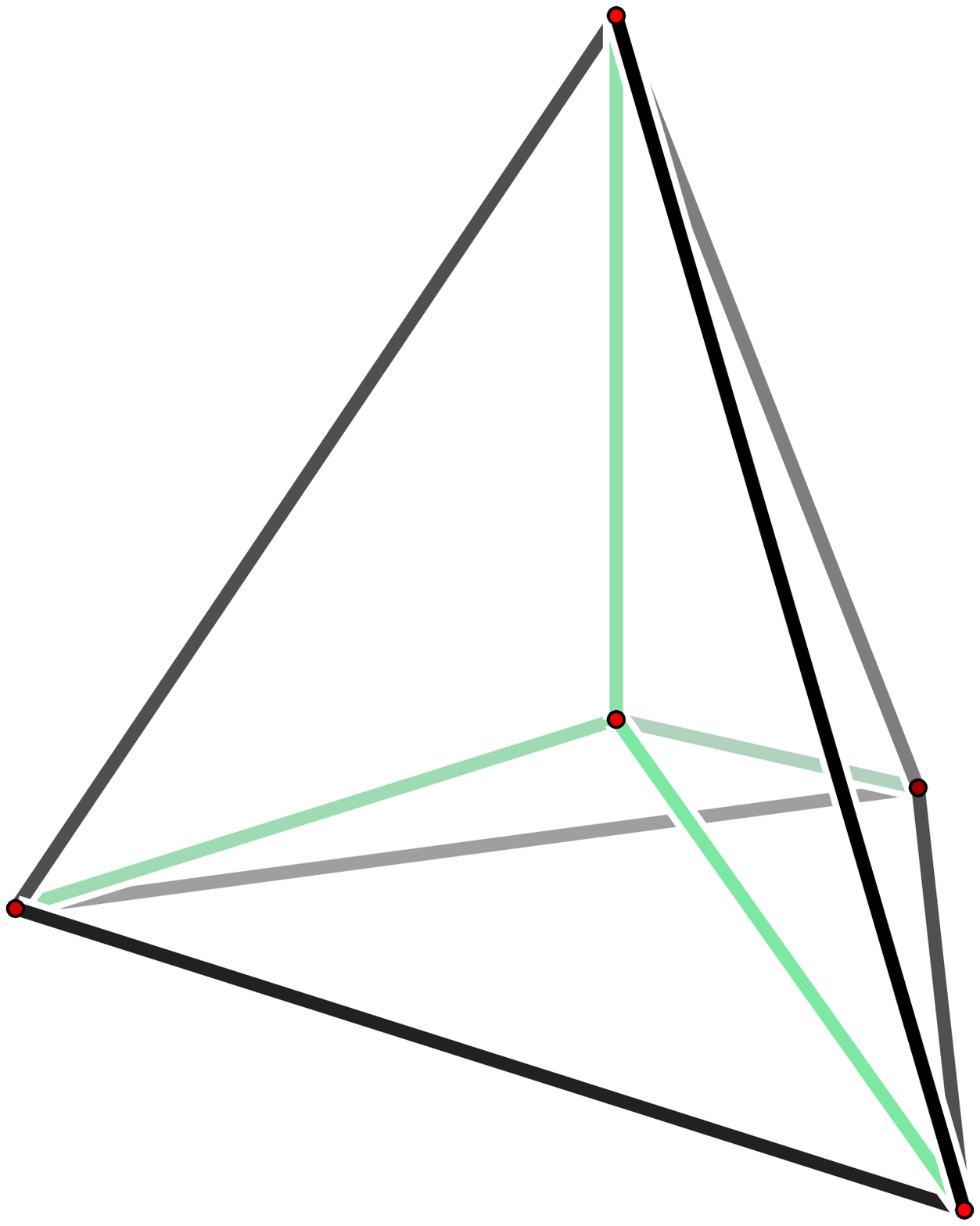}\\[-5mm]
\mbox{\includegraphics[height=60mm]{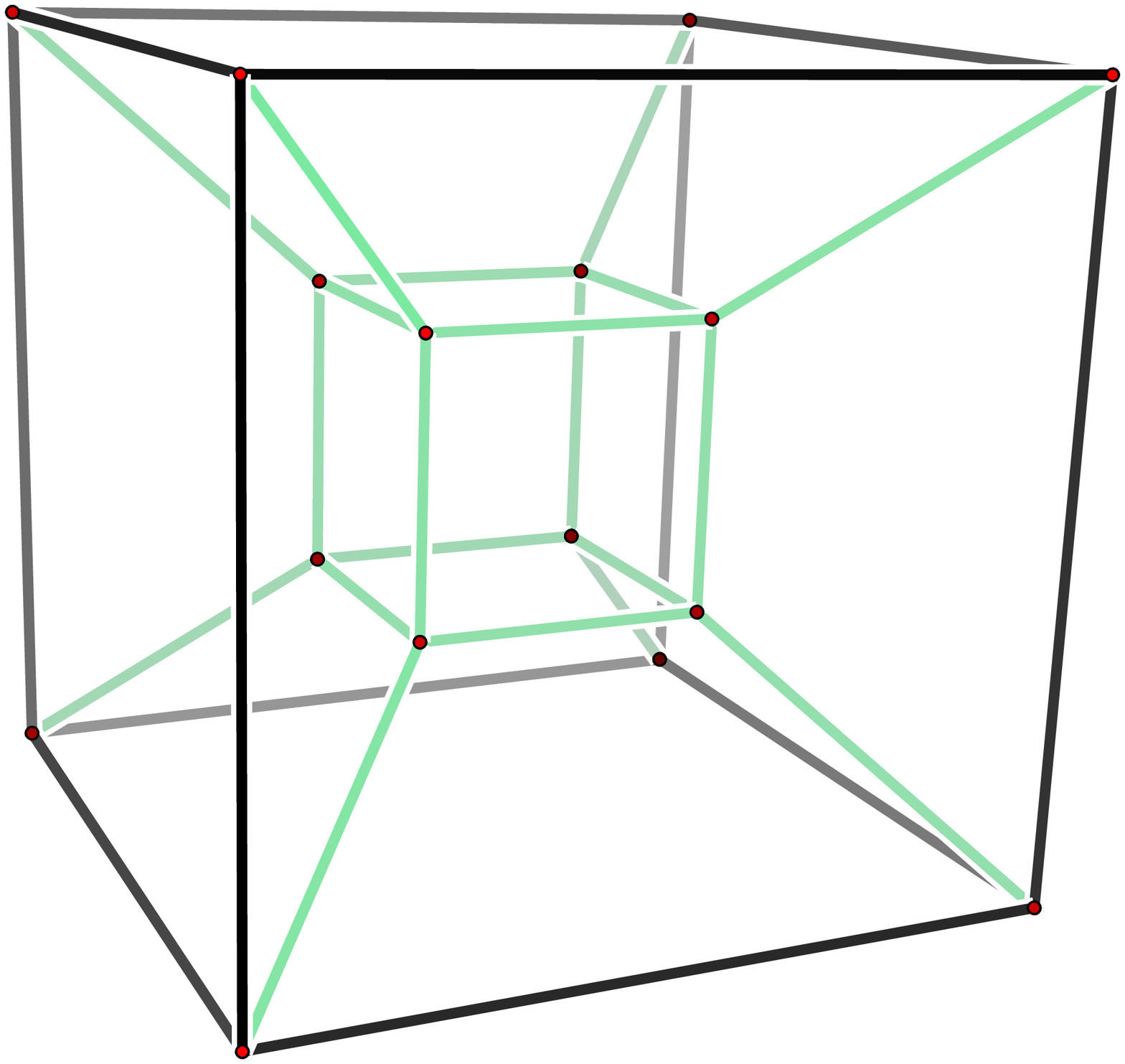}\hspace{-10mm}
\includegraphics[height=70mm]{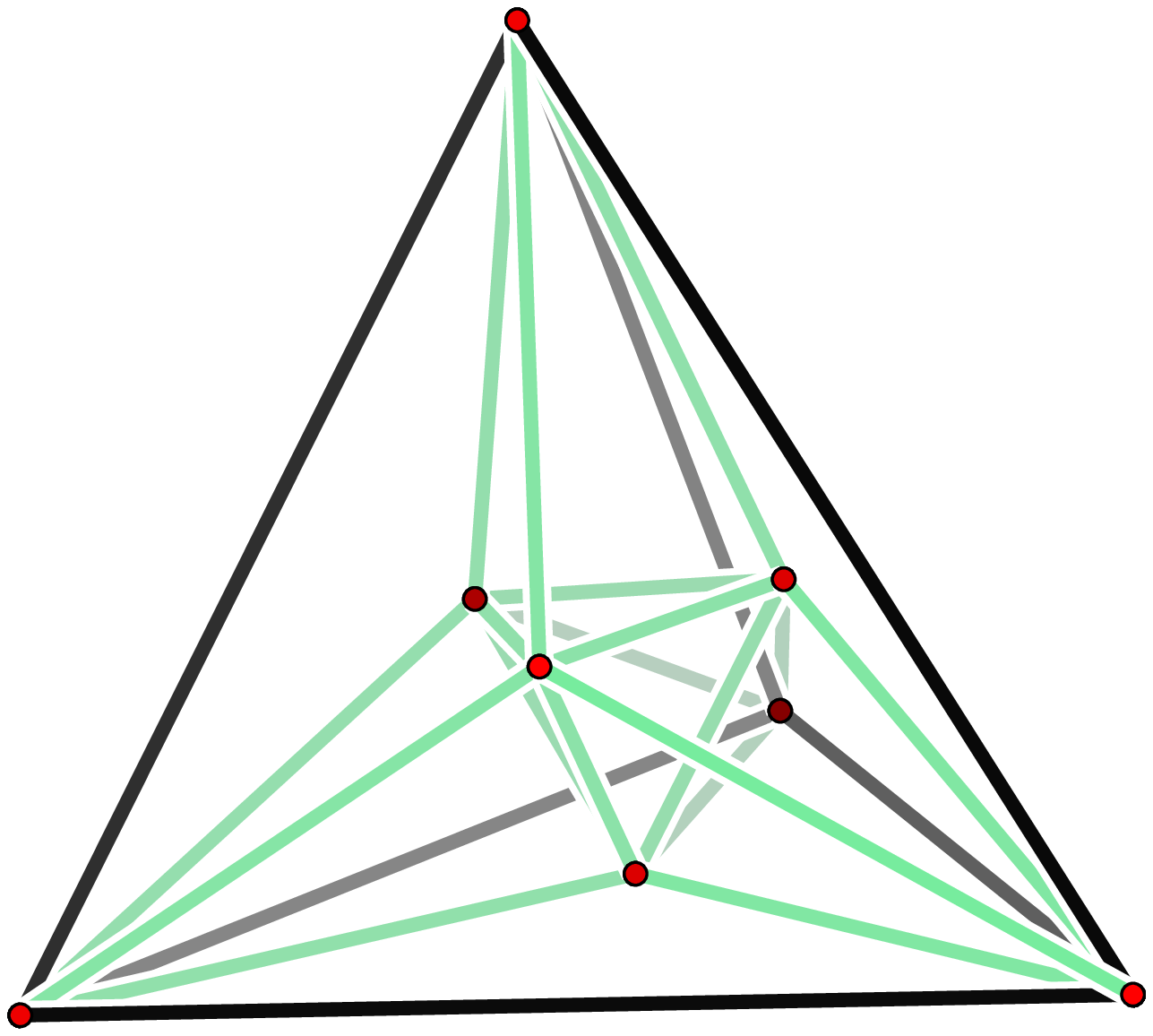}}
\end{center}
\caption{Schlegel diagrams for the $4$-dimensional simplex, cube, and 
cross polytope}
\label{figure:Schlegel1}
\end{figure}

\begin{example}[Simplex, cube, and cross polytope]
  Schlegel diagrams of the $4$-simplex, the $4$-cube and
  the $4$-dimensional cross polytope appear in Figure~\ref{figure:Schlegel1}.
  You should read off the $f$-vectors from this figure:
$f(\Delta_4)=(5,10,10,5)$, 
$f(C_4)=(16,32,24,8)$, and 
$f(C_4^*)=(8,24,32,16)$.

The simplex and cube are simple, so $f_1=2f_0$, while the
simplex and cross polytope are simplical, so $f_2=2f_3$.
\end{example}

\begin{figure}[ht]
\begin{center}
\includegraphics[width=100mm]{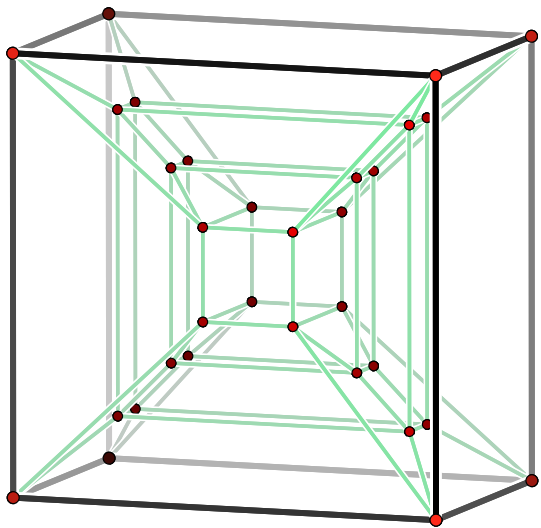}
\end{center}
\vskip-5mm
\caption{A cubical $4$-polytope with the graph of the $5$-cube}
\label{figure:Schlegel2}
\end{figure}

\begin{example}[A cubical $4$-polytope with the graph of a $5$-cube \cite{Z62}]
The construction 
\[
P\ :=\ \conv((2Q\times Q)\cup (Q\times 2Q)), 
\]
for a square such as $Q=[-1,1]^2$, 
yields a $4$-polytope whose Schlegel diagram is displayed in
Figure~\ref{figure:Schlegel2}. This polytope is \emph{cubical}:
All its facets are combinatorially equivalent to the $3$-cube $[-1,1]^3$.

The $f$-vector $(32,80,72,24)$ may be derived from the figure, but indeed
it may also be deduced just from the information that this
is a cubical $4$-polytope with the graph of a $5$-cube.
(The latter yields $f_0$ and $f_1$, the ``cubical'' property
implies $2f_2=6f_3$ by double counting, and then there is the
Euler--Poincar\'e equation \cite[Sect.~8.2]{Z35},
which for $4$-polytopes reads
$f_0-f_1+f_2-f_3=0$.
See also \exercise{ex:compute-f}.)%
\end{example}

\begin{figure}[ht]
\begin{center}
\includegraphics[height=70mm]{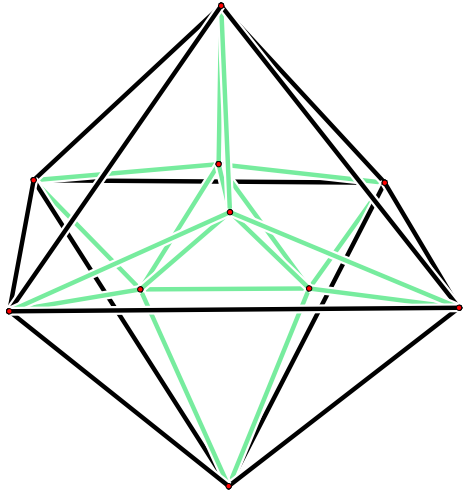}
\end{center}
\caption{A Schlegel diagram of the hypersimplex}
\label{figure:SchlegelHypersimplex}
\end{figure}

\begin{example}[The hypersimplex]
The hypersimplexes form a $2$-parameter family $\Delta_{d-1}(k)$ of remarkable
polytopes; as Robert MacPherson said in his PCMI lectures, they
have by far not received the attention, study, and popularity
that they deserve. They do appear, for example, 
as $K^d_k$ in \cite[p.65]{Gr1-2}, 
as $\Delta^{k,\ell}$ in \cite[Sect.~1.6]{GGL} %
(where apparently the name ``hypersimplex'' appeared first),
in \cite{gelprimefand82:_geomet_grass},
in \cite[p.~207]{GKZ}, and in \cite{LST};
but also elsewhere they appear under disguise, for example,
as the cycle polytopes of uniform matroids (see e.g.~\cite{Groetschel2004b}).

The hypersimplex $\Delta_{d-1}(k)$ may be defined as the 
convex hull of all the $0/1$-vectors of length $d$
that consist of $k$ ones and $d-k$ zeroes.
This is a $(d-1)$-dimensional polytope with $\binom dk$ vertices.
In the special case $k=1$ and $k=d-1$ we obtain simplices.

What we call \emph{the hypersimplex} is a $4$-dimensional
polytope $\Delta_4(2)$ that appears in this family.
It may be defined, lying on a hyperplane in~$\R^5$, as
\[
\big\{x\in[0,1]^5:\sum_{i=1}^5 x_i=2\big\}\ =\ 
\conv\{e_i+e_j:1\le i<j\le5\}, 
\]
or equivalently, after projection to $\R^4$ by ``deleting the last
coordinate,'' as
\[
\big\{x\in[0,1]^4:1\le \sum_{i=1}^4 x_i\le2\big\}\ =\ 
\conv\big(\{e_i:1\le i\le4\}\cup\{e_i+e_j:1\le i<j\le4\}\big).
\]
The first representation is more symmetric: It yields ``by
inspection'' that all $\binom52=10$ vertices of this polytope are
equivalent (under symmetries that permute the coordinates), but that
there are two types of facets, five simplices and five octahedra,
which appear in vertex-disjoint pairs, ``opposite to each other,'' in
parallel hyperplanes.  In particular, all the facets are simplicial,
that is, all the $2$-faces are triangles, so the polytope is
\emph{$2$-simplicial}.  

The second representation has the advantage of
being full-dimensional, and it supplies us with a Schlegel diagram
(using an octahedron facet as a ``window''), as displayed in
Figure~\ref{figure:SchlegelHypersimplex}.  In the figure we may see
that the (ten, equivalent) vertex figures are triangular prisms, so
they are simple; thus in this $4$-polytope, each edge is in exactly
three facets, so the polytope is \emph{$2$-simple}.  So we have seen our
first example (other than the $4$-simplex) of a $2$-simple, $2$-simplicial
$4$-polytope.

From the data given it is easy to compute the $f$-vector of
the hypersimplex: It is $f=(10,30,30,10)$.
\end{example}

\section{\boldmath$2$-simple $2$-simplicial $4$-polytopes}

\begin{definition}\label{def:2s2s}
  A $4$-polytope $P\subseteq\R^4$
  is \emph{$2$-simple $2$-simplicial}
  (``\emph{$2$s$2$s}'' for short) if 
  all $2$-faces of~$P$, and of~$P^*$, are triangles.
\end{definition}

\begin{figure}[ht]
\begin{center}
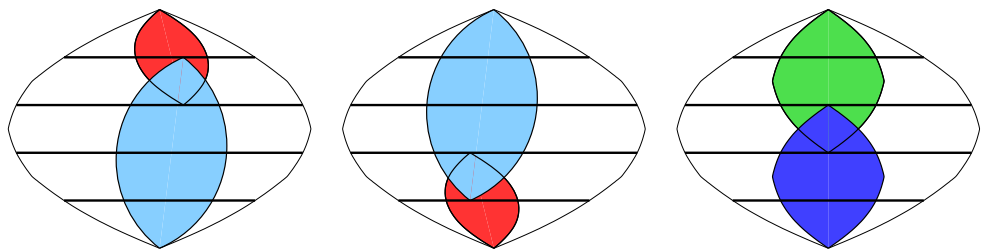
\end{center}
\caption{Simplicial, simple, and 2s2s $4$-polytopes in terms
of their face lattices: The shaded intervals, and all the other
intervals between the same rank levels, must be boolean.}
\label{figure:SimpleSimplicial2S2S}
\end{figure}

The definition given here has the nice feature of being
self-dual: Clearly, $P$ is 2s2s if and only if its dual $P^*$ is~2s2s.
A more explicit version is that a $4$-polytope is 2s2s if and only if 
\begin{compactitem}[~$\bullet$~]
\item every $2$-face has the minimal number $3$~of vertices, and if 
\item every $1$-face (edge) lies in the minimal number $3$~of facets.  
\end{compactitem}
Still equivalently, this is if and only if
\begin{compactitem}[~$\bullet$~]
\item for every $2$-face $G$ the lower interval $[\emptyset,G]$ in the
  face lattice of~$P$ is boolean, and if
\item for every $1$-face $e$ the upper interval $[e,P]$ in the
  face lattice of~$P$ is boolean.
\end{compactitem}
Thus the 2s2s property may be pictured in analogy with the properties
of being simple, or being simplicial.
For this we note that, for example, $P$ is simplicial if
\begin{compactitem}[~$\bullet$~]
\item for every $3$-face $F$ (facet) the lower interval $[\emptyset,F]$ in the
  face lattice of~$P$ is boolean, and if
\item for every $2$-face $R$ (ridge) the upper interval $[R,P]$ in the
  face lattice of~$P$ is boolean.
\end{compactitem}
(The first property just says that the facets should be
simplices; the second property is automatically satisfied: Every ridge
lies in two facets.)
And similarly for simple $4$-polytopes --- see
Figure~\ref{figure:SimpleSimplicial2S2S}.
 
Of course all this suggests generalizations, to ask for 
\emph{$h$-simple $k$-simplicial $d$-polytopes},
apparently introduced by Grünbaum \cite[Sect.~4.5]{Gr1-2}.
For $h+k>d$ these don't exist (other than the $d$-simplex),
but also for small $h$ and $k$ they are hard to construct.
Indeed, are there any $5$-simple $5$-simplicial $d$-polytopes
that are not simplexes? Not a single example is known.
Compare \cite{Z89} for more information.
Here we will restrict ourselves to the $4$-dimensional case
of 2s2s polytopes. Let's note one interesting
property that is specific for the $4$-dimensional case, 
and which also confirms the impression that 
2s2s $4$-polytopes form a ``diagonal'' case.

\begin{lemma}\label{lemma:symmf}
Every $2$s$2$s $4$-polytope has a symmetric $f$-vector: $f_0=f_3$, $f_1=f_2$.
\end{lemma}

\begin{proof}
  If $P$ is $2$-simplicial, then each $2$-face has three edges.
  Thus the number of incidences between $2$-faces and edges,
  denoted $f_{12}$, is 
  $f_{12}=3f_2$. If it is $2$-simple, then each edge lies in three
  $2$-faces, that is, the number of indicences is
  $f_{12}=3f_1$.  Combination of the two
  conditions forces $f_1=f_2$.  With this, Euler's equation
  yields $f_0=f_3$.
\end{proof}

\begin{figure}[ht]
\begin{center}
\includegraphics[height=70mm]{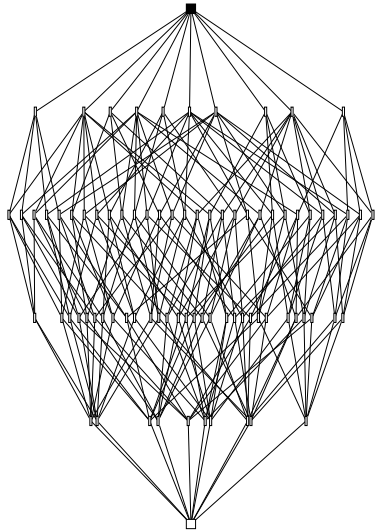}
\end{center}
\caption{The face lattice of the hypersimplex}
\label{figure:facelattice:hypersimplex}
\end{figure}

This proof may be rephrased in terms of the face lattice: For
$4$-polytopes the
2s2s conditions force the two middle rank levels of the face
lattice to form a bipartite cubic graph --- which as any other regular
bipartite graph has to have the same number of vertices on each shore.
You should identify this
bipartite cubical graph in the face lattice of
the hypersimplex, as displayed in 
Figure~\ref{figure:facelattice:hypersimplex}, and thus verify the
2s2s property for this face lattice.
The symmetry of the $f$-vector $(10,30,30,10)$ is explained
by Lemma~\ref{lemma:symmf}; nevertheless, the hypersimplex and its
face lattice are not self-dual: There are two types of
facets, but only one symmetry class of vertices.

The fact that the dual of any 2s2s $4$-polytope is again 2s2s (by
definition), and the symmetry property for the
$f$-vector, might suggest that 2s2s polytopes live
in some sense ``between'' simple and simplicial.
This is not true, as we will see in the next
lecture, when we locate their $f$-vectors
in the cone of all $f$-vectors of $4$-polytopes.
Indeed, the 2s2s polytopes are so interesting because they
form a class of extremal polytopes
in terms of the flag vector: A $4$-polytope is 2s2s if
and only if the valid inequality 
\[2f_{03}\ \ge\ (f_1+f_2)+2(f_0+f_3)\]
holds with equality. (Compare \exercise{ex:CF1}.)

\section{Deep vertex truncation}

The idea for ``deep vertex truncation'' is very easy: 
Cut off all vertices of a polytope --- but don't
just truncate the vertices, but cut them off by
``deep cuts,'' that is, so deeply 
that exactly one point remains from each edge.

All that is said and done about ``deep vertex truncation''
in the following works and makes sense for $d\ge3$.
Nevertheless, the pictures will primarily represent the 
case $d=3$, while the most interesting results appear
for $d=4$.

\begin{definition}[Deep vertex truncation]
Let $P$ be a $d$-polytope, $d\ge2$.
  
A \emph{deep vertex truncation} 
\[
\DVT(P)\ =\ P\ \cap\ \bigcap_{v\in V(P)} H^-_v
\]
of~$P$ is obtained by cutting off all the vertices $v\in V(P)$
of~$P$ (by closed halfspaces $H_v^-$, one for each vertex~$v$) in such
a way that from each edge~$e$ of~$P$, exactly one (relative interior)
point $p_e$ remains.

Equivalently, a deep vertex truncation is obtained 
as the convex hull 
\[
\DVT(P)\ =\ \conv\big\{p_e: e\in E(P)\big\}
\]
of points $p_e$ placed on the edges $e\in E(P)$ of~$P$ in
such a way that for each vertex of~$P$, 
the points $p_e$ chosen on the edges adjacent to~$v$ 
lie on a hyperplane $H_v$.
\end{definition}

It is quite obvious that a deep vertex truncation $\DVT(P)$
can be constructed for each simple polytope~$P$, 
but we will be particularly interested in the case
of simplicial polytopes: For these it is not so clear that
the cutting can be performed so that all constraints
are satisfied simultaneously. 

\begin{lemma}
  Every $3$-polytope has a realization for which deep vertex truncation
  can be performed.
\end{lemma}

\begin{proof}
  Take an edge-tangent Koebe--Andreev--Thurston representation
  (according to Lecture~1). Then $p_e$ can be taken as the tangency
  points, and the cutting hyperplanes $H_v$ are spanned by the vertex
  horizon circles.
\end{proof}

\begin{figure}[h]
\begin{center}
  \includegraphics[height=40mm,bb=66 120 488 371,clip]{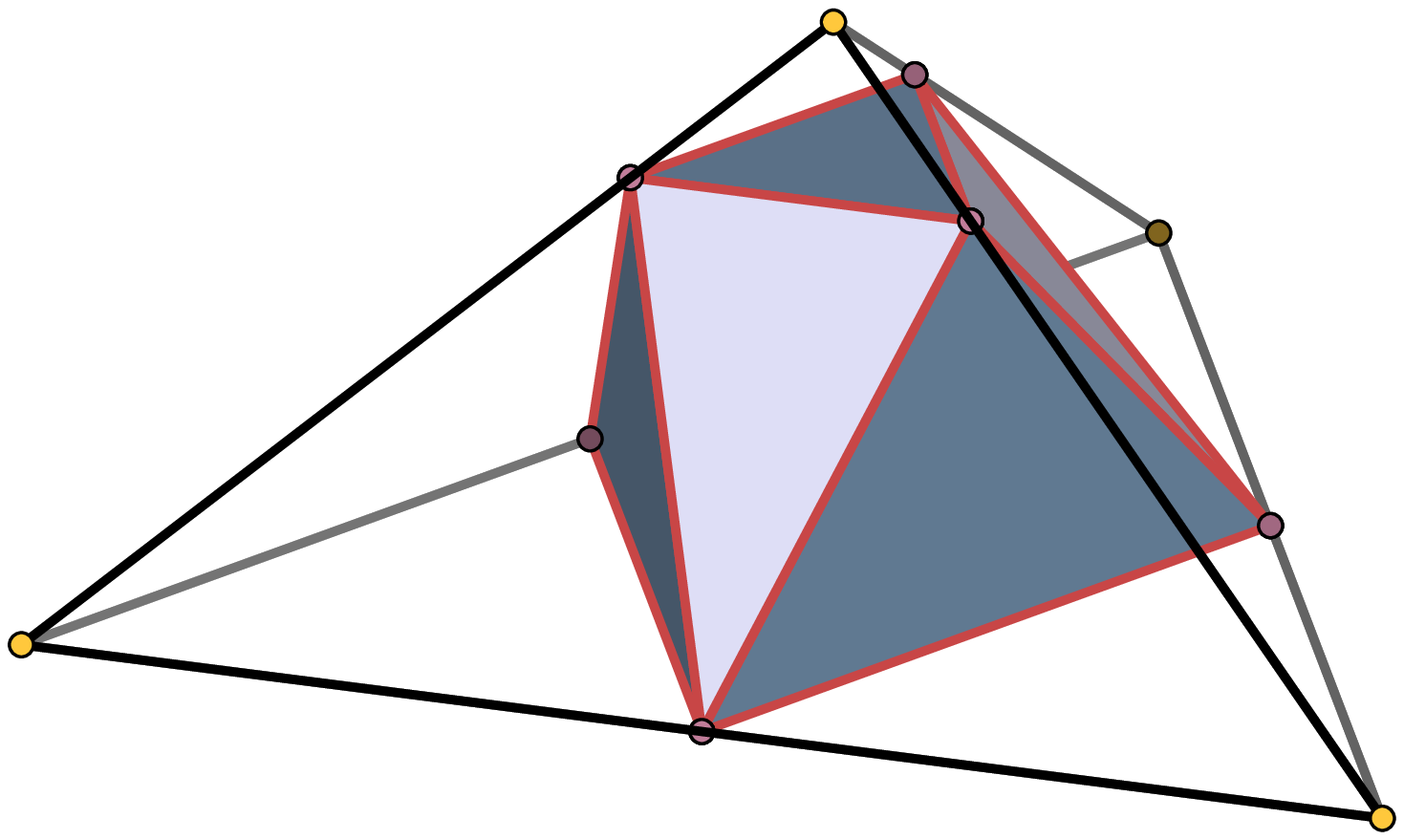}
  \includegraphics[height=40mm,bb=36 111 469 418,clip]{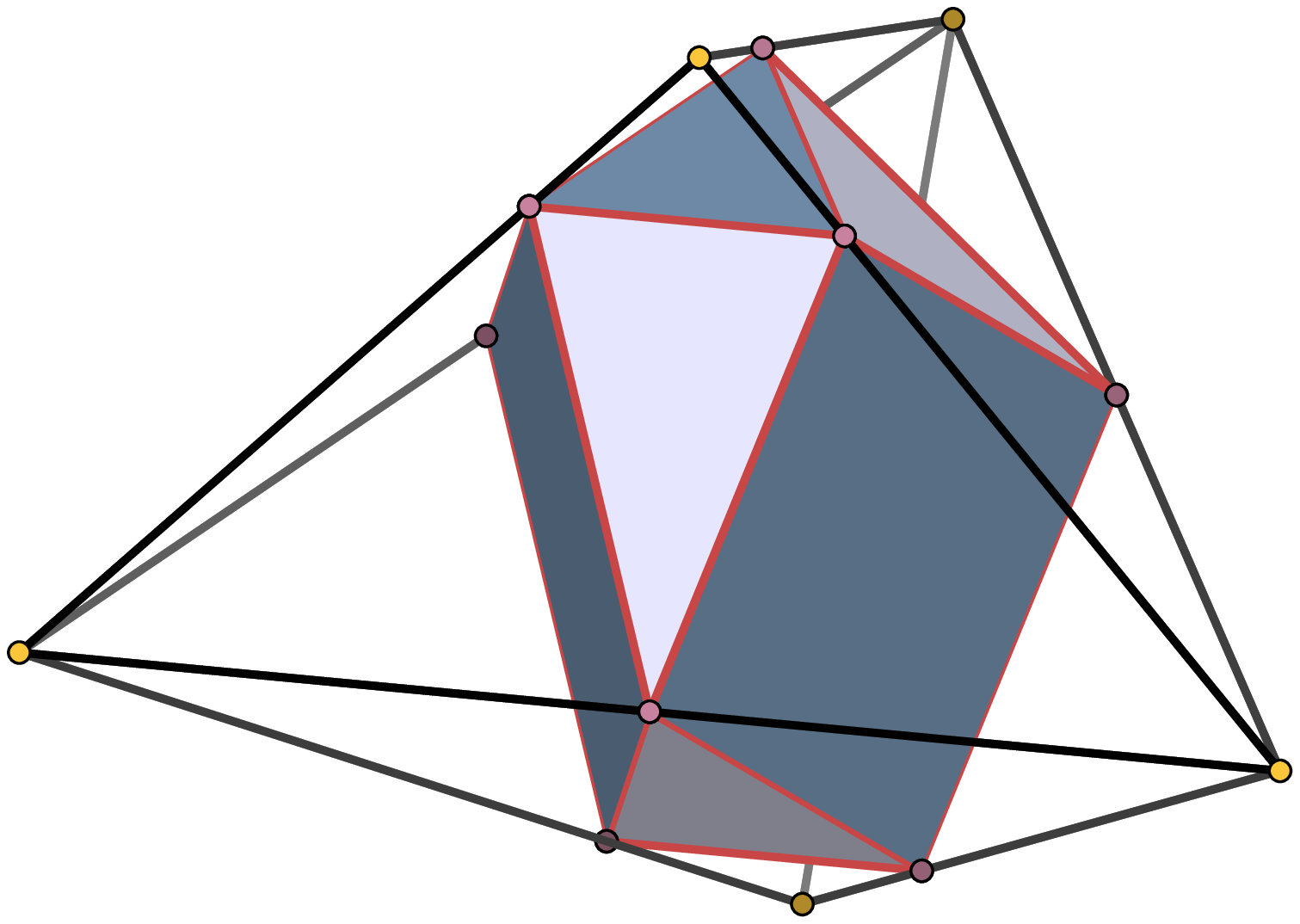}
\end{center}
\caption{Deep vertex truncation of a simplex and of a bipyramid 
yields an octahedron, and a polytope that is ``glued'' from
two octahedra. (Pictures from \cite{Z89})}
\label{figure:DVT1}
\end{figure} 

For $d\ge3$, every deep vertex truncation polytope $\DVT(P)$ 
has two types of facets:
\begin{compactitem}[~$\bullet$~]
  \item deep vertex truncations $\DVT(F)$ of the facets $F$ of~$P$, and
  \item the vertex figures $P\cap H_v=\conv\{p_e: e\ni v\}$ of~$P$.
\end{compactitem}

\begin{proposition}[Paffenholz \& Ziegler \cite{Z89}]\label{prop:DVT2s2s}
  If $P$ is a simplicial $4$-polytope, then any deep vertex truncation
  $\DVT(P)$ is $2$-simple and $2$-simplicial.
\end{proposition}

\begin{proof}
  The two types of facets of~$\DVT(P)$ are the octahedra $\DVT(F)$,
  for the tetrahedron facets $F$ of~$P$, and the
  vertex figures of~$P$, which are simplicial. Thus $\DVT(P)$
  $2$-simplicial.

  Since all edges of~$P$ are reduced to points by deep vertex truncation,
  all the edges of~$\DVT(P)$ are ``new,'' they arise
  by deep vertex truncation from the $2$-faces (that is, the ridges) of~$P$.
  Each such ridge lies in two facets $F_1,F_2$ of~$P$,
  so the edge we are looking at lies in two facets 
  $\DVT(F_1)$ and  $\DVT(F_2)$ of the first type, and in 
  one facet of the second type. Thus each 
  edge of~$\DVT(P)$ lies in exactly three facets, that is,
  $\DVT(P)$ is $2$-simple.
\end{proof}

So we have that $\DVT(P)$ is 2s2s for any simplicial $4$-polytope $P$
\ldots\ if it exists. And that's the problem: In general it is not at
all guaranteed that deep vertex truncation can be performed.  One
would try to realize cyclic $4$-polytopes in such a way that deep
vertex truncations can be performed, but it seems that this is not
possible.  Similarly, if a sum $P_m\oplus P_n$ is realized ``the
obvious way,'' with regular polygons in orthogonal subspaces, then
deep vertex truncation is not possible except for very special cases
(such as $\frac1m+\frac1n\ge\frac12$): It is quite surprising that the
sums of polygons \emph{do} have a realization such that deep vertex
truncation is possible, as proved by Paffenholz~\cite{paffenholz-pc}.
On the other hand, there does not seem to be a single example
of a simplicial polytope for which it has been \emph{proved}
that deep vertex truncation is impossible for all realizations.

\begin{figure}[ht]
\begin{center}
\includegraphics[height=84mm]{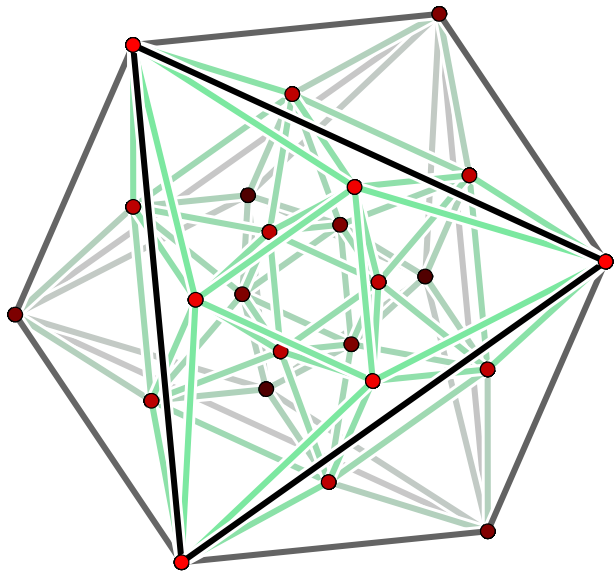}
\end{center}
\vskip-6mm
\caption{The 24-cell}
\label{figure:24cell}
\end{figure}

However, in special cases deep vertex truncation can indeed be performed.
In particular, any \emph{regular} polytope admits a deep vertex
truncation --- just take the edge midpoints for $p_e$.
From this we get the following three examples of 2s2s $4$-polytopes:
\begin{compactitem}[~$\bullet$~]
\item Deep vertex truncation of a simplex, $\DVT(\Delta_4)$,
       yields the hypersimplex.
\item Deep vertex truncation of the $4$-dimensional cross polytope,
   \[ C_4^*\ =\ \conv\{\pm e_i:1\le i\le 4\}\ =\ 
               \{x\in\R^4: |x_1|+|x_2|+|x_3|+|x_4|\le1\},\]
   yields Schläfli's 24-cell (see Figure~\ref{figure:24cell}):
   \begin{eqnarray*}
\hspace{5mm}\DVT(C_4^*)&{=}& 
       \conv\{\pm \tfrac12e_i\pm\tfrac12e_j:1\le i<j\le 4\}\\ 
 &{=}&         \{x\in\R^4: |x_i|\le1\textrm{ for }1\le i\le 4,\ 
                         |x_1|+|x_2|+|x_3|+|x_4|\le1          \}. 
   \end{eqnarray*}
 \item Deep vertex truncation of the regular 600-cell (which has 600
   regular tetrahedra as facets) yields a 2s2s $4$-polytope with 720
   vertices, whose vertex figures are prisms over regular pentagons;
   its facets are 600 octahedra, and 120 regular icosahedra.  
   It seems that this remarkable polytope, with 
   $f$-vector $(720,\,3600,\,3600,\,720)$,
   first occured in the literature in 1994,   
   as the dual of the ``dipyramidal 720-cell'' 
   constructed by G\'evay~\cite{Gevay}.
   See also \exercise{ex:compute-f}.
\end{compactitem}

\section{Constructing \boldmath$\DVT(\Stack(n,4))$}

The stacked polytopes form an infinite family of 
simplicial polytopes which can quite easily be 
realized in such a way that deep vertex truncation 
can be performed.

For this, we denote by $\Stack(n,d):=\stack^n(\Delta_d)$ 
any combinatorial type of a $d$-polytope, $d\ge3$, which is obtained by 
$n$ times stacking a pyramid onto a simplex facet,
starting at a $d$-simplex. This is a 
simplicial $d$-polytope with $d+1+n$ vertices
and $d+1+n(d-1)$ facets; see~\exercise{ex:compute-f2}.
Note that the notation ``$\Stack(n,d)$'' does not 
specify a combinatorial type; many different types may
be obtained by stacking onto different sequences of facets
(cf.~\exercise{ex:exp-stacked}).

\begin{theorem}[Paffenholz \& Ziegler \cite{Z89}]\label{thm:CutStack}
  Any combinatorial type of a stacked $d$-polytope $\Stack(n,d)$ can be
  realized so that it admits a deep vertex truncation.
\end{theorem}

\begin{proof}
  We proceed by induction on~$n$, starting at $n=0$,
  with a $d$-simplex, and a deep vertex truncation that
  takes the convex hull of the edge midpoints.

  Assume now that $\Stack(n,d)$ has been realized as $P\subset\R^d$
  such that $\DVT(P)$ can be obtained by a suitable choice
  of points $p_e$ on the edges $e\subset P$.
  Assume that $\Stack(n+1,d)$ arises by stacking onto a facet of
  $\Stack(n,d)$ that is realized by the facet $F\subset P$
   with vertex set $\{v_1,\dots,v_d\}$.
  The ``new'' vertex $w$ is now chosen 
  ``beyond'' the facet $\DVT(F)$ of~$\DVT(P)$, and
  ``beneath'' all other facets of~$\DVT(P)$.
  That is, addition of~$w$ to $\DVT(P)$ would mean stacking a pyramid
  onto the facet $\DVT(F)$ of~$\DVT(P)$.
  In particular, $w$ lies 
  ``beyond'' the facet $F$ of~$P$, and
  ``beneath'' all other facets of~$P$,
  so $P':=\conv(\{w\}\cup P)$ is a stacked polytope
  realizing $\Stack(n+1,d)$, as required.

  The facet hyperplanes $H_{v_i}$ of~$\DVT(P)$
  cut the edges $[v_i,w]$ of~$P'$ in points $p_i$: 
  This is since $w$ is beneath $H_{v_i}$, while $v_i$ is cut off by~$H_{v_i}$.
  Thus we obtain points $p_i$ on the new edges
  of~$P'$, and the hyperplane $H_w:=\aff\{p_1,\dots,p_d\}$
  may be taken to cut off the new vertex $w$ of~$P'$.
  This new truncation plane is determined uniquely by the $d$
  intersection points, because the new vertex $w$ of~$P'$ is simple.
\end{proof}

\begin{figure}[ht]
\begin{center}
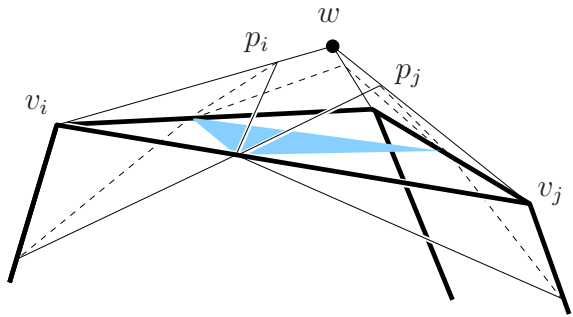
\end{center}
\caption{The induction step in Theorem~\ref{thm:CutStack}, for $d=3$.
$\DVT(F)$ is drawn shaded.}
\label{figure:InductionStack}
\end{figure}

This theorem is valid for all $d\ge3$; in particular, 3D-pictures
work. (Figure~\ref{figure:InductionStack} is a feeble attempt.)
However, the construction produces by far the most interesting results
for $d=4$.

\begin{corollary}[\cite{Z89}]
  For each $n\ge0$, and for every type of stacked $4$-polytope 
  $\Stack(n,4)$ with $f$-vector
  $(5+n, 10+4n, 10+6n, 5+3n)$, there is a corresponding
   $2$-simple $2$-simplicial $4$-polytope $\DVT(\Stack(n,4))$,
   with $f$-vector 
\[
f(\DVT(\Stack(n,4)))\ =\ (10+4n, 30+18n, 30+18n, 10+4n).
\]
\end{corollary}

In particular, this yields infinitely many combinatorial types
of $2$-simple $2$-simplicial $4$-polytopes. Moreover,
with a bit of care 
the proof of Theorem~\ref{thm:CutStack} yields 
these polytopes with rational vertex coordinates.
See \cite{Paffenholz:Webpage2s2s} for explicit examples
of such coordinates.

\begin{corollary}[\cite{Z89}]
  The number of combinatorial types of $2$-simple $2$-simplicial
  $4$-polytopes with $10+4n$ vertices grows exponentially in~$n$.
\end{corollary}

See Paffenholz \& Werner \cite{PaffenholzWerner:many} for further
constructions of $2$-simple $2$-simplicial $4$-polytopes with interesting 
$f$-vectors. In particular, they describe the ``smallest'' example
of such a polytope (other than the simplex), which has only
$9$ vertices.

\subsection*{Exercises}
\begin{compactenum}[\thechapter.1.]
\item Show that any simple or simplicial $d$-polytope with $f_0=f_{d-1}$
  must be a simplex, or $2$-dimensional.
\item\label{ex:compute-f}
Compute the full $f$-vectors, as well as the number
$f_{03}$ of vertex-facet incidences, for the following $4$-polytopes, based
only on the information given here:
\begin{compactenum}[(a)]
\item The $24$-cell: a 2s2s polytope whose facets
  are $24$ octahedra;
\item The $600$-cell: a simple polytope whose facets are $120$ dodecahedra;
\item The $720$-cell: a 2s2s $4$-polytope whose facets are $720$ bipyramids
   over pentagons;
\item A neighborly cubical polytope NCP$^n_4$, a cubical polytope
  with the graph of the $n$-cube ($n\ge4$).
\end{compactenum}
\item\label{ex:compute-f2}
Compute the full $f$-vectors of the stacked $d$-polytopes $\Stack(n,d)$.
\item Show that if a $4$-polytope $P$ is not simplicial, then
  $\DVT(P)$ cannot be $2$-simplicial.
\item Find coordinates for $\DVT(\Stack(1,4))$. Check them with
  \texttt{polymake}. 
\\
  (This is Braden's ``glued hypersimplex'' \cite{braden97}.)
\item\label{ex:exp-stacked}
  Show that there are exponentially many distinct combinatorial 
  types of stacked $d$-polytopes with $d+1+n$ vertices, for any $d\ge3$.
  Derive that there are exponentially many types
  of $2$-simple $2$-simplicial $4$-polytopes with the same $f$-vector.
\item\label{ex:CF1}
Show that 
$f_{13}=f_{03}+2f_2-2f_3$, and dually
$f_{02}=f_{03}+2f_1-2f_0$,
holds for the flag vector of each $4$-polytope.\\
(Hint: Sum the Euler equations for the facets, which are $3$-polytopes.)
\\
Derive from this that the inequality
$2f_{03}\ge(f_1+f_2)+2(f_0+f_3)$
is valid for all $4$-polytopes, and that it is tight 
exactly for the $2$-simple $2$-simplicial $4$-polytopes.
\item
Show that there is no $f$-vector inequality
(not involving $f_{03}$) that
characterizes the 2s2s $4$-polytopes.
\item If $P$ is a $d$-dimensional simplicial polytope, and if
  $\DVT(P)$ exists, is $\DVT(P)$ then $2$-simple? $2$-simplicial?
\end{compactenum}

\lecture{{\itshape f}\/-Vectors of 4-Polytopes}

The $f$-vector of a $4$-polytope is a quadruple of integers
$f(P)=(f_0,f_1,f_2,f_3)$,
but due to the Euler-Poincar\'e relation 
the set of all $f$-vectors of $4$-polytopes is a $3$-dimensional
set: It lies on the ``Euler-Poincar\'e hyperplane''
in~$\R^4$, given by
\[f_0-f_1+f_2-f_3\ =\ 0.\]
The task we are facing is to describe the set of all $f$-vectors,
\[ 
\mathcal{F}_4\ :=\ 
\{f(P)\ =\ (f_0, f_1, f_2, f_3)\in\Z^4: P\textrm{ a convex 4-polytope}\}.
\]
Here ``describe'' may mean a number of different things:
Probably one should not hope for a complete
description (as Steinitz got for the $3$-dimensional case),
since the set of $f$-vectors is way more complicated
in the $4$-dimensional case. 

\begin{figure}[ht]
\begin{center}
\input{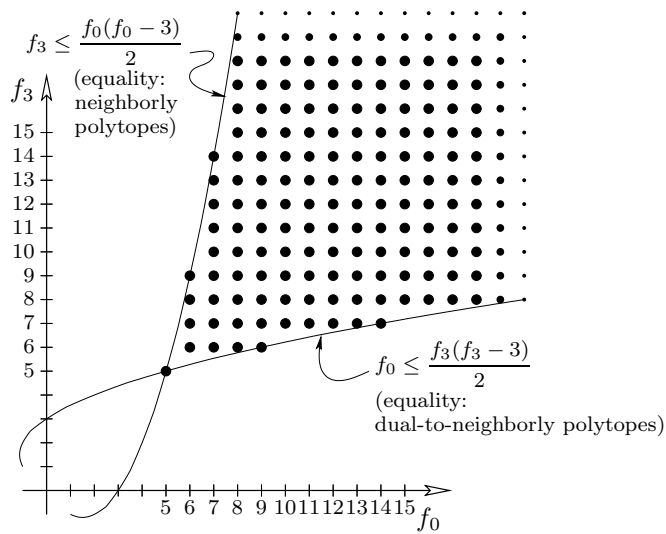}
\end{center}
\caption{The $(f_0,f_3)$-pairs of convex $4$-polytopes, according
to Grünbaum \cite[Sect.~10.4]{Gr1-2})}
\label{figure:f0f3plane}
\end{figure}

Indeed, $\mathcal{F}_4$ 
is \emph{not} the set of all integral points in a polyhedral
cone, or even in a convex set. This may be seen
from the characterizations
of the projections of~$\mathcal{F}_4$ to the coordinate
$2$-planes in~$\R^4$, by Grünbaum, Barnette, and Reay
\cite[Sect.~10.4]{Gr1-2} \cite{barnette74:_e_s} \cite{barnette73:_projec},
which show non-convexities and holes (see Figure~\ref{figure:f0f3plane}).
Or you just note that some of the rather basic, tight inequalities,
such as the upper bound inequality $f_1\le\binom{f_0}2$, are concave.
For example,
\begin{eqnarray*}
f(C_4(5))&=& (5,10,10,\ 5),\\
f(C_4(7))&=& (7,21,28, 14),\\
f(C_4(9))&=& (9,36,54, 27).
\end{eqnarray*}
The midpoint of the segment between $f(C_4(5))$ and~$f(C_4(9))$
is the integral point $(7,23,32,16)$: It violates the
upper bound inequality, and indeed a $4$-polytope with
$7$ vertices cannot contain more than the $21=\binom72=f_1(C_4(7))$ edges.
(See also Bayer \cite{Bay}, Höppner \& Ziegler \cite{Z59}.)

In the following, we will head for a complete description of
the $f$-vector cone for $4$-polytopes, $\cone(\mathcal{F}_4)$.
This seems to be a challenging but realistic goal.
Once that is achieved  (the ``2006 project''), a
logical next goal might be a description
of the ``large'' $f$-vectors, that is, of
\[ 
\{f(P)=(f_0, f_1, f_2, f_3)\in\Z^4: P \textrm{ a convex 4-polytope},
\ f_0+f_3\ge M\}
\]
for some large~$M$. But let's not get too ambitious too fast.

\section{The \boldmath$f$-vector cone}

\begin{definition}[$f$-vector cone]\label{def:fcone}
  The \emph{$f$-vector cone} of $4$-polytopes, $\cone(\mathcal{F}_4)$,
  is the topological closure of the convex cone with apex
  $f(\Delta_4)=(5,10,10,5)$ that is spanned by the $f$-vectors of
  $4$-polytopes,
\[
\Big\{f(\Delta_4)+\sum_{i=1}^N\lambda_i\,\big(f(P_i)-f(\Delta_4)\big):
P_1,\dots,P_N\textrm{ 4-polytopes},\ 
\lambda_1,\dots,\lambda_N\ge0\Big\}.
\]
  Equivalently, $\cone(\mathcal{F}_4)\subset\R^4$ is the solution set
  to all the linear inequalities that are valid for all $f$-vectors
  for $4$-polytopes, and that are tight at the $f$-vector of the simplex.
\end{definition}

The equivalence between the two versions of the definition
rests on basic facts about closed convex sets, which you should
put together yourself (\exercise{ex:defsequiv}). 
You are also asked to verify that the cone generated
by the $f$-vectors is not closed, so we do have to take the topological
closure (\exercise{ex:notclosed}.)

\begin{figure}[ht]
\begin{center}
\input{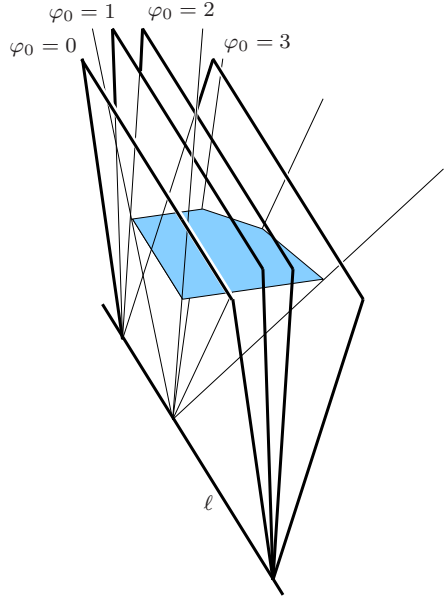}
\end{center}
\caption{The function $\varphi_0$ is constant on certain planes that contain
  the apex of the cone. It is not defined on the line $\ell$ where all
  those planes intersect. (In terms of $(f_0,f_1,f_2)$-coordinates,
  $\ell$ is defined by $f_0=5$ and $f_1+f_2=20$.)}
\label{figure:phi0oncone}
\end{figure}

The closed convex cone we are looking at is $3$-dimensional,
so we may view it as the cone over a $2$-dimensional convex figure,
which might be just a pentagon or hexagon.
Instead of looking at a $2$-dimensional section
(say intersecting by $f_1+f_2=100$), we may equivalently
introduce homogeneous (``projective'') coordinates,
which are rational linear functions, normalized to yield
``$\frac00$'' at the $f$-vector of a simplex (compare Lecture~1).
There is no unique best way to do this; we choose
\[
\varphi_0\ :=\ \frac{f_0-5}{f_1+f_2-20}\qquad\textrm{and}\qquad
\varphi_3\ :=\ \frac{f_3-5}{f_1+f_2-20}
\]
as our homogeneous coordinates. (Figure~\ref{figure:phi0oncone}
illustrates the geometry of such a rational function on a
cone.)
So we are trying to describe $\proj(\mathcal{F}_4)\subset\R^2$, the
closure of
\[
\conv\{(\varphi_0(P),\varphi_3(P))\in\R^2:P\textrm{ a convex 4-polytope}\}.
\]

Any $4$-polytope yields a (rational) point 
in the $(\varphi_0,\varphi_3)$-plane.
Any valid linear inequality, tight at the $4$-simplex, translates
into a linear inequality in $\varphi_0$ and~$\varphi_3$.
So let's look at some families of polytopes
and of linear inequalities that we know, and let's see what they buy us.
\medskip

\noindent\textbf{Some 4-polytopes we know:}
\[\begin{tabular}{lccc}
Stacked:   & $(5+\hphantom{3}n, 10+4n, 10+6n, 5+3n)$  & $\longrightarrow$ & 
                       $(\frac{1}{10}, \frac{3}{10})$\\[1pt]
Truncated: & $(5+3n, 10+6n, 10+4n, 5+\hphantom{3}n)$  & $\longrightarrow$ & 
                       $(\frac{3}{10}, \frac{1}{10})$\\
Cyclic:    & $(n,\frac{n(n-1)}2,n(n-3),\frac{n(n-3)}2)$ & 
$\stackrel{n\rightarrow\infty}\longrightarrow$ &
                       $(0,            \frac{1}{3})$\\
Dual-to-cyclic: &  
             $(\frac{n(n-3)}2,n(n-3),\frac{n(n-1)}2,n)$ & 
$\stackrel{n\rightarrow\infty}\longrightarrow$ &
                       $(\frac{1}{3},  0          )$.
\end{tabular}\]
The truncated polytopes are the duals of the stacked polytopes,
so they are simple. Similarly, the duals of cyclic polytopes are simple.
Thus we find the four points 
$(\frac{1}{10}, \frac{3}{10})$,
$(\frac{3}{10}, \frac{1}{10})$, 
$(0,            \frac{1}{3})$, and
$(\frac{1}{3},  0          )$,
which span a quadrilateral subset of $\proj(\mathcal{F}_4)$.
This quadrilateral also represents the 
$f$-vectors of simple and of simplicial polytopes and ``everything in between.''
(Note that duality interchanges the coordinates $\varphi_0$ and~$\varphi_3$,
and thus 
$\proj(\mathcal{F}_4)$ is symmetric with respect to the main diagonal.)
\medskip

\noindent\textbf{Five linear constraints we know:}
\[\begin{tabular}{llcc}
``Few Vertices'': & $f_0\ge 5$  & $\Longleftrightarrow$ & $\varphi_0\ge 0$,\\
``Few Facets'':   & $f_3\ge 5$  & $\Longleftrightarrow$ & $\varphi_3\ge 0$,\\
``Simple'':       & $f_1\ge 2f_0$&$\Longleftrightarrow$ & 
                                            $3\varphi_0+\varphi_3 \le 1$,\\
``Simplicial'':   & $f_2\ge 2f_3$&$\Longleftrightarrow$ & 
                                            $\varphi_0+3\varphi_3 \le 1$,\\
``Lower bound'':  & $2f_1+2f_2\ge 5f_0+5f_3-10$&$\Longleftrightarrow$ & 
                                            $\varphi_0+\varphi_3  \le\frac25$.
\end{tabular}\]
The first four inequalities are quite trivial, and we have named them
by the polytopes that satisfy them with equality, at least asymptotically.
The translation into $(\varphi_0,\varphi_3)$-inequalities, using
the Euler-Poincar\'e relation, poses no problem.
There is no polytope with $\varphi_0=0$,
but the condition is satisfied asymptotically by any
family of $4$-polytopes with far more vertices than facets.
For example, the products of $n$-gons, with
$f(P_n\times P_n)=(n^2,2n^2,n^2+2n,2n)$, yield
$(\varphi_0,\varphi_3)=(\frac{n^2-5}{3n^2+2n-20},\frac{2n-5}{3n^2+2n-20})
\in\proj(\mathcal{F}_4)$,
which in the limit $n\rightarrow\infty$ yields $(\frac13,0)$.

The one non-trivial inequality in our table above is the last one, a
``Lower Bound Theorem.'' It may be derived quite easily \cite{Bay}
from the inequality $f_{03}\ge3f_0+3f_3-10$, which was first
established by Stanley \cite{Sta7} in terms of the so-called toric
$g$-vector (it is the inequality ``$g_2^{\mathrm{tor}}(P)\ge0$''); a
proof via rigidity theory was later given by
Kalai~\cite{kalai87:_rigid_i}.

\begin{figure}[b]
\begin{center}
\input{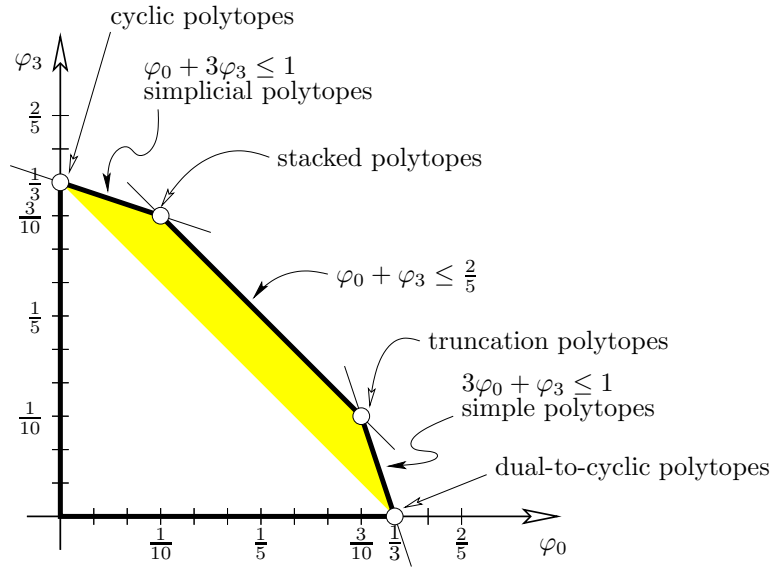}
\end{center}
\caption{Projective representation of the $f$-vectors of $4$-polytopes, in the
  $(\varphi_0,\varphi_3)$-plane. The convex set $\proj(\mathcal{F}_4)$
  is contained in the bold pentagon; it contains the shaded trapezoid.}
\label{figure:phiplane}
\end{figure}

Figure~\ref{figure:phiplane} summarizes our discussion up to this
point: We are interested in~$\proj(\mathcal{F}_4)$, the closure of
the set
\[
\conv\{(\varphi_0(P),\varphi_3(P)): 
P\textrm{ is a 4-polytope, not a simplex}\,\}\ \subset\ \R^2.
\]
This set is contained in the pentagon cut out by the five linear
inequalities discussed above, and it contains the shaded trapezoid,
which represents ``everything between simple and simplicial
polytopes.''  Indeed, simple and simplicial polytopes satisfy the
additional linear inequality $\varphi_0+\varphi_3\ge\frac13$.

\noindent
Thus we are left with the following ``upper bound problem'':

\begin{UBT}
Are there $4$-polytopes with $\varphi_0+\varphi_3\rightarrow 0$\,?
\end{UBT}

The inequality $\varphi_0+\varphi_3\ge\frac13$ is certainly not valid
for all (possibly non-simple non-simplicial) $4$-polytopes: Already for the
hypersimplex we get $(\varphi_0,\varphi_3)=(\frac18,\frac18)$.\break

\noindent
However, currently it is not clear how small $\varphi_0+\varphi_3$ 
can be for convex polytopes. 
Thus the Upper Bound Problem is the key remaining problem
in the description of the $f$-cone for $4$-polytopes.
\begin{compactitem}[~(!)~]
\item
If the answer is YES to the problem as posed above,
then the five inequalities above constitute
a complete linear description of~$\cone(\mathcal{F}_4)$. 
\item
If the answer is NO, then this is also exciting, since it 
means that the answers for cellular spheres
and for convex polytopes are distinct!
Indeed, cellular $3$-spheres with arbitrarily small $\varphi_0+\varphi_3$
have been constructed by Eppstein, Kuperberg \& Ziegler \cite{Z80};
see our discussion in Section~\ref{sec:LBP}.
\end{compactitem}

\section{Fatness and the Upper Bound Problem}

We prefer to rephrase the Upper Bound Problem
in terms of a somewhat more graphic quantity,
which we call ``fatness.''

\begin{figure}[ht]
\begin{center}
\input{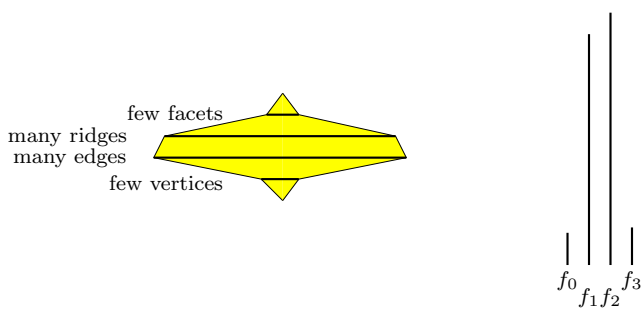}
\end{center}
\caption{Fatness for a $4$-polytope face lattice, and for an $f$-vector}
\label{figure:fatness}
\end{figure}

\begin{definition}[Fatness]
The \emph{fatness} of a $4$-polytope is the quotient
\[ F(P)\ :=\ \frac{1}{\varphi_0+\varphi_3}\ =\            
\frac{f_1+f_2-20}{f_0+f_3-10}.
\]
\end{definition}

The fatness of a $4$-polytope is large if both $\varphi_0$ and~$\varphi_3$
are small. This happens if the polytope has relatively few vertices
and facets, but many edges and $2$-faces. Thus, graphically, the
face lattice and the $f$-vector are ``fat in the middle,'' whence the
name (see Figure~\ref{figure:fatness}).

\begin{UBT}
Can the fatness of a $4$-polytope be arbitrarily large?
\end{UBT}

Here are a few explicit values to start with:
For stacked and truncated $4$-polytopes we have
$F(P)=\frac{5}{2}$ exactly.
For cyclic polytopes we get $F(C_4(n))\rightarrow3$
for $n\rightarrow\infty$, and the same for the duals ---
fatness is a self-dual quantity, that is, any $4$-polytope and its dual
have the same fatness.
Moreover, it is easy to compute 
(or to derive from Figure~\ref{figure:phiplane}) that all simple
and simplicial polytopes satisfy $\frac{5}{2}\le F(P)<3$.

But how large can fatness be? The attempts to answer this question 
have led to a multitude of interesting examples and constructions,
and to a fast succession of record holders
for ``the fattest examples found so far.''
Many of them can be obtained by deep vertex truncation
of simplicial polytopes,
so they are $2$-simple and $2$-simplicial by Proposition~\ref{prop:DVT2s2s}:
\begin{compactitem}[--]
\item The hypersimplex, which is the dual of $\DVT(\Delta_4)$, has fatness $4$.
\item Schäfli's $24$-cell \cite{Schla}, $\DVT$(cross polytope), 
                                          has fatness~$4.526$.  
\item G\'evay's $720$-cell~\cite{Gevay}, the dual of $\DVT$($120$-cell), 
has $720$ facets that are bipyramids over regular pentagons.
It has fatness~$5.020$. 
\item 
Eppstein, Kuperberg \& Ziegler \cite{Z80} used hyperbolic geometry
arguments to achieve a fatness of~$5.048$ by
the ``E-construction,'' which is dual to deep vertex truncation.
\item
Paffenholz \cite{paffenholz-pc}
has very recently shown that there are realizations
for any sum of an $n$-gon and $m$-gon such that the
deep vertex truncation $\DVT(P_m\times P_n)$ can be obtained.
For $m=n\rightarrow\infty$ this yields fatness approaching~$6$.
\end{compactitem}
However, we'll go a different route. In the next and final lecture
we will present a construction that generalizes and extends
the construction of ``neighborly cubical'' $4$-polytopes of
Joswig \& Ziegler~\cite{Z62}, to achieve fatness arbitrarily
close to~$9$, the latest record (as far as I know at the time of
writing). I would have been happy to have a ``note added in proof''
about this~\dots

\begin{figure}[ht]
\begin{center}
\input{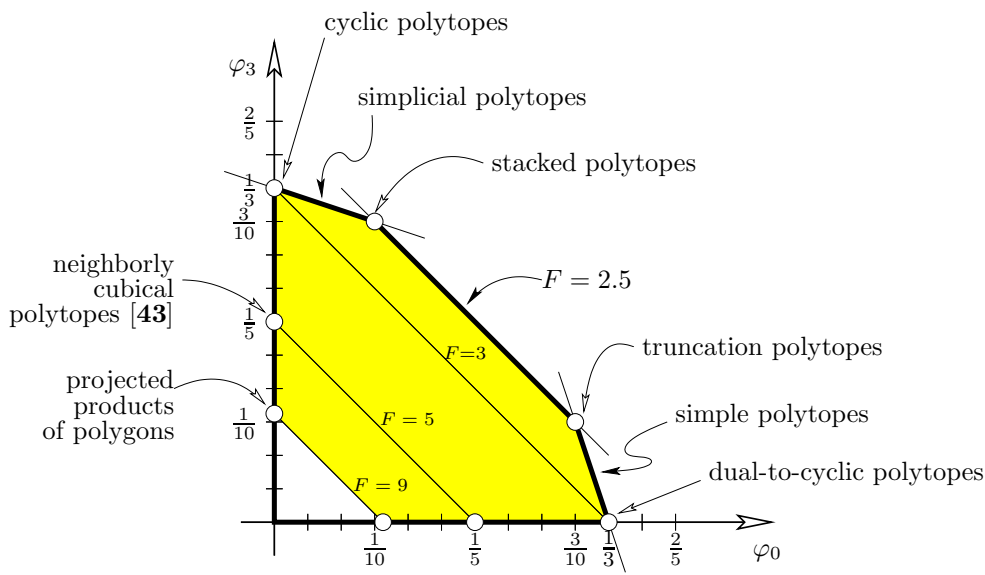}
\end{center}
\caption{$4$-polytopes in the
  $(\varphi_0,\varphi_3)$-plane.
  The shaded hexagon is spanned by the
  $(\varphi_0,\varphi_3)$-pairs of known $4$-polytopes.}
\label{figure:phiplane2}
\end{figure}

What do polytopes ``of very high fatness'' look like?
You can verify (via \exercise{ex:complexity}) 
that they have two properties:
\begin{compactenum}[~(1)~]
\item 
The facets have many vertices (on average).
\item
The vertices are in many facets (on average).
\end{compactenum}
Either of these properties are easy to satisfy ---
just look at the products $P_n\times P_n$ for the first property,
and at their duals, the free sums $P_n\oplus P_n$, for the second one.
The key question is whether they can \emph{simultaneously} be
satisfied.

Finally, here is a problem on $3$-dimensional polytopal tilings that
is ``essentially'' equivalent to the fatness problem:
Consider face-to-face tilings of~$\R^3$
(cf.~\cite{SchattschneiderSenechal})
that satisfy some
regularity properties, e.g.\ one of the following
(each implies the next):
\begin{compactitem}[~--~]
\item the tiling is triply periodic (that is, there are three linearly
independent translational symmetries),
\item
there are only finitely many distinct congruence classes of tiles,
\item
in- and circumradius of the tiles are uniformly bounded.
\end{compactitem}
For such tilings, we may define notions of ``average'' vertex degrees,
face numbers, etc.
The question is whether there is such a tiling where
the tiles have lots of vertices on average, \emph{and}
the vertices are in many tiles on average.
Again, either property is easy to achieve
(look at tilings by Schlegel diagrams),
but can they be \emph{simultaneously} satisfied?

\section{The Lower Bound Problem}\label{sec:LBP}

The upper bound problem discussed here has a natural 
``lower bound'' counterpart. It arises if we don't restrict ourselves to the
geometric model of convex polytopes, but consider
the larger class of cellular spheres that
are ``regular'' in the sense that their cells have no 
identifications on the boundary, and that satisfy
the ``intersection property'' that any two faces should
intersect in a single cell (which may be empty).
These are the regular CW spheres \cite{CookeFinney} whose face poset is
a lattice (where the meet operation corresponds to intersection of faces).

\begin{LBT}
Does $\varphi_0+\varphi_3\le\frac25$ hold for the cellular spheres
 that satisfy the intersection property?
\end{LBT}

This problem seems crucial in terms of the separation of the 
``geometric'' model of convex polytopes from the
``topological'' model of cellular spheres/balls.
\begin{compactitem}[~(!)~]
\item
If the answer to the problem is NO, then this would establish such
a separation, which would be quite remarkable.
\item
If the answer is YES, then this would imply a complete characterization
of the $f$-vector cone for cellular $3$-spheres, by the
five linear inequalities given above; indeed,
Eppstein, Kuperberg \& Ziegler \cite{Z80} have constructed
cellular spheres for which fatness is arbitrarily large,
that is, $\varphi_0+\varphi_3$ is arbitrarily small.
\end{compactitem}
We will not discuss this here further, but refer to \cite{Z80}
and~\cite{Z82}.

\subsection*{Exercises}
\begin{compactenum}[\thechapter.1.]
\item\label{ex:defsequiv} Show that the two definitions of the
  $f$-vector cone
  given in Definition~\ref{def:fcone} are indeed equivalent.\\
  Hint: You need a separation lemma; see for 
  example Matou\v{s}ek~\cite[p.~6]{matousek02:_lectur_discr_geomet}.
\item\label{ex:notclosed}
  Show that the union of the
  line segments $[f(\Delta_4),f(C_4(n))]\subset\R^4$ has the 
  whole ray $\{(5,10+t,10+2t,5+t):t\ge0\}$ in its closure.
  Note that $f_0\ge5$ is a valid linear inequality, which is tight at
  $f(\Delta_4)$, but for no other $f$-vector.\\ 
  Conclude that the cone with apex $f(\Delta_4)$ spanned by the $f$-vectors
  of $4$-polytopes is not closed.
\item
Compute the fatness and the $(\varphi_0,\varphi_3)$-pair
for the hypersimplex, the $24$-cell, and for $\DVT(600$-cell). 
\item
Compute the fatness of the 2s2s polytopes
$\DVT(\Stack(n,4))$, and show that it lies in the interval $[4,\,4.5)$.
\newline
Show that for any simplicial $4$-polytope $P$, the fatness of
$\DVT(P)$ is smaller than~$6$.
\newline
Where would the $f$-vectors of the polytopes $\DVT(P)$ lie
in $\proj(\mathcal{F}_4)$, as graphed in Figure~\ref{figure:phiplane2}?
\item
If $C^n_4$ is a cubical $4$-polytope with the graph of an $n$-cube 
(see Exercise~3.\ref{ex:compute-f}), compute the fatness and
the pair $(\varphi_0,\varphi_3)$.
\item\label{ex:complexity}
Define the \emph{complexity} of a $4$-polytope to be the quotient
\[
C(P)\ :=\ \frac{f_{03}-20}{f_0+f_3-10}.
\]
\begin{compactenum}[(a)]
\item
Show that $F(P)\le 2C(P)-2$, with equality if and only if $P$ 
is $2$-simple and $2$-simplicial.
\item
  Show that $C(P)\le 2F(P)-2$, with equality if and only if
if all facets of~$P$ are simple,
or equivalently, if all vertex figures are simple.
\item
Derive from this that fatness is high if and only if both
the average number of vertices per facet, $f_{03}/f_3$, and 
the average number of facets per vertex, $f_{03}/f_0$,
are large.
\end{compactenum}
\end{compactenum}

\lecture{Projected Products of Polygons}

In this lecture we present a construction
of very recent vintage, ``projected products of polytopes.''
We will not have the ambition 
to work through all the technical details for this; these
appear in~\cite{Z97}, see also \cite{Z102};
rather, our main objective is here to identify the structural
features of the construction which lead to fat polytopes, and to outline 
(possibly ``for further use'') some interesting
components that go into the construction.
In the following version of the main result,
some concepts that will be explained below are highlighted 
by quotation marks.

\begin{theorem}[Ziegler~\cite{Z97}]\label{thm:PPP}
  For each $r\ge 2$, and even $n\ge 4$, there is a realization
  $P^r_n\subset\R^{2r}$ of a product of polygons $(P_n)^r$ (a
  ``deformed product of~polygons'') such that the vertices and edges
  and all the ``$n$-gon $2$-faces'' of~$P^r_n$ ``survive'' the
  projection $\pi:\R^{2r}\rightarrow\R^4$ to the last $4$ coordinates.
\end{theorem}

A number of nice tricks go into the 
construction that proves the theorem --- see below.
Before we look into these we want to derive the enumerative consequences:
The $f$-vector of~$\pi(P^r_n)$ can be derived purely from the information
given in the theorem, not using details about the combinatorics 
of the resulting polytopes (which were worked out only 
recently \cite{dipl-Sanyal} \cite{Z102}).

\section{Products and deformed products}

We have discussed the construction and main properties of
products of polytopes already in Example~\ref{ex:products}. 
A key observation is that the non-empty faces of a product are 
the products of non-empty faces of the ``factors.''
Now we specialize to the case of products of (several) polygons:
We consider products of~$r$ $n$-gons --- and later we
will be looking at polytopes that just have the combinatorics
of such polytopes.

If $P_n$ is an $n$-gon, then $(P_n)^r$ is a simple polytope
of dimension $2r$. 
It has
\begin{compactitem}[~$\bullet$~]
  \item $f_0=n^r$ vertices (of the form ``vertex $\times$ vertex $\times$
    \ldots \ldots \ldots $\times$ vertex''), and
  \item $f_1=rn^r$ edges (of the form 
 ``vertex $\times$ \ldots $\times$ edge $\times$ \ldots $\times$ vertex).
\end{compactitem}
The products of polygons have two different types of $2$-faces,
``quadrilaterals'' and ``polygons,'' that we need to distinguish:
\begin{compactitem}[~$\bullet$~]
  \item $\binom r2 n^r$ \emph{quadrilateral} $2$-faces (which arise
    as products of two edges, and vertices from the other factors), and
  \item $rn^{r-1}$ \emph{polygon} $2$-faces (arising as a product of
    one $n$-gon factor with vertices from the other factors).
\end{compactitem}
Thus we get $f_2=\binom r2 n^r+rn^{r-1}$.

Finally, let's note that $(P_n)^r$ has $nr$ facets,
which arise as a product of one edge (from one of the factors)
with $n$-gons as the other factors.

The ``deformed products'' of polygons $P^r_n$ considered below
are combinatorially equivalent to the
``orthogonal products'' $(P_n)^r$. Thus, the $f$-vector count given
here is valid for deformed products as well.

\section{Computing the \boldmath$f$-vector}

It is remarkable that for the proof of the following corollary to
Theorem~\ref{thm:PPP} we don't need detailed combinatorial 
information about $\pi(P^r_n)$; it is sufficient to know
that it is a generic projection of a deformed product $P^r_n$,
and that all the polygon $2$-faces survive the projection.
The facets of~$\pi(P^r_n)$ are $3$-faces of~$P^r_n$ that 
``have survived the projection,'' that is, they are 
combinatorial cubes and prisms over $n$-gons.

\begin{corollary}
  $\pi(P^r_n)\subset \R^4$ is a $4$-polytope with $f$-vector
\begin{eqnarray*}
&&
\textstyle\big(n^r, rn^r, 
\frac{5}{4}rn^r-\frac{3}{4}n^r+rn^{r-1},
\frac{1}{4}rn^r-\frac{1}{2}n^r+rn^{r-1}\big)\\
&=&
\textstyle\big(0+\frac{4}{r}, 4, 
5- \frac{3}{r}+ \frac{4}{n}, 
1- \frac{2}{r}+\frac{4}{n}             \big)\cdot\frac{1}{4}rn^r.
\end{eqnarray*}
In particular, for $n,r\rightarrow\infty$ the fatness of 
the projected products of polygons $\pi(P^r_n)$
gets arbitrarily close to $9$.
\end{corollary}

\begin{proof}
  From the combinatorial information above we see that 
  each vertex, and each edge, of a product of polygons $(P_n)^r$ is contained
  in a polygon $2$-face: So if all the polygon $2$-faces
  survive the projection, then in particular all vertices and edges
  do. Thus, for $(f_0,f_1,f_2,f_3):=f(\pi(P^r_n))$ 
  we get $f_0=n^r$, and $f_1=rn^r$.
 
  Furthermore, we know that $\pi(P^r_n)$ has 
  $p:=rn^{r-1}$ polygon $2$-faces (and a yet unknown number of 
  quadrilaterals).
  Moreover, the facets of~$\pi(P^r_n)$ are prisms over polygons, and
  cubes.
  Each prism has two polygon faces (and $n$ quadrilateral faces),
  while each cube has $6$ quadrilateral $2$-faces, but no polygon
  faces.
  So, each polygon (ridge!)\ lies in two prism facets, and
  each prism facet contains two polygons: Double counting yields
  that there are exactly $p=rn^{r-1}$ prism facets.

  Denote by $c$ the number of cube facets. 
  Then we have 
  \begin{compactitem}[~$\bullet$~]
    \item $f_3=c+p$ (all facets are cubes or prisms),
    \item $2f_2=6c+(n+2)p$ (double counting the facet-ridge
      incidences), and
    \item $f_2-f_3=(r-1)n^r$ (Euler's equation).
  \end{compactitem}
  These three linear equations can now be solved for the three 
  unknowns, $f_2$, $f_3$, and~$c$.
\end{proof}

\section{Deformed products}

A number of different polytope constructions 
have been studied in attempts to produce interesting
polytopes. For example, \emph{any} polytope may 
be obtained by projection of a simplex;
on the other hand, any  projection of an orthogonal
cube is a zonotope, which
is a polytope of very special structure (see \cite[Lect.~7]{Z35}).
The projections of orthogonal products of centrally symmetric polygons
are still zonotopes, and the projections of 
orthogonal products of arbitrary polygons
are only a bit more general.
However, it has been noted since the seventies (in the context of
linear programming, trying to construct ``bad examples'' for the
simplex algorithm; cf.\ \cite{KlMi} \cite{Z51}) that some
projections of ``deformed products'' have very interesting extremal
properties.
For a very simple example, look at a $3$-cube
(which is a product of three $1$-polytopes).
Any orthogonal cube projected to the plane will produce a
hexagon (at best), while a deformed cube can be projected to the
plane to yield an octagon: All the vertices ``survive the
projection'' (see Figure \ref{figure:Projectprism}).

\begin{figure}[ht]
\begin{center}
\input{EPS/deformedcube.pstex_t}
\end{center}
\caption{A combinatorial 
$3$-cube, in a deformed realization such that all eight vertices
``survive'' the projection to the plane}
\label{figure:Projectprism}
\end{figure}

And indeed, it is not so hard to show that one can realize
an $n$-cube in such a way that all its vertices ``survive the 
projection'' to the plane --- this was first proved
by Murty \cite{murty80:_comput} and Goldfarb~\cite{goldfarb94};
in the context of linear programming it yields exponential examples for
the ``shadow vertex'' pivot rule for linear programming.

Similarly, if you project an orthogonal product of polygons to
$\R^4$, you cannot expect that all the edges survive the projection,
but with deformed products, this is possible (although hard to 
visualize --- proofs are mostly based on linear algebra criteria rather
than on geometric intuition).

Here are linear algebra descriptions of the polytopes
we'll be looking at:

\subsubsection*{Polygons:}
If $V$ is any $(n\times 2)$-matrix whose
rows are non-zero, (w.l.o.g.)\ ordered in cyclic order,
and positively span $\R^2$,
then for a \emph{suitable} positive right-hand side vector $b\in\R^n$ 
the system $Vx\le b$ describes a convex
$n$-gon $P_n\subset\R^2$. 

\subsubsection*{Product of polygons:}
Given such a $V$, we immediately get the system
\[
\left(\raisebox{-39mm}{\input{EPS/matrix0a.pstex_t}}\!\right)x 
\ \ \le\ \ 
\left(\raisebox{-39mm}{\input{EPS/vector0a.pstex_t}}\!\right),
\]
with a block-diagonal matrix of size $rn\times 2r$,
which describes the orthogonal product of polygons $(P_n)^r\subset\R^n$.

The combinatorics of the product of polygons is reflected
in the facet-vertex incidences, as follows:
Each inequality defines a facet;
each of the $n^r$ vertices is the unique
solution of a linear system of equations that is
obtained by requiring that from each block,
two cyclically adjacent inequalities are tight.

\subsubsection*{Deformed product of polygons:} 
Our Ansatz is as follows: 
\begin{equation}\label{eq:Ansatz}
\left(\raisebox{-55mm}{\input{EPS/matrix2.pstex_t}}\right)x
\ \ \le\ \ 
\left(\raisebox{-55mm}{\input{EPS/vectorgrey.pstex_t}}\right).
\end{equation}
If the diagonal 
$V^\eps$-blocks $V^\eps\in\R^{n\times2}$ satisfy the conditions above,
then the blocks $U,W\in\R^{n\times2}$ below the diagonal blocks can be
arbitrary --- the right-hand sides can be adjusted in such
a way that we get a \emph{deformed product}, which is
combinatorially equivalent to a product $(P_n)^r$.

For this, we only have to verify that we get the correct
combinatorics: If for each of the $r$ blocks we choose
two cyclically adjacent inequalities and
require them to be tight, then this should yield a 
linear system of equations with a unique solution,
for which all other inequalities are satisfied, but not tight.
It is now easy to prove (by induction on~$r$) that
our system satisfies this --- if we just choose the
right-hand sides suitably; in particular, $b_i:=M^{i-1}b\in\R^n$
works, for large enough~$M$, if
$V^\eps x\le b$ (as above) defines an $n$-gon.

\begin{figure}[ht]
\begin{center}
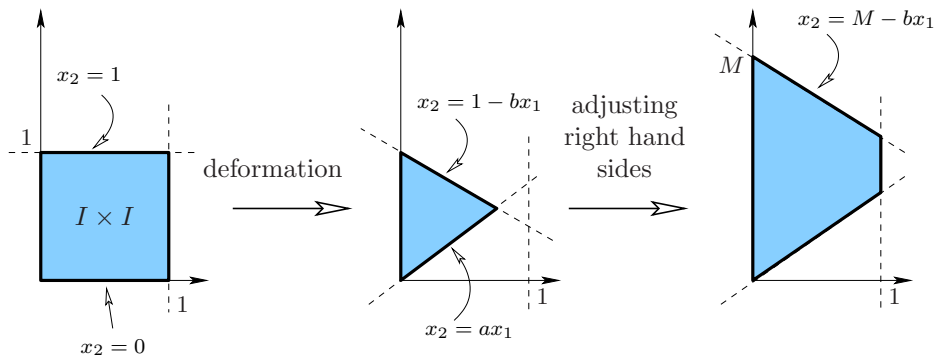
\end{center}
\caption{Construction of a deformation of $I\times I$}
\label{figure:deformedItimesI}
\end{figure}

\begin{example}
  To illustrate this in a simple case, let's consider deformed
  products in the low-dimensional case of~$I\times I$, a square, where
  $I$ denotes an interval (a $1$-dimensional polytope) such as
  $I=[0,1]$.

  In this case $I$ can be written as
\[\textstyle
I\ =\ \big\{x\in\R:
\left(\begin{array}{r}-1\\1\end{array}\right)x
\ \le\ 
\left(\begin{array}{r}0\\1\end{array}\right)\big\}.
\]
  Consequently, the product $I\times I$ may be represented
  by 
\[\textstyle
I\times I\ =\ \Big\{x\in\R^2:
\left(\begin{array}{rr}-1&0\\1&0\\0&-1\\0&1\end{array}\right)x
\ \le\ 
\left(\begin{array}{r}0\\1\\0\\1\end{array}\right)\Big\}.
\]
  Now changing the matrix into 
\[\textstyle
\left(\begin{array}{rr}-1&0\\1&0\\a&-1\\b&1\end{array}\right)x
\ \le\ 
\left(\begin{array}{r}0\\1\\0\\1\end{array}\right)
\]
leaves the first two inequalities (and thus the first $I$ factor)
intact, but it changes the slopes of the other two inequalities ---
and if you are unlucky (that is, for $a+b>1$)
the resulting polytope will not be equivalent to $I\times I$ any
more. This situation is depicted in the middle part of
Figure~\ref{figure:deformedItimesI}.
However, it can be remedied by increasing right-hand sides:
For any given $a$ and $b$, a suitably large $M$, namely $M>a+b$, in
\[\textstyle
\left(\begin{array}{rr}-1&0\\1&0\\a&-1\\b&1\end{array}\right)x
\ \le\ 
\left(\begin{array}{c}0\\1\\0\\M\end{array}\right)
\]
 will result in 
a product again (as in the right part of Figure~\ref{figure:deformedItimesI}). 
\end{example}

\section{Surviving a generic projection}

We are looking at \emph{generic projections} of \emph{simple polytopes}
 --- that is,
the direction of projection is in general position with
respect to all the edges of the polytope to be projected,
and thus a small perturbation of the inequalities of the polytope,
and of the direction of projection, will not change the
combinatorics of the projected polytope.

\begin{figure}[ht]
\begin{center}
\input{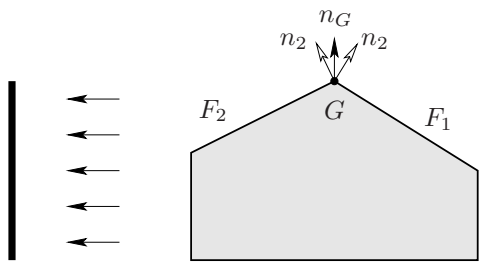}
\end{center}
\caption{A vertex $G$ surviving a projection, a normal vector $n_G$
orthogonal to the projection, and the normal vectors $n_1$ and $n_2$ of facets
$F_1$ and $F_2$ it can be combined from}
\label{figure:projectlemma}
\end{figure}

It should thus seem plausible that in this setting the following
holds (compare Figure~\ref{figure:projectlemma}) for a polytope
projection $\pi: P\to \pi(P)$:
\begin{compactitem}[~$\bullet$~]
\item
  The normal vectors to a face $G\subset P$ are the positive 
  linear combinations
  of the (defining) normal vectors $n_F$ to all the facets $F\supset G$
  that contain $G$.
\item The proper faces of~$\pi(P)$ are isomorphic copies of the faces
  of~$P$ that \emph{survive} the projection.
\item A face $G\subset P$ survives the projection exactly if it has
  a normal vector that is orthogonal to the direction of projection.
\end{compactitem}

\noindent
For Theorem \ref{thm:PPP} we have to specify
a deformed product realization for $(P_n)^r$ of the type
given in Ansatz~(\ref{eq:Ansatz}), such that
all $n$-gon $2$-faces are strictly preserved by the projection. That is, 
\begin{compactitem}[~$\bullet$~]
\item if we choose two cyclically adjacent rows from each block
      except for one, and
\item truncate these rows to the first $2r-4$ coordinates,
\end{compactitem}
then the resulting $2r-2$ vectors must be positively dependent and span.

\section{Construction}

\begin{figure}[ht]
\begin{center}
\includegraphics{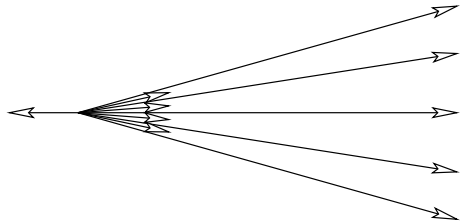}
\end{center}
\caption{A configuration of $10$ vectors, which positively span the
  plane, and which (in cyclic order) alternate between 
  vectors close to $(1,0)$, and close to~$(0,0)$}%
\label{figure:cyclic}
\end{figure}

Now we want to specify the lower-triangular block
matrix in our Ansatz (\ref{eq:Ansatz}) so that it satisfies
the following two main properties:
\begin{compactenum}[~~(1)]
  \item the diagonal blocks have ``rows in cyclic order,'' and
  \item any ``choice of two cyclically adjacent rows''
    from all but one of the blocks, truncated to the
    first $2r-4$ components, yields a positively-spanning
    set of vectors.
\end{compactenum}
Five observations (you may call them ``tricks'') help us to achieve this:
\begin{compactenum}[(i)]
\item Condition (2) is stable under perturbation.
  So, we first construct a matrix that satisfies (2), then 
  perturb it in order to achieve~(1).\\
  (The diagonal blocks of the matrix that we construct to
  satisfy (2) are denoted~$V$;  after perturbation, they will be~$V^\eps$.)

\item The submatrices $V$, $W$, and $U$ of size $n\times 2$
  are constructed to have alternating rows: So if you 
  choose two cyclically adjacent rows from a block,
  you know what you get!

  Specifically, we will let matrix $V$ have
  rows that alternate between $(1,0)$ and $(0,0)$,
  matrix $W$ gets rows $(0,1)$ and $(a,b)$, and
  matrix $U$ gets rows $(c,d)$ and $(e,f)$, with the six
  parameters $a,b,c,d,e,f$ to be determined.

\item To make sure that you get the positive linear
  dependence for (2), we specify a positive coefficient sequence
  and compute the matrix entries to satisfy them.

\item Rather than admitting that from one of the
  blocks no row is chosen, we will prescribe coefficient sequences
  that could have zeroes on any one of 
  the blocks, which yields linear dependencies for which the
  vectors from one block are ``not used.''

  Specifically, we take coefficient sequences of the form
  $\alpha_k:= (2^{k-t}-1)^2$ and 
  $\beta_k:= (2^{k-t}-1)(2^{k-t}-\frac32)$
  for the odd resp.\ even-index rows.
  These coefficients are clearly positive for integral $k$,
  except they vanish at $k=t$.

  Moreover, they are linear combinations of the three exponential
  functions $2^{k-t}$, $4^{k-t}$, and~$1$.
  If we write out the condition that
  ``the rows chosen should be dependent, with coefficients
  $\alpha_k$ (for the even-index row chosen from the $k$-th block) and
  $\beta_k$  (for the odd-index row chosen from the $k$-th block),
  then this leads to a system of six linear equations,
  in six unknowns $a,b,c,d,e,f$ --- solve it!

\item The properties ``alternating rows'' and
  ``rows in cyclic order'' are, of course, 
  incompatible --- but a matrix with rows in cyclic
  order can be obtained as a perturbation $V^{\eps}$
  of the matrix $V$ that has rows $(1,0)$ and $(0,0)$
  in alternation. Figure~\ref{figure:cyclic} suggests a way to do this.
\end{compactenum}
\medskip

\noindent
This completes our sketch of ``how to do it'' --- it should be
sufficient to let you construct the polytopes $P^r_n$ and thus prove
Theorem~\ref{thm:PPP} (but \cite{Z97} provides these details, too.)
\bigskip

If you think about this construction, can you perhaps manage to
simplify some of the details, or even to improve the construction?
After all, on the way to the characterization of the
$f$-vector cone for $4$-polytopes, this construction 
of polytopes of fatness up to~$9$ should be taken just as an intermediate
step. Can you get further? There is a long way to go
from $9$ to infinity~\ldots

\subsection*{Exercises}
\begin{compactenum}[\thechapter.1.]
\item Show that every polytope arises as a projection of a simplex.
\item Show that if $\pi:\R^n\rightarrow\R^d$, $P\rightarrow\pi(P)$
  is a polytope projection, for an $n$-polytope $P$, and $n>d$, then 
  the strictly preserved faces have dimension at most $d-1$.
  In particular, no facet of~$P$ is strictly preserved by the
  projection.
\item Give examples of polytope projections where no faces are
  strictly preserved.
\item Realize the prism over an $n$-gon ($n\ge3$) in such a way
  that the projection $\R^3\rightarrow\R^2$ strictly preserves all $2n$
  vertices.
\item Show that the product $\Delta_2\times\Delta_2$ of two triangles
  cannot be realized in such a way that all $9$ vertices are strictly
  preserved in a projection to~$\R^2$.
\item How large do we have to choose $r$ and $n$ in order
  to obtain $4$-polytopes of fatness $F(\pi(P^r_n))>8$?
\item Any projection of a (non-deformed) product of 
  centrally symmetric polygons is a \emph{zonotope}.
  For these, one knows that $f_1<3 f_0$
  (for such inequalities, for the dual polytope, 
   see~\cite[pp.~198/199]{Z10-2}).
  Deduce from this information that the fatness is smaller than~$5$.
\end{compactenum}

\lecture*{Appendix:\\ A Short Introduction to \ {\tt polymake}}%
\label{appendix:POLYMAKE}

\vskip-5mm\rightline{by Thilo Schröder and Nikolaus Witte}

\vspace{16mm}
\noindent
The software project \url{polymake}~\cite{GawrilowJoswig} has been
developed since 1997 in the Discrete Geometry group at TU Berlin by
Ewgenij Gawrilow and Michael Joswig, with contributions by several
others.  It was initially designed to work with
convex polytopes. Due to its open design the \url{polymake}
framework can also be used on other types of
objects; the current release includes a second
application, \url{topaz}, which treats simplicial complexes.

\url{polymake} is designed to run on any Linux or Unix system,
including Mac~OS~X.
It runs in a shell using command line input. This introduction is for
\url{polymake} versions~2.0 and~2.1.
\url{polymake} is free software and
you can redistribute it and/or modify it 
under the terms of the GNU General Public License as published by the
Free Software Foundation.

This introduction aims at getting you to work on the computer rather than
explaining the details about the machinery of \url{polymake}. Therefore,
there will be only a short description of the software design,
in Section~A.2. 
For further reading we refer to Gawrilow \& Joswig
\cite{GawrilowJoswig2} \cite{GawrilowJoswig4}. 
On the \url{polymake} website
\begin{center}
  \url{http://www.math.tu-berlin.de/polymake}
\end{center}
you will find extensive \emph{online documentation} as well as an
introductory \emph{tutorial}.

\section*{A.1. Getting Started}

For \url{polymake}, every 
polytope is treated as an \emph{object} (the file storing the data)
with certain \emph{properties}
such as its $f$-vector, its Hasse diagram, etc. The $\mathcal{V}$- and
$\mathcal{H}$-representation are also regarded as properties, and any
one of them may be used to define the polytope in the first place. If
you are only interested in the combinatorial structure you may also
input the vertex-facet incidences. For more details concerning the
\url{polymake} file format see
Section~A.1.3. Once a polytope is defined in terms of
some property, you may ask
\url{polymake} to compute further properties.

\url{polymake} also provides standard constructions for polytopes. You can
either construct polytopes from scratch (e.g.\ the $d$-cube) or by applying constructions
to an existing polytope (e.g.\ the pyramid). In both cases \url{polymake}
will produce a new file defining the new polytope.

This section explains the command line syntax of \url{polymake} for
constructing, analyzing and visualizing polytopes. The \url{polymake}
file format is briefly described at the end of this section.

\subsection*{A.1.1. Constructions of Polytopes}

If you want to use one of \url{polymake}'s constructions, you have to
call a \emph{client program} to create the polytope. The clients have
appropriate names, e.g.\ the client producing the $d$-cube is called
\url{cube}. All clients are documented at~\cite{GawrilowJoswig}. If in doubt
about the exact syntax of a client program just type the client's name,
press return and you will get a \emph{usage} message.
\medskip

\noindent\textbf{Clients producing polytopes from scratch}\\
    The basic syntax to create a polytope from scratch is the 
    command line
\begin{verbatim}
  <client> <file> [ <options> ]
\end{verbatim}
    For example, to produce a 4-cube type
\begin{verbatim}
  cube c4.poly 4
\end{verbatim}
    Some other clients which produce polytopes you might have heard of
    are
\begin{compactitem}[\ ]\item
    \url{simplex}, \url{cyclic}, 
    \url{cross}, \url{rand_sphere}, \url{associahedron}, and
    \url{permutahedron}.
\end{compactitem}
  \medskip

\noindent\textbf{Clients producing polytopes from others}\\
To construct a polytope from existing one(s), use
\begin{verbatim}
  <client> <out_file> <in_file> [ <options> ]
\end{verbatim}
    For example, to produce the pyramid over our 4-cube \url{c4.poly} type
\begin{verbatim}
  pyramid c4.pyr.poly c4.poly
\end{verbatim}    
    Some other constructions to get interesting polytopes are obtained
    via the clients 
\begin{compactitem}[\ ]\item
    \url{bipyramid}, \url{prism}, \url{minkowski_sum}, 
    \url{vertex_figure}, \url{center}, \url{truncation}, and 
    \url{polarize}.
\end{compactitem}
  
\subsection*{A.1.2. Computing Properties and Visualizing}

To compute a property of a polytope, just ask \url{polymake} using the
following syntax%
\smallskip

\begin{verbatim}
  polymake <file> <PROPERTY_1> <PROPERTY_2> ...
\end{verbatim}
\smallskip

\noindent
The properties are written in capital letters. To compute 
the numbers of facets and vertices of the pyramid over the 4-cube
constructed above, type
\smallskip

\begin{verbatim}
  polymake c4.pyr.poly N_FACETS N_VERTICES
\end{verbatim}
\smallskip

\noindent
Some other useful properties are
\smallskip

\begin{compactitem}[\ ]\item
  \url{GRAPH}, \url{DUAL_GRAPH},
  \url{HASSE_DIAGRAM}, \url{CENTERED}, 
  \url{VERTICES_IN_FACETS}, \url{F_VECTOR}, \url{H_VECTOR}, 
  \url{SIMPLE}, \url{SIMPLICIAL}, \url{CUBICAL},
  \url{SIMPLICIALITY}, and \url{SIMPLICITY}.
\end{compactitem}
\smallskip

\noindent
The visualizations use the same syntax. For example, to take
a look at the graph of a polytope, just use
\smallskip

\begin{verbatim}
  polymake my.poly VISUAL_GRAPH
\end{verbatim}
\smallskip

\noindent
Here are some more visualizations
\smallskip

\begin{compactitem}[\ ]\item
  \url{VISUAL}, \url{SCHLEGEL}, \url{VISUAL_FACE_LATTICE}, and
  \url{VISUAL_DUAL_GRAPH}.
\end{compactitem}

\subsection*{A.1.3. File Format}%

Let's have a brief look at the \url{polymake} file format. If you
have a look at the file \url{c4.pyr.poly}, you will find a paragraph for each
property. The paragraphs are headed by the property's name in capital
letters, followed by the data. The \url{VERTICES}
and \url{POINTS} are represented in \emph{homogeneous} coordinates,
where the first coordinate is used
for homogenization. The inequality
$a_0 + a_1 x_1 + \ldots + a_d x_d \geq 0$ is encoded as
the vector $(a_0 \; a_1 \; \ldots \; a_d)$.

To define a polytope with your own $\mathcal{V}$- or
$\mathcal{H}$-description you should create a file that contains a 
\url{POINTS} or \url{INEQUALITIES} section. 
\medskip

\begin{center}
  \psfrag{D}{$\Delta$}
  \psfrag{x0}{$x_0 \geq 0$}
  \psfrag{x1}{$x_1 \geq 0$}
  \psfrag{x0x1}{$1 - x_0 - x_1 \geq 0$}
  \psfrag{1}{1}
  \epsfig{figure=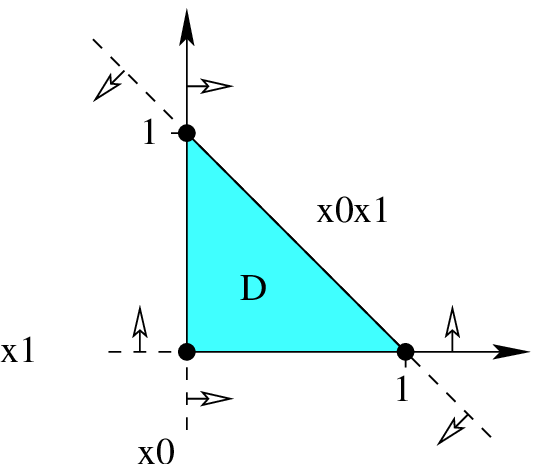,height=6cm}
\\[3mm]
  {\small\textbf{Figure A.1.} The triangle $\Delta$ and its $\mathcal{V}$- and $\mathcal{H}$-description.}
\end{center}
\medskip

\begin{example*}
If you want to construct the triangle $\Delta$ (cf.\ Figure~A.1) that is
given as
\[
\Delta\ \ :=\ \ \conv \{ (0,0),(0,1), (1,0) \}
\]
your input file should look as follows
(homogeneous coordinates):
\begin{verbatim}
  POINTS
  1 0 0
  1 0 1
  1 1 0
\end{verbatim}
If you prefer to enter the same example by the inequality description
\[
\Delta\ \ :=\ \  \{ x \in \mathbb{R}^2 :
1-x_0-x_1 \geq 0, x_0 \geq 0, x_1 \geq 0 \},
\]
this should be your \url{polymake} file:
\begin{verbatim}
  INEQUALITIES
  1 -1 -1
  0 1 0
  0 0 1
\end{verbatim}
\smallskip

\noindent
The difference between \url{POINTS} and \url{VERTICES}
is that the former may contain\break redundancies.
So you should enter \url{POINTS}, but ask for \url{VERTICES}.
Similarly, one should enter \url{INEQUALITIES} and
ask for \url{FACETS}.
\end{example*}

\section*{A.2. The \texttt{polymake} System}

The world seen through the eyes of the \url{polymake} system
consists of \textit{objects} with \textit{properties}. For a given
object you may ask for one of its properties and \url{polymake}
will compute it. Yet \url{polymake} acts only as a framework for
the computation, knowing little about the math involved. The actual
computations are delegated to client programs.
The open design, keeping the system and the math independent,
allows for the use of external software as clients.
It also makes \url{polymake} extremely flexible with respect to the
kind of objects you want to examine.  For example the applications
\url{polytope} and \url{topaz} examine different classes of
geometric objects.
\bigskip

\begin{center}
  \psfrag{clients}{\hspace{-6mm}\url{polymake} clients}
  \psfrag{object}{object}
  \psfrag{object(s)}{object(s)}
  \psfrag{request for a}{request for a}
  \psfrag{property of}{property of}
  \psfrag{an object}{an object}
  \psfrag{rule}{rule}
  \psfrag{base}{base}
  \psfrag{polymake}{\url{polymake}}
  \psfrag{server}{server}
  \psfrag{property}{property}
  \psfrag{others}{others}
  \psfrag{scratch}{scratch}
  \psfrag{producing}{producing}
  \psfrag{from}{from}
  \psfrag{external}{external}
  \psfrag{software}{software}
  \epsfig{figure=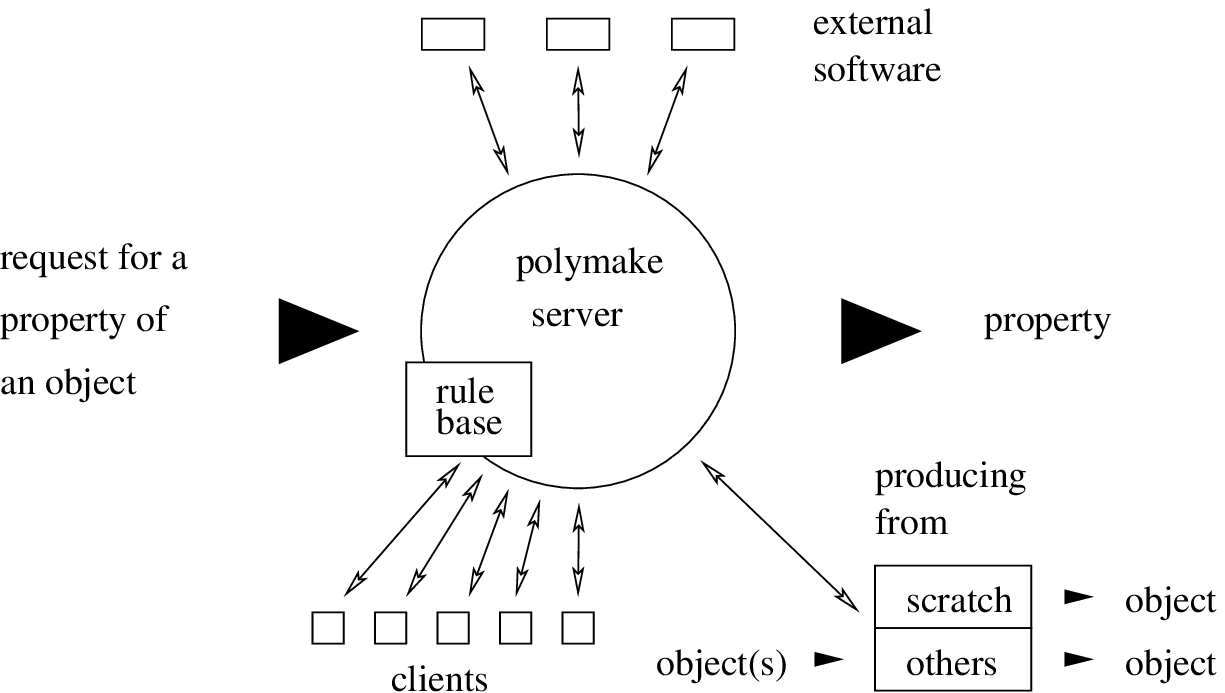,height=7cm} 
\\[3mm]
  {\small\textbf{Figure A.2.} Main components of the \texttt{polymake} system.}
\end{center}
\bigskip

The \url{polymake} system consists of the following three
components as illustrated in Figure~A.2:
\medskip

\begin{compactitem}[~$\bullet$~]\itemsep=1.5mm
\item The \url{polymake} \textit{server} together with the 
    \textit{rule base}.
\item  \url{polymake} \textit{client} programs computing properties
  and new objects.
\item External software, such as the \url{cdd}~\cite{Fuku} convex hull
  algorithm or the \url{JavaView} \cite{polthier-javaview04}
  visualization package.
\end{compactitem}
\medskip

\noindent
The rule base is a collection of rules, each rule containing a set of
input and output properties and an algorithm (that is a client or
external software) which computes the output properties \emph{directly} from
the input properties. If you request a property of an object, the server has to
determine how to compute the requested property from the ones which
are already known. There might not be an algorithm computing the
requested property directly, so other properties might have to be
computed first. Therefore, the server has to compose a sequence of rules
(from the rule base) to be executed in order to compute the requested property.

\subsection*{Examples and Exercises}
\begin{compactenum}[{A}.1.]
\item Construct the $d$-simplex and the $d$-cube for $d = 3,4,5$.
  Visualize the face lattice of the polytope and its Schlegel diagram
  (see Section~\ref{sec:Examples}).
\item Construct the cyclic polytopes $C_3(7)$ and $C_4(8)$. Visualize
  and check Gale's evenness criterion (see
  Example~\ref{example:cyclic}). Alternatively, try constructing
  $C_4(8)$ using the client \url{cyclic_caratheodory}. What is the
  difference?
\item Produce the dual polytopes of the polytopes above. Watch out, they have to
  be \url{CENTERED}. (Check the documentation~\cite{GawrilowJoswig} for the property \url{CENTERED})
\item Construct an octahedron using as many different ways as possible.
\item Build the bipyramid over a square, truncate both apexes and polarize.
  What does the resulting polytope look like?
\item Construct the product of a 5-gon and a 7-gon and visualize it.
\item Construct a 4-polytope with the \url{pcmi} logo as its
  Schlegel diagram.
\medskip

\[
    \includegraphics{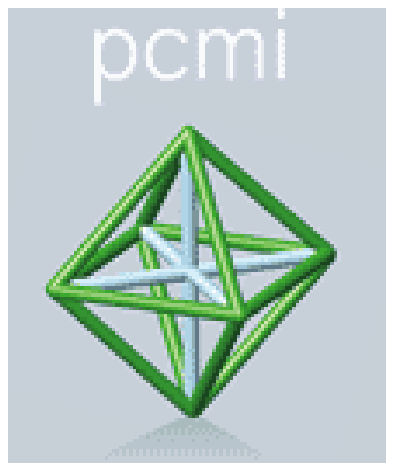}
\]
\smallskip

\item Use the client \url{rand_sphere} to create a random polytope by
  uniformly distributing 1000 points on the unit 2-sphere, visualize
  it and compare it to Figure~\ref{figure:random3dpoly}. Take a look at its
  \url{VERTEX_DEGREES} and cut off a vertex of maximal degree.
\item Take a 3-polytope and truncate all its vertices. Is the
  resulting polytope always simple?
\item
  Truncate the vertices of the $4$-dimensional cross polytope,
  and let \url{polymake} compute the $f$-vector, and 
  whether the resulting polytope is simple.
\newline
  Can you justify the results theoretically?
\item For a given planar 3-connected graph $G$ containing
  a triangular face, produce a 3-polytope with $G$ as its graph (see
  Section~\ref{SteinitzTheorem}). Use the client \url{tutte_lifting}.
\item Visualize the Schlegel diagrams of the \url{dwarfed_cube} 
  that you get by using different projection facets.
\item Visualize the effect of standard constructions (such as truncation,
  stellar subdivision) on the Schlegel diagrams of 4-polytopes.
\end{compactenum}
\label{appendix:POLYMAKE2}


\end{document}